\documentclass{article}
\usepackage{amssymb}
\usepackage{amsmath}
\usepackage{amscd}
\usepackage{oldgerm}
\usepackage{amsfonts}
\usepackage{html}

\def\Re{{\mathrm{Re}}}
\def\Im{{\mathrm{Im}}}
\def\ra{\rightarrow}

\def\la{\leftarrow}
\def\rras{\rightrightarrows}

\def\longra{\longrightarrow}

\def\mpo{\mapsto}

\def\bsh{\backslash}

\def\sbs{\subset}

\def\wh{\widehat}
\def\wt{\widetilde}

\def\a{\alpha}
\def\b{\beta}
\def\d{\delta}
\def\D{\Delta}
\def\e{\epsilon}
\def\g{\gamma}
\def\G{\Gamma}

\def\lb{\lambda}

\def\o{\omega}

\def\s{\sigma}

\def\z{\zeta}

\def\cA{{\cal A}}

\def\cC{{\cal C}}
\def\cD{{\cal D}}
\def\cE{{\cal E}}
\def\cF{{\cal F}}

\def\cI{{\cal I}}

\def\cL{{\cal L}}
\def\cM{{\cal M}}

\def\cO{{\cal O}}

\def\cS{{\cal S}}

\def\cX{{\cal X}}

\def\cZ{{\cal Z}}

\def\tc{\tilde{C}}

\def\ti{\tilde{I}}

\def\tr{\tilde{R}}

\def\tu{{\tilde{U}}}

\def\tx{\tilde{X}}
\def\ty{\tilde{Y}}
\def\tz{\tilde{Z}}

\def\wa{\tilde{a}}
\def\wb{\tilde{b}}

\def\dx{\frac{\partial}{\partial x}}
\def\dz{\frac{\partial}{\partial z}}
\def\dy{\frac{\partial}{\partial y}}

\def\kf{K_{\cal F}}
\def\kxf{K_{X/{\cal F}}}
\def\bkxf{\bar{K}_{X/{\cal F}}}

\def\kx{K_{\cal X}}

\def\ba{\mathbb A}
\def\bc{\mathbb C}

\def\bg{\mathbb G}

\def\bn{\mathbb N}
\def\bp{\mathbb P}
\def\bq{\mathbb Q}
\def\br{\mathbb R}

\def\bz{\mathbb Z}

\def\iy{{\cI}_{Y}}
\def\oxf{\Omega _{X /\cF}}
\def\otxf{\Omega _{\tilde{X} /\cF}}
\def\ox{\Omega _{X}}

\def\cdx{{\cD}_{X}}

\def\pab{\bar{\pa}}
\def\dm{d\mu}
\def\dn{d\nu}

\def\pa{\partial}

\def\iso{\build\ra_{}^{\sim}}

\def\un{{\rm 1\mkern-4mu I}}

\def\build#1_#2^#3{\mathrel{
\mathop{\kern 0pt#1}\limits_{#2}^{#3}}}

\def\gh{{\mathfrak{H}}}

\def\nx{\Vert x \Vert}
\def\mt{|t|}

\def\mx{|x|}
\def\my{|y|}
\def\mw{|w|}
\def\ms{|s|}
\def\mz{|z|}
\def\cd{dx\otimes d\bar{x}}
\def\mmu{|\dm|}
\def\mte{\mt=\e}
\def\mxe{\mx=\e}

\def\mze{\mz=\e}
\def\mtd{\mt=\d}

\def\mzd{\mz=\d}

\def\mzld{\mz\leq\d}
\def\mtle{\mt\leq\e}
\def\mxle{\mx\leq\e}

\def\mzle{\mz\leq\e}

\def\nzs{\Vert z \Vert ^{2}}
\def\nz{\Vert z \Vert }

\def\uny{\un_{Y}}

\def\Aut{\mathrm{ Aut}}
\def\End{\mathrm{ End}}
\def\Der{\mathrm{Der}}
\def\H{\mathrm{H}}
\def\Log{\mathrm{Log}}
\def\Hol{\mathrm{Hol}}

\def\Ext{\mathrm{ Ext}}

\def\Hom{\mathrm{Hom}}

\def\Ker{\mathrm{Ker}}
\def\Coker{\mathrm{Coker}}

\def\mod{\mathrm{mod}}

\def\ghi{\Aut_{-\infty}(\gh)}

\def\RES{\mathrm{RES}}

\def\Spf{\mathrm{Spf}}
\def\Pic{\mathrm{Pic}}

\def\c1{\mathrm{c}}

\def\gmi{{\mathfrak{m}}}

\def \bfp {\mathbf{p}}

\def \bfq {\mathbf{q}}
\def \bft {\mathbf{t}}

\def \tm {{\mathbb{Z}(1)}}

\def \szm {\mathrm{s}_{Z,\dm}}
\def \szmd {\szm(\d)}
\def \szme {\szm(\e)}
\def \lem {\e^{-1}d\e}

\def \imx {\Im(x^{-r})}
\def \imz {\Im(\z)}
\def \imzt {\Im(\z(t))}
\def \rex {\Re(x^{-r})}
\def \rez {\Re(\z)}
\def \rezt {\Re(\z(t))}

\def \imxr {\Im(x^{-r})=R}

\def \re {R^{-1/r}}

\def\ti{{\tau_i}}
\def\tj{{\tau_j}}
\def\tie{{\tau_i(\e)}}

\def\si{{\s_i}}
\def\sj{{\s_j}}
\def\sie{{\si(\e)}}

\def\mch{|\chi|}
\def\sks{{\s_k^*}}
\def\sjs{{\s_j^*}}
\def\skes{{\sks(\e)}}
\def\sjes{{\sjs(\e)}}

\def\tit{{\tilde{\tau}_i}}

\def\tiet{{\tilde{\tau}_i(\e)}}

\def\sit{{\tilde{\s}_i}}

\def\sjst{{\tilde{\s}_j^*}}

\def\sjest{{\sjst(\e)}}

\usepackage{psfrag}

\renewcommand \emptyset {\varnothing}

\title{Two $\Ext$ groups and a residue }
\author{Michael McQuillan}
\date{{\it In memoria di Marco Brunella}}
\begin{document}

\maketitle

\begin{abstract}

The present constitutes the key lemma
which was hitherto missing in the author's
proof of, inter alia, the Green-Griffiths
conjecture for surfaces with enough 2-jets,
{\it e.g.} $13\mathrm{c}_1^2> 9\mathrm{c}_2$.
As to why this is so is not the immediate
concern of this article since it is a lemma
within a whole. By way of what
is essentially an appendix to the current
paper, a guide to the whole is provided at \cite{mysite}.

\end{abstract}


\section*{Introduction}

\subsection*{Prooemium}
A foliation by curves $\cF$ on a 
complex space $X$, or, indeed champ
de Deligne Mumford analytique, is
most usefully viewed as a fibration
$X\ra [X/\cF]$ over the champs classifiant. 
Supposing $X$ smooth for ease of exposition,
then away from the
singularities of $\cF$ this has
its obvious sense, cf. \S I.1,
otherwise, while one can make
definitions at the singularities,
the sense of this arrow is ambiguous,
and it's role is largely for the
memonic purpose of emphasising
the dynamical nature of the study.
What is unambiguous, however is
that there is an exact sequence
of coherent sheaves on $X$,
$$\begin{CD}
0@>>>\oxf@>>>\ox@>>>\kf\iy@>>>0
\end{CD}$$
with $\kf$ a line bundle, which
one may think of as the canonical
bundle along the leaves so strictly:
$K_{X/[X/\cF]}$, whence the abbreviation,
and $\iy$ an ideal supported in
co-dimension 2, {\it i.e} the
singularities. This leads to many
invariants, known as Baum-Bott
residues, \cite{bb}, which express
this ambiguity to first order,
{\it e.g.} an invariant measure,
$\dm$ is an unambiguous notion
away from the singularities which
one usually supposes (I.1.1) to
extend across the singularities
as a closed positive current, which
is not supposing very much, whence
there is a residue,
$$
\begin{CD}
\dm\in\Ext^{n-1}(\iy,\kxf)@>>{\RES}>
\Ext^n(\cO_Y, \kxf)
\end{CD}
$$
where $\kxf=\kx-\kf$ may be thought of
as the canonical
bundle of $[X/\cF]$, albeit even this 
statement is
ambiguous at the singularities, and
is subject to other residues, cf. I.1.5. 
Some limited insight may
be had when the foliation is, say,
a fibration $X\ra [X/\cF]$ of a surface
over an algebraic curve, so everything
extends un-ambiguously over the
singularities, and we find singular
fibres,
$$X_b=\sum_i n_i C_i$$
with multiplicities $n_i$ along the
reduced components $C_i$. According to
the definition of an invariant measure,
every combination,
$$\sum_i \nu_i C_i, \,\,\,\,\, \nu_i \in \br_{\geq 0}$$
with at least one $\nu_i\neq 0$ is to
be considered such. However $\RES$
vanishes on such a  combination iff
it is parallel to $X_b$ iff the
measure descends to the quotient.
Similarly straightforward remarks
apply in the ``universal case" $\cM_{g,1}\ra\cM_g$.
Already, however, for a slight perturbation,
cf. III.4.1 (c), of such algebraic
examples, even locally $[X/\cF]$
is quite complicated, {\it i.e.}
the action of a diffeomorphism which
is no better than conjugate to a
rotation up to some finite order,
and $\RES$ sees only a very small
part of the dynamics. Nevertheless,
the part which it does see is exactly
what is pertinent to algebraic geometry.

To understand this relation, we require
a brief parenthesis, \cite{uu} \S III, on closed positive
currents $T$ of dimension $(p,p)$
(so acting on $2p$ forms)
defined in a neighbourhood
$U$ of some proper sub-variety $V$.
An example of which would be integration
over a holomorphic $p$-cycle, and in
this case one has the Fulton/Macpherson specialisation,
\cite{fulton},
of the cycle to a $p$-cycle on the normal 
cone $C_{V/U}$, and indeed without any
hypothesis of the properness of $V$.
In general, however, properness appears
to be a necessary assumption to define
a specialisation of $T$ to $C_{V/U}$,
which again is closed positive (non-negative
would be more accurate) of dimension $(p,p)$.
This in turn defines, albeit only
unambiguosly on closed forms, a closed positive
$(p-1,p-1)$ current $\mathrm{s}_{V,T}$ on the
projectivised cone which yields exactly the
Segre class of op. cit. when $T$
is an analytic $p$-cycle. This is perhaps
clearer in the immediate case of relevance: 
$V$ is a divisor and $T$ of type $(1,1)$,
then $T$ splits into closed positive
currents,
$$\un_{U\bsh V} T\, +\, \un_V T$$
for $\un_{\bullet}$ the characteristic
function. The Segre class of $T$ around
$V$ is the same as that of $\un_{U\bsh V} T$,
it is a measure on $V$, and its total
mass is given by,
$$\int_V \mathrm{s}_{V,T} \, =\, V\cdot\un_{U\bsh V} T$$
with the latter product just being that
between $\H^2_c$ and $\H^{2(n-1)}$. Plainly
for $T$ of dimension $(1,1)$ the general
case can always be reduced to this case
by blowing up. In particular the Segre
class $\mathrm{s}_{V,T}$ is, functorially with
respect to the ideas, the winding number
of $\un_{U\bsh V} T$ around $V$. 

In contrast to $\RES$ the Segre class
of an invariant measure, or indeed
of any closed positive current, around
a sub-variety, which here will be the
singular locus $Y$ of the foliation
is highly amenable to estimation by
algebraic techniques such
as counting sections of linear systems
vanishing to high order on $Y$.
Whence the following relating
the two is central to the global study
of foliations,

\noindent{\bf Lemma} {\it Let $U\ra [U/\cF]$
be a foliated 3-fold 
with (foliated) canonical singularities
in which the singular
locus $Y$ is proper (so, in practice a tubular
neighbourhood of the same) and $\dm$
an invariant measure then,}
$$\mathrm{s}_{Y,\dm}=0\, \iff\, \RES(\dm)=0$$

Here, by way of clarification, let us
note that by hypothesis/definition
$\un_Y\dm$ is 0 otherwise neither side has
sense. The notion of canonical
singularities is functorial with
respect to the ideas, and the
existence of such resolutions has
been proved in \cite{mp2}.
Plainly the content of the
lemma is to reduce the study of a
dynamical invariant, {\it i.e.} something
about $[X/\cF]$ to a much simpler algebraic
one. Before overviewing its proof, let
us note an illustrative corollary,





\noindent{\bf Corollary}\cite{uu} {\it Let $S$ be an algebraic
surface of general type with enough 2-jets, e.g. $\mathrm{c}_1^2(S)>9/13\mathrm{c}_2(S)$,
then the Green-Griffiths conjecture holds.} 

In fact, way more is true for such surfaces, {\it e.g.} an optimal version
of Gromov's isoperimetric inequality, boundedness of curves
of genus $g$ in moduli, a Bloch type dichotomy for convergence
of discs in the Gromov-Hausdorff sense, {\it etc., etc.}. One
can also apply it to Mordell for surfaces over characteristic
zero function fields too, but apart from the 2-jet condition on
the generic point, one needs related conditions on the fibres
of bad reduction. This is, however, purely illustrative, it is
really a lemma about foliated 3-folds of which it almost implies
a complete description in the spirit of \cite{canmod}. Indeed,
what remains to be done to achieve this is the construction of
$\kf$-flops, which in one of life's ironies is more difficult
than flips which were done in all dimensions in \cite{ss}, and
to improve the lemma so that it implies something when the
numerical kodaira dimension, $\nu(\cF)$, of $\kf$ is 1. More precisely 
not only for $\nu(\cF)\neq 1$ are the hypothesis of the lemma adequate,
but, quite generally the lemma addresses the most difficult
part of the computation of $\RES$, to apply it however to
the case of $\nu(\cF)=1$ one needs a more delicate hypothesis
than vanishing Segre class close to the beast of \S III.5. This
is a particular type of foliation singularity with very wild
behaviour even formally which the vanishing Segre class hypothesis
renders benign, but which will certainly show its teeth
as soon as one seeks a more quantitative
statement. As a result, even with all of the machinery
at our disposal we still
cannot quite claim non-trivial analogues of the above
corollary for 3-folds.

All of which, apart from the side stepping of a subtle
dynamical issue at the beast, is irrelevant except as
motivation since the proof of the corollaries is
elsewhere, {\it i.e} \cite{uu}. Let us, therefore,
turn to the question of the proof. The lemma is
stable under blowing up in centres inside
the singular locus, so we can, and will, blow
up as much as we can to improve the situation. 
Any sequence of such operations has the following
effect: we get a modification $\rho:\tx\ra X$,
with a total exceptional divisor $E$ invariant by
the induced foliation, again denoted $\cF$ with
no change in the canonical, so an exact sequence,
$$
\begin{CD}
0@>>>\otxf (\log E) @>>> \ox(\log E)
@>>> \kf I_{L}@>>> 0
\end{CD}
$$
The resulting ideal $I_L$ is supported in
what \cite{uu} describes, correctly, as
points where the foliation is not log-flat.
Log-flatness is a stable condition under
blowing up, and the lemma is trivial when
it holds.
This leads to a finite list,
\cite{uu} I.3, \cite{mp1} V.1.8, VI.2.1
 of final forms
for the singularities where the lemma has
to be proved. There is only one isolated
case which can be done in 
a couple of lines in
the same way that the surface
version of the lemma was proved 
in \cite{ihes}, which we'll here by
refer to as the ``baby lemma''. This is, however,
relevant and the proof is recalled in \S I.3.
Otherwise there is a formal (2 dimensional)
centre manifold, $Z$, in the completion of the
singular locus, which itself has dimension 1.
If this were to converge, then there is a
rather good trick in \cite{uu+} for reducing
the lemma to a measure supported in the
centre manifold. Of course, even for saddles
on surfaces the centre manifold does not
converge, but there are large sectorial
domains where it does, and smaller, albeit
only relatively, sectors, where the holonomy
on transversals has basins of attraction.
The relevant meaning here of large is
large enough to define the centre manifold
uniquely. As such the motivation for \cite{mp1},
was to prove similar behaviour in dimension
3, and deduce the lemma from the trick of
\cite{uu+}. Although the results of \cite{mp1}
are optimal, this particular assertion is
hopelessly false, and the basic conclusion
of \cite{mp1} is that even locally the
general singularity that we must worry about
is insanely complicated in a way that makes
the theory of 2-dimensional saddles look like
a game for age 5 and under. An exception to
this is when locally many leaves adhere in
the exceptional divisor, at which point one
can make statements very similar to those
found in dimension 2, and, basically, \S II, 
invariant
measures don't exist, which does make the
computation rather easy. The generic situation,
however, is that locally there are no
obstructions to existence, and, worse,
where the leaves are asymptotically as
if they factored through the centre
manifold no local obstruction to the
measure being absolutely arbitrary. This,
and the fact that the injectivity of
distributions implies that locally $\RES$  
is meaningless anyway is explained in
more detail in \S I.2. We therefore have
the following problems,

\noindent {\bf (a) } A computation that only has meaning
globally.

\noindent{\bf (b)}  A dynamical structure, $[X/\cF]$ which
is so complicated that even locally it's rare that
anything
pertinent can be said.

\noindent{\bf (c)}  A centre manifold that probably 
doesn't converge, but it might, and, regardless,
two completely unrelated kinds
of behaviour for leaves ``outside'' it-
generally escaping from our neighbourhood,
and ``inside'' it- going nowhere.

The final point is perhaps the most
important to understand, and it contains
the following ``toy lemma'': suppose the
centre manifold did converge, and the
measure was supported in it, how do we
prove the lemma ? Well, the singular
locus has two kinds of components,
those left invariant by the induced
foliation in the centre manifold,
and those transverse to it. The latter
are easy, because on blowing up they
can be supposed everywhere transverse,
and the zero Segre class hypothesis
implies the measure doesn't exist, albeit 
conversely,
II.1.6, as Kontsevich pointed out to
me, the lemma at such components is
best possible, {\it i.e.} it has no
quantitative variant. Otherwise the toy
lemma does not follow from the baby lemma,
since 
the residue here has poles that are
much worse.
It almost follows from a
much better lemma, \cite{marco}, of Marco Brunella
in dimension 2, to the effect that
the completion of the holonomy around
each invariant curve must take values
in $S^1$. This, together with standard
existence results about the convergence
of branches through foliation singularities
on surfaces is, in fact, enough if the
dual intersection graph formed by these
invariant components has no cycles. Otherwise,
one can get a super attracting phenomenon
around such cycles which a priori buggers
up the calculation without excluding
invariant measures, I.4.2. This is due
to what might be called logarithmic
holonomy, which is introduced in \S I.4,
but on combining this with Brunella's
theorem, one would find that were there
to be super attraction around cycles
in the graph then the measure would be radially
symmetric on transversals, and this
permits the computation to be done.

All of which looks like we're moving
forward, but remember this was a toy
problem, and the only
step in the above reasoning that
is a priori meaningful
in our situation is the first one
about the transverse components, \S II.1. 
At the same
time, the above might be our situation,
so we must be able to do the toy 
problem with what we have at our
disposal. One thing we have at our
disposal is the formal holonomy in the formal
centre manifold,
and, in the context of the toy
problem, one realises that 
this can be used to re-prove Brunella's theorem.
Similarly, one may prove approximate radial
symmetry, {\it i.e.} for transversal
discs of radius $\e$, radial symmetry
up to an error $O(\e^n)$,
in
the presence of (formally) super attracting
cycles, and
being carefull not to use  facts
about convergence of branches through 
foliated surface singularities that
are false in dimension 3 yields a
solution to the toy problem which doesn't
use anything we don't have except the
toy hypothesise themselves.

This definitely constitutes progress,
and merits a re-examination of the
baby lemma. The only case of the
baby lemma that is non-trivial is
the saddle node, and we now have 3
approaches,

\noindent{\bf (1)} The original. A quick 
2 line trick using nothing but Stokes' theorem,
and the first few terms in the formal
normal form of the singularity,
but, unfortunately, utterly reliant on 
the singularity
being isolated, or, more accurately
the complete absence of curvature.
In particular it does not imply the
toy lemma.

\noindent{\bf (2)} A 
theorem by Brunella that the measure
itself must be integration over a
convergent curve through the singularity,
of which there are at most 2. 
In particular, when there is zero Segre
class there is no measure outside the
exceptional divisor. This renders the
baby lemma trivial, but requires
an entire
series of theorems about conjugation
to normal form on sufficiently large sectorial
domains, which are not only much more
difficult in dimension 3, but generically,
cannot  meet the sufficiently large condition.
A variant proves the toy lemma.

\noindent{\bf (3)} A via media. Just take anything
which agrees with the normal form to sufficiently
high order. The measure is approximately
invariant by its holonomy, and this affords
a proof of Brunella's theorem using only
the formal structure of the singularity. 
This generalise to a proof of the
toy lemma.

To which it should be added that whether
it be the invariance of the measure by the
holonomy on transversals employed in (2),
or the approximate invariance in (3), either is
most usefully seen by way of Stokes, and
so, post factum, (1) is really just a
special case of (3) in which just enough
and no more of the Taylor expansion was
considered.

At this point, however, all we've done is
to understand the obvious, {\it viz:} without
any need for analytic conjugation to
normal form, the latter alone will give
whatever bound we want on transversals
whenever a naive inspection of the
normal form suggests there should be 
transversals with no measure. This might
reasonably be described as {\it approximate
holonomy}. However. under the conditions where
it would be applicable, 
not by coincidence,
one can, \cite{mp1},
obtain actual analytic conjugation to 
normal form on large sectors, and actually
conclude from the holonomy that the measure is zero in a
way that is more convenient than approximate
holonomy: an infinite sum of
zeroes is zero, but an infinite sum of
small still needs to be estimated.
Whence, there is still a
difficulty to overcome, and this is most easily grasped by
considering the particular, and improbable,
sub-case, \S IV.1, where the singular locus has
exactly one component and the induced
foliation in the centre manifold is
everywhere smooth. Let us call this the
``warm up lemma''. In such circumstances, on
completion in the singular locus (not
just at some point, which would
be useless) one 
finds on a small open formal  
coordinates $x,y,z$, such that the
foliation is given by,
$$z\dz + x^p\dy$$
for some $p\in\bn$. The exceptional
divisor is $x=0$, and the centre 
manifold is $z=0$. The curves $x=0$,
$y=$const. $z\neq 0$ converge, but the holonomy
around them isn't very interesting,
whereas the good candidate for a
curve which governs the dynamics
is the singular locus itself, which,
as it happens, would be an elliptic curve.
On the other hand, holonomy around a
curve in the singular locus is
meaningless. Certainly there is the
formal holonomy in the formal centre
manifold, but a priori this has
no relation whatsoever with the
structure of $[X/\cF]$, {\it i.e.}
the dynamics. A
local analogy, is the
pleasing example II.2.1,
where the centre manifold exists
on large sectors, associated to which, 
cf.
the proof of II.2.3, there is
an \'etale groupoid $R\rras T$ 
isomorphic to an action of $\tm$
on $\bc$ whose completion
is the formal holonomy, nevertheless
this association is simply by way
of an abstract gluing operation,
and it is not the case that $[T/R]$
maps to $[X/\cF]$. In the presence
of an invariant measure, however,
II.2.3, this abstract structure
becomes relevant, and it is how
measure in the centre manifold is
eliminated. For what pertains to
the warm up lemma, there is no
comparable analytic conjugation
to normal form. This can only be
done under the condition that the
argument of $x$ does not turn
through more than $\pi/p-\d$. Were
it possible to do $\pi/p +\d$, the
centre manifold would be unambiguously
defined in such sectors, and on
the transversals $y=$const. one
would find the mass bound 
$\e^p$
for
balls of radius $\e$ off the
centre manifold, whence, off the
centre manifold, it would be
possible to justify calculating $\RES$
locally, and otherwise one would
aim to conclude the warm up lemma
by a global analogue of II.2.3
to bring the formal holonomy in the
formal centre manifold into the game.
We know, however, by \cite{mp1}
that this cannot be justified,
but like the local example II.2.3,
the measure itself brings the formal
holonomy in the centre manifold into
play and leads to the {\it almost holonomy}
estimate, III.1.2. Unlike a conjugation
to normal form on large sectors, it
need not (and does not in the 
above example) provide local information on the
mass of transversals. It does, however,
provide mass bounds on transversals
to the singular locus as a function of
the formal holonomy in the centre
manifold. In the particular case
of the warm up lemma, \S IV.1,  it's quite
easy to see that these bounds are
exactly what one needs, and if one
is prepared to assume the useful
fact the the singular locus is an
elliptic curve in the warm up lemma,
or the corollary that the formal
holonomy in the formal centre manifold
is commutative, it's particularly
easy to see how it solves all of
the above problems (a)-(c).  
Thus, there is a good case for
working through the warm up lemma
before attempting anything else.

Once one has the almost holonomy
estimate, it's reasonably clear
that everything is going to work.
Nevertheless the calculation is
long, and is organised as follows:

\noindent($\a$) Holonomy/Approximate holonomy
at certain singularities. More
precisely, 
apart from the warm up lemma,
there will be singularities in the
induced foliation in the formal
centre manifold, which, in turn
lead to complicated singularities
in the ambient 3-fold. In particular,
there are some where the leaves can
be locally very large and adhere in the singular
locus. Until we can exclude that there
is no measure here this is an obstruction
to making use of the almost holonomy.
The obstruction can be dealt
with by looking either at the holonomy
or the approximate holonomy on transversals
to curves through the singularities.
As we've said \cite{mp1} applies here,
and this has been done, \S II \& III.3.2, 
by holonomy, but only to avoid repetition
of calculations that one finds in
op. cit.. In the case, III.3.1.(a),
however of a generic node in the
centre manifold the vanishing of
the measure in the region defined
in III.3.2 is only a paranthesis,
and is not actually used.
This is because such nodes have
many formal invariants, and certain
combinations of these yield a region
of nil measure which appears to be
a bit too small to be useful. As such,
we use a variant of the original proof
of the ``baby lemma'' to get appropriate local
estimates which more properly should
be ascribed to approximate holonomy.

\noindent($\b$) Almost Holonomy. Off the singularities,
component by component this works exactly like
the warm up lemma, and is the basic technique
for addressing the non-local nature of the
calculation. It also provides the local
strategies at singularities, \S III, where
holonomy/approximate holonomy does not apply,
or, II.3.3 \& III.3.2, provides only partial
information, and, how to glue the local strategies
to the global one. Nevertheless, as we've
attempted to explain, almost holonomy is
not approximate holonomy, but rather an
estimate that manages to address simultaneously
the distinct dynamical features of problem (c).
In particular this estimate has little room
for error, and has nowhere near the robustness
of approximate holonomy. One might think post
\S III.2 that it'll be plain sailing nevertheless,
but when the almost holonomy estimate is provided by holonomy tangent 
to a rational rotation there is need for caution.

\noindent($\g$) The topology of the intersection graph. The
plan, cf. \S I.4,  was to do this by  ``logarithmic almost
holonomy'', but I couldn't get the
right shape of estimate. Instead a more
algebro-geometric strategy has been adopted,
which leads to a wholly different difficulty arising
from ludicrously improbable combinations of 
singularities in the centre manifold with 
a rational eigenvalue.

On a philosophical level, 
the combination of approximate
and almost holonomy addresses
what the author has always
considered to be the most
difficult part of his programme
to prove everything you want
to know about algebraic surfaces,
but were afraid to ask
{\it i.e} the computation of $\RES$,
and one might even say, the ``right''
proof of the baby lemma. 
It only uses
the natural tools at one's disposal:
the measure regularity (the lemma
is rubbish for 
invariant distributions) and the
formal structure of the singularities,
which, given that an invariant  measure 
(rather than the restriction
of a 2 extension)
on 
completion 
in the singularities
is meaningless is reason
to be cheerful.
Whence, it's a bit disappointing that
the calculation is so long. The
nature of ($\a$) is clear, and it
should be possible to produce quite
a clean theory based on approximate
holonomy in all dimensions. The
problem posed by ($\g$) is real,
and it would have been preferable
to have sub-ordinated this to ($\b$).
Otherwise, the use of the almost
holonomy is rather uniform, and
could arguably have been presented in a less
ad-hoc way. Some
singularities, however, have a much more
complicated local structure than others, 
so much so that they have to be looked
at closely, and only yield to a 
combination
of approximate and almost holonomy,
so, 
there is an inductive structure here
that remains to be organised. 

By way of thanks, I am indebted to
Adam Epstein for an illuminating
tutorial on the dynamics of 
diffeomorphisms tangent to the
identity, to Kontsevich for avoiding
a wild goose chase after a more
quantitative lemma, and Sibony for
a precision on the (non-) uniqueness
of the Segre class as a current, albeit
in this article it's highly unique, {\it i.e.} 0.
It should, however,
be plain from the introduction that
this article owes much to Marco Brunella.
In our last mathematical conversation,
I posed to him the problem which is
the intersection of the toy lemma with
the warm up lemma. I was already fairly
sure that even such a simple case had
to be done by way of his structure
theorem for invariant measures on
surfaces, but I still hoped that I
was overlooking some simple point, which
was necessarily much the same point as
to why \cite{mp1}, which had just been
completed, did not imply the lemma.
Consequently, I gave him no hints,
and after 15 minutes or so, he replied,
well the holonomy is $S^1$, ``ma questo
non ti piace''. I had learned holonomy
from him, and his memory was of a
suspicious amateur. Such days, however,
were in the past, and I replied that
I agreed that the $S^1$ holonomy argument
was the right one, but the real problem
was way more difficult, and I had no
idea how a $S^1$ holonomy argument could
be  organised, or even what it might mean. My worst fears thus confirmed,
yet loathing to accept that I'd spent years
thinking in totally the wrong terms, I 
pleaded for an
answer different to what I had just been
given, and asked ``Sei sicuro che non c'\`e
altro modo ?''. He replied, ``Sono sicuro,
sono sicuro''. Requiescat in pace. 

\newpage

\subsection*{De radice -1}

Let $\bc$ be the complex plane,
with $o_\bc$ its orientation sheaf,
and $M$, 
$\cA^{2}$, and $\cA^{1,1}$
the sheaves of measures,
2-forms, and
tensors (not differential
forms) of type $(1,1)$, all
with real values, then there
are canonical inclusions,
$$i:\cA^{1.1}\hookrightarrow M
\hookleftarrow \cA^{2}\otimes o_\bc :j_\bc$$
the subscript on $j$ being used to
indicate that it is in fact, a
natural transformation, to which
there is, furthermore,
a unique canonical isomorphism,
such that the following
diagram commutes,
$$
\begin{CD}
\cA^{1,1}@>{\sim}>>  \cA^{2} \otimes o_\bc\\
@VViV @VV{j_\bc}V\\
M@= M
\end{CD}
$$
Consequently for
$\br(1)$ the imaginary numbers there
is an isomorphism,
$$a: \cA^2\otimes o_\bc \iso \cA^2(1):
dz\otimes d\bar{z}\mpo -2\pi dzd\bar{z}$$ 
which does not depend on the choice of
the square root of $-1$, 
albeit that this requirement, plainly,
does not determine
$a$ uniquely.
Irrespectively, the choice of $a$, so,
again independently of a choice of
the square root of -1,
defines a map,
$$\int:\G_c (\bc, \cA^2(1))\longra\br$$
Similarly for $\pa U$ the border of a
bounded domain $U\sbs \bc$,
there is the Stokes isomorphism, $\s_U:
o_\bc\ra o_{\pa U}$, which although conventional,
has the pleasing feature of relating two
canonically defined quantities in a convenient way, {\it viz}:
$$\int_U j_\bc (d(\o)) \,=\, \int_{\pa U} j_{\pa U}(\s_U \o),
\,\,\,\,\, \o\in \G(\bar{U}, \cA\otimes o_\bc)$$
It therefore follows, still independently
of any choice of the root of -1, that
there is a map,
$$\oint_{\pa U}: \G( \pa U, \cA^1(1) ) \longra \br$$
depending on the
Stokes isomorphism, and the isomorphism $a$.
The agreeable nature of the former having
been observed, the latter is chosen so
that for $U$ the unit disc $\D$,
$$\oint_{{\pa {\D}}}\dfrac{dz}{z}\,=\, 1$$
which together with functoriality under
conformal mapping serves to define $a$
uniquely. In particular, by way of notation,
this is how the contour integral symbol 
will be understood. In addition, since
$a$ is a choice, albeit an elegant one,
when integrating an actual measure we
shall eschew it, and work with tensors
rather than differential forms.

Related to this we have the exponential
mapping,
$$\exp: \bc\longra \bc^{\times}$$
which affords a canonical and base
point free isomorphism,
$$\pi_1(\bc^{\times})\longra\tm\,\,\, (=\bz 2\pi \sqrt -1)$$
using which we will freely identify oriented
loops $\g$ around a punctured with their image
in $\tm$ which will equally be denoted $\g$.
Being oriented, $\g$ defines a distribution
via integration, which in the case that $\g$
borders the unit disc $\D$ is related to
the previous discussion by way of,
$$\oint_{{\pa {\D}}} \,=\, \dfrac{1}{\g} \int_\g$$ 
and, again, no choice of a root of $-1$
is necessary.

Unfortunately there will be choices that may
be ``forced''. One way this can arise is when
we have a loop $\g$, identified to its image
in $\tm$ and a number $\lb\in\br(1)$, typically
an eigenvalue, so the ratio $\lb/\g\in\br$,
and its notationally more convenient to express formulae
when the sign of this number is, say, positive
rather than negative. In such circumstances
we will typically make the choice of a square
root of $-1$ in the form of an imaginary part
function $\Im$. Another example which is
worse still, and it will occur, is that we find
ourselves in a domain which contains one root
rather than the other, and again ease of
notation forces a choice of $\Im$, and if both
of these occur at once, then we should be
carefull.

\newpage

\section*{I. Residues}

\subsection*{I.1 Measures and co-homology}

Let $X\ra [X/\cF]$ be a foliation by curves of a smooth
complex space, or, indeed smooth champ de Deligne-Mumford
analytique, $X$ of dimension $n$, then by definition we
have a short exact sequence,
$$0\longra\oxf ^{1}\longra\ox ^{1}\longra\kf\iy\longra 0$$
where $\kf$ is a line bundle, $\oxf$ is reflexive, and the
support of the singular locus $Y$ has co-dimension at least
2. Away from the singularities an invariant measure $\dm$
may be unambiguously identified with a closed positive
$(n-1,n-1)$ current lying in $\kxf\otimes\bkxf$ for
$\kxf$ the double dual of the top power of $\oxf$. More
usefully, there is a smooth groupoid,
$$
F  \, \textstyle\build\genfrac{}{}{0pt}{}{\longra}{\longra}_{s}^{t} \, X\bsh Y
$$
which may be sliced along a not necessarily connected
transversal $T$ to yield an \'etale groupoid,
$$
R \, \textstyle\build\genfrac{}{}{0pt}{}{\longra}{\longra}_{s}^{t} \, T
$$
while in addition we have a fibre square,
$$
\begin{CD}
X\bsh Y @<<{\rho}< L=\build\coprod_{t}^{} L_{t} \\
@V{\pi}VV @VV{\pi}V  \\
[X\bsh Y/F]  @<<{\rho}< T
\end{CD}
$$
where $\pi$ is the projection to the champ classifiant
of the foliation, $\rho$ is \'etale, and $L_t$ is the
leaf through $t$. Whence,
\begin{itemize}
\item $\rho^*\dm$ is a closed $(n-1,n-1)$ current lying in
$\pi^*K_T\otimes\pi^*\bar{K}_T\iso\kxf\otimes\bkxf$, and
thus descends to a measure $\dm(t)$ on $T$, {\it i.e.} 
closure is not just equivalent to a local descent datum,
but also a global one since $\rho^*\dm$ is itself global,
and by definition $\rho^*\dm=\pi^*\dm(t)$.
\item By stokes there is even a descent datum for $R \rras T$.
Indeed if $f\in R$ is an arrow with source $\tau=t(f)$, and
sink $\sigma=s(f)$ then we may find a small neighbourhood
$V\ni f$ such that $s,t$ are homeomorphisms about $f$, then
connect points of $t(V)$ to $s(V)$ by way of real paths
inside leaves to form a real co-dimension $1$ manifold
$\tilde{V}$ in $L$ so that for any function $g$ on $t(V)$,
$$0=\int_{\tilde{V}} d(\pi^*g) \dm = \int_{s(V)} s_*t^*g\dm
-\int_{t(V)} g\dm = \int_T (s_* h - t_* h) \dm$$
for $t^*g$ identified with the function $h$ on $V\sbs R$.
The final integral will be referred to as a {\it co-equaliser}.
It is well defined whenever both $s_* h$ and $t_* h$ are
absolutely integrable, {\it e.g.} bounded Borel functions.
\item Equivalently, but less informatively, $R \rras T$
is \'etale, so $\kxf$, indeed any power of $\oxf$,
is a bundle on $[X\bsh Y/F]=[R/T]$. The class $\dm$ is
global and locally closed, so it descends to,
$$\dm \in \Gamma ([X\bsh Y/F], M\otimes \kxf\otimes\bkxf)$$
where $M$ indicates that the coefficients have measure
regularity on the given (whence any) \'etale transversal
$T\ra [X\bsh Y/F]$.
\end{itemize}

Around the singularities there is no similar discussion.
Indeed the groupoid $R \rras T$ may only very rarely
be completed across the singularities (basically probability
zero with respect to Lebesgue measure for any moduli space
of foliations) in a way that $s,t$ remain discrete. There
are ways to do it by the addition of $B_{\bg_m}$'s, but
this results in points with negative dimension, which
is demonstrably meaningless from the point of view of
measure theory, whence we'll take the usual ad-hoc approach,
{\it viz}:

\noindent{\bf I.1.1 Definition} By an invariant measure
$\dm$ for the foliation is to be understood a closed
positive $(n-1,n-1)$ current with $\uny\dm= 0$ (where
here, and elsewhere, $\un_*$ will be the characteristic
function of a set), which lies in 
$\kxf\otimes\bkxf$ away from the singularities. In particular
the only hypothesis at the singularities is that it
extends across the same with finite mass.

There is a co-homological shadow of this ambiguity
by way of,
$$A_c^{0,1}\otimes\kf\iy\longra A_c^{0,0}\otimes\kf\iy\otimes\bar{K}_{\cF}\bar{I}_Y
\build\longra_{\dm}^{}\bc$$
whence $\dm$ defines a class in $\Hom (A_c^{0,1}\otimes\kf\iy, \bc)$,
which by way of the Verdier isomorphism is a closed class
in $\Hom_X( \iy, \cdx ^{0,n-1}\otimes\kxf)$, for $A, D$,
respectively $\cA, \cD$,
smooth functions, and distributions, respectively their
sheafifcation. The complex,
$$\cdx ^{0,*}:\cD _{\cX}^{0,0}\build\longra_{\pab}^{}
\cdx ^{0,1}\build\longra_{\pab}^{}
\hdots\hdots\build\longra_{\pab}^{}\cdx ^{0,n}$$
is an injective resolution of the structure sheaf, so
we obtain:
$$\dm\in\Ext_X^{n-1} (\iy,\kxf)$$
As such the obstruction to $\dm$ being a co-homology class
is purely local, {\it i.e.}
$$\H^{n-1} (X,\kxf) \longra \Ext_X^{n-1} (\iy,\kxf)  
\build\longra_{\RES}^{}\Ext_X^{n} (\cO_Y,\kxf)$$
For $Y$ compact the final group is in duality 
with $\H^0(\kf\otimes\cO_Y)$, and $\RES$ may be
made explicit as follows: take $U$ any neighbourhood
of $Y$- in practice small and tubular- with $\o_{\a}$
local liftings to $\kf$ of some $\tau$ in $\H^0(\kf\otimes\cO_Y)$,
and $\rho_{\a}$ compactly supported bump functions forming a partition
of unity on some smaller neighbourhood $V\sbs\sbs U$
of $Y$, then,
$$\tilde{\tau} =\sum_{\a} \o_{\a}\pab\rho_{\a}\in A_c^{0,1}\otimes\kf(U)$$
may be identified with an integral $(1,1$ form, and:
$$\RES(\dm) (\tau) =\int_U \tilde{\tau}\dm = \lim_{\e\ra 0}\int_{\pa Y_{\e}}\o\dm,
\,\,\,\,\,
\o= \sum_{\a}\o_{\a}$$
where $\pa Y_{\e}$ is any continuous boundary around $Y$
at some distance $\e$, tending to zero in any sequence whatsoever
since $\tilde{\tau}$ is integrable. All of which is
equally true for any $\Ext$ class with measure regularity.
Otherwise, identifying the connecting homomorphisms, 
$$\Ext_X^{n-1} (\iy,\kxf)  \build\ra_{\RES}^{\d}\Ext_X^{n} (\cO_Y,\kxf),
\,\,\, \H^0(\kf\otimes\cO_Y) \build\ra_{}^{\d^\vee} \H_c^1 (U,\kf\iy)$$
with $\pab$, {\it i.e.} $\d^\vee =\pab\o$, as above,
and $\d\phi=\pab\theta$ for $\phi$ an $\Ext$ class,
$\theta\in D^{0,n-1}\otimes\kxf (U)$, one has,
$$\phi(\pab\o)=\theta(\pab\o)=(\pab\theta)(\o)=(\pab\theta)(\tau)$$
which is rather less than being able to write the 
obstruction as an unambiguous limit of an honest
integral, even if it remains a ``residue'', {\it i.e.}
concentrated around the singularities. Consequently
although our interest is confined to invariant
measures we do quite generally have a pairing,

\noindent{\bf I.1.2 Definition} For $U$ any neighbourhood of a
supposed compact connected singular locus $Y$,
$$\RES_Y:\Ext_X^{n-1} (\iy,\kxf)\times \H^0(\kf\otimes\cO_Y):
(\phi,\tau)\mpo \RES_Y (\phi,\tau) = \RES(\phi) (\tau)$$
and this pairing is independent of the neighbourhood
$U$ of $Y$.

The principle utility of this pairing is for $X$
compact, albeit its calculation is local around $Y$.
More precisely, we have a diagram with exact rows:
$$
\begin{CD}
\H^0 (X, \kf) @.\ra \Ext_X^1(\frac{\Omega_X^{n-1}}{\kxf}, K_X) 
@>>> \H^{1,1}(X) @>>> \H^1(X,\kf) \\
@AAA @AAA @| @AAA \\
\H^0 (X, \kf\iy) @.\ra \H^1(X, \oxf) 
@>>> \H^{1,1}(X) @>>> \H^1(X,\kf\iy)
\end{CD}
$$
which is in duality with,
$$
\begin{CD}
\H^n (X, \kf) @.\la \H^{n-1}(X, \frac{\Omega_X^{n-1}}{\kxf}) 
@.\la  \H^{n-11,n-1}(X) @.\la \H^{n-1}(X,\kxf) \\
@VVV @VVV @| @VVV \\
\Ext_X^n (\iy, \kxf) @.\la \Ext_X^{n-1}(\oxf, K_X) 
@.\la  \H^{n-1,n-1}(X) @.\la \Ext_X^{n-1}(\iy,\kxf)
\end{CD}
$$
In particular by way of their images in $\H^{1,1}$,
respectively $\H^{n-1,n-1}$, we have a pairing,
$$\Ext_X^1(\frac{\Omega_X^{n-1}}{\kxf}, K_X) \times
\Ext_X^n (\iy, \kxf) \longra \bc$$
The functoriality of duality ensures that the
pairing depends only on the images (defined
via the obvious diagram chase) in:
$$\Coker: \H^0 (X, \kf)\ra \H^0 (\kf\otimes\cO_Y),
\,\,\,
\Ker: \Ext_X^n(\cO_Y.\kxf) \ra \H^n (X, \kf)$$
Consequently, irrespectively of how we lift, say:
$$\tilde{\c1}: \Ext_X^1(\frac{\Omega_X^{n-1}}{\kxf}, K_X) \longra \H^0 (\kf\otimes\cO_Y)$$
from the co-kernel to $\H^0$, the functoriality of
duality ensures that this latter pairing may be
expressed in terms of I.1.2 by way of,
$$\Ext_X^1(\frac{\Omega_X^{n-1}}{\kxf}, K_X) \times
\Ext_X^n (\iy, \kxf): (L,\phi)\mpo\sum_Y \RES_Y (\phi\tilde{\c1}(L))$$
where the sum is taken over all connected components
of the singular locus to emphasise the local nature of
the latter. The mild ambiguity in the definition of
$\tilde{\c1}$ is genuine as the following example
should illustrate:

\noindent{\bf I.1.3 Example} Away from the singularities
an infinitesimal descent datum for a coherent sheaf $\cE$
on $X$ to the champ classifiant $[X/\cF]$ of the 
foliation is simply a connection along the leaves,
{\it i.e.} a map,
$$\nabla:\cE\longra\cE\otimes_{\cO_X}\kf$$
satisfying the Leibniz rule with respect to $\pa:
\cO\ra\kf$. For any vector bundle the obstruction
to the existence of such a connection on all
of $X$ lies in
a relative Atiyah class,
$${\mathrm at}_\cF (E)\in \H^1(X, \End(E)\otimes\kf)$$
As such a line bundle $L$ admits a leafwise
connection as soon as the image from
$\H^{1,1}$ to $\H^1(X, \kf)$ of its chern
class vanishes. Given such a connection
we obtain a global section,
$$\tilde{\c1}(L,\nabla)\in \H^0(X, \kf\otimes\cO_Y):\ell_{\a}\mpo\frac{\nabla\ell_{\a}}{\ell_{\a}}$$
where $\ell_{\a}$ is a local generator of $L$ around 
a geometric point $y$ in some opens $U_{\a}$ which
cover $U$. Of course the connection is only unique
up to,
$$\H^0 (X, \End(E)\otimes\kf)$$
for a vector bundle, whence the above ambiguity 
modulo $\H^0(X, \kf)$.

Since $X$ is smooth, general local co-homology
considerations imply that to give such a
connection is equivalent to giving it in 
co-dimension 2, whence:

\noindent{\bf I.1.4 Definition/Fact} Denote by $\Pic^0 (\cF)$ 
the group of line bundles on $X$ with infinitesimal
descent data to $[X/\cF]$ away from the singularities,
{\it i.e.} bundles with leafwise connection as per
I.1.3. In particular while such bundles are flat 
along the leaves they may descend only locally
rather than globally to $[X/\cF]$. In addition let
$M(X/\cF)$ be the (possibly empty) space of invariant
measures without support on $Y$, then the images
of these groups in $\H^2$, respectively $\H_2$, admit
the residue formula,
$$(L,\dm)\mpo \int_X  \c1_1 (L)\dm =\sum_Y \RES_Y(\tilde{\c1}(L,\nabla),\dm)
=-\sum_Y \lim_{\e\ra 0}\int_{\pa Y_{\e}}\tilde{\c1}(L,\nabla)\dm$$
which is independent of the connection $\nabla$.

To which we may conclude with the most pertinent
example,

\noindent{\bf I.1.5 Further Example} Away from the
singularities $\kxf$ is genuinely invariant. Dually
to I.1.1 the fact that $[X/\cF]$ will most likely
have no \'etale neighbourhoods at the singularities
implies that even for this bundle it is ambiguous
as to whether it really is a bundle on $[X/\cF]$
(depending on how we define this at the singularities)
or not. It does, however, have a canonical connection,
$$d=\nabla:\kxf\longra K_X=\kxf\otimes\kf$$
induced by exterior differentiation. In turn, around
the singularities $\pa:\cO_X\ra\kf\iy$ affords a
linear map,
$$\ox\otimes\cO_Y\longra\ox\otimes\cO_Y(\kf)$$
whence Tr$\pa$ is a section of $\H^0(X, \kf\otimes\cO_Y)$
and,
$$\tilde{\c1}(\kxf,\nabla)=-{\mathrm Tr}(\pa)$$.

\subsection*{I.2 Local vs Global}

There follows a series of remarks intended to
illustrate the difficulty of calculating the
residue I.1.2 as soon as the singular locus
is positive dimensional. Of themselves they
are in-essential, but should aid understanding.
To this end recall, \cite{sga2}, the particularly satisfactory
description of the duality theorem for a
complete regular local ring $\hat{\cO}$ over
a field $k$. Supposing the dimension is $n$ we
may suppose $\hat{X}=\Spf \hat{\cO}$ is the
completion of $X=\ba_k^n$ in the origin $x$,
and we have a duality,
$$\H^0(\hat{X},  \cO_{\hat{X}})\times \H_x^n(X,\o_X)\longra k$$
The latter group is,
$$\lim_{n\ra\infty} \Ext_X^n (\cO/\gmi^n,\o_X)$$
and the pairing is the Grothendieck residue.
In terms of coordinates, $z_1,\hdots , z_n$
over $\bc$ its elements may be identified with,
$$f(z_1^{-1},\hdots ,z_n^{-1})dz_1\hdots dz_n\otimes\g_1\otimes\hdots\otimes\g_n$$
for $f$ a polynomial in $n$ indeterminants, and $\gamma_i$
the $\d$-function of a loop around $x$ in the
$i$th direction. As such the duality is the
Cauchy residue,
$$(g,f)\mpo\oint_{\g_1\times\hdots\times\g_n}
\tilde{g}(z_1,\hdots , z_n) f(z_1^{-1},\hdots ,z_n^{-1})dz_1\hdots dz_n$$
where $\tilde{g}$ is any convergent ({\it e.g.}
polynomial) function which agrees with the
formal function $g$  up to the maximal degree
of $f$. Alternatively one may express this
as an integral over a co-dimension 1 real
hypersurface, {\it cf.} \cite{gh} \S 5.

In contrast to this the first non-local
case of I.1.2 occurs in dimension 3. As
it happens (and we'll see in the next
section), one can quite quickly
reduce to a slightly different sub-scheme
of $Y$ which is LCI, albeit since we're
just sketching the difficulty one may
as well assume this for $Y$. In any case
if $Y$ is LCI of co-dimension $d$, then
the only local Ext groups are in 
dimension $d$, {\it i.e.}
$$\cE xt^q(\cO_Y, E)= 0,\,\, q\neq d$$
for any locally free $E$. In particular,
the local global spectral sequence for
Ext is degenerate, and for $U$ as per
the last section:
$$\Ext_U^n(\cO_Y.\kxf)= \H^{n-d} (U, \cE xt^{d}(\cO_Y.\kxf))$$
In dimension $n=3$, $d=2$ this is particularly
simple, and doesn't involve more than
a couple of lines diagram chasing, {\it i.e.}
by the injectivity of distributions we have
a diagram with exact rows:
$$
\begin{CD}
 \Hom_U (\iy, \kxf\otimes\cD_X^{0,2}) 
@.\twoheadleftarrow \Gamma (U, \kxf\otimes\cD_X^{0,2}) 
@.\hookleftarrow \Hom_U (\cO_Y, \kxf\otimes\cD_X^{0,2})  \\
@VV{\pab}V @VV{\pab}V @VV{\pab}V \\
\Hom_U (\iy, \kxf\otimes\cD_X^{0,3} )
@.\twoheadleftarrow \Gamma (U, \kxf\otimes\cD_X^{0,3}) 
@.\hookleftarrow \Hom_U (\cO_Y, \kxf\otimes\cD_X^{0,3}) 
\end{CD}
$$
where one starts with $\dm$ in the top left,
lifts to some $T$, so that $\pab T$ is the
global class in $\Ext^3$. On the other hand,
there is no local $\cE xt^3$, so on a cover
$\coprod_{\a} U_{\a}\ra U$, we find on $U_{\a}$
distributions $S_\a$ supported in $Y$ such
that $T_\a = T + S_\a$ is $\pab$ closed, and
$S_{\a\b}=S_\a-S_\b$ is the desired class in
$\H^1 (U, \cE xt^{2}(\cO_Y,\kxf))$. In particular,
there is absolutely no local obstruction to 
lifting $\dm$ in a way which is $\pab$-closed.
As such, it is not at all a question of how
one does the integral (in an indefinite sense) locally, 
but rather how these liftings patch, {\it i.e.}
in the notation of the last section, the 
residue, supposing spt$\rho_\a\sbs U_\a$ compact, is:
$$\sum_\a \pab\rho_\a\o_\a T = \sum_\a \pab\rho_\a
(\o_\a T_\a-S_\a)= -\sum_\a \o\pab (\rho_\a S_{\a\b})$$

Alternatively, let us bring the foliation into play.
A non-trivial case, albeit extremely special and
improbable, {\it cf.} IV.1, occurs when on some
cover $U_\a$ we can find coordinates such that
the foliation has the form,
$$\pa_\a = z_\a\frac{\pa}{\pa z_\a} + x_\a^p \frac{\pa}{\pa y_\a}$$
and we suppose (actually it follows) for further simplicity that the
$z_\a$ patch to an invariant centre manifold
$Z$, so $\o_\a=z_\a^{-1}dz_\a$ is a perfectly
good lifting of the trace of I.1.5.  Observe
that while the sum $\tilde{\tau}$ over $\a$ of the $\pab\rho_\a\o_\a$ is integrable,
in fact its even smooth, the individual
$\pab\rho_\a\o_\a$ are not necessarily integrable,
and one cannot interchange the order of 
summation with integration. As such, to
fix ideas, say with respect to a different
index set $i$ (possibly the same) subordinate to
a cover as above so that we have good coordinates
$x_i,y_i,z_i$ we triangulate by way of $V_i$
along $Y$ rather than use bump functions, then
we have a triangulation $\tilde{V}_i$ of $U$, and:
\begin{align*}
\int_U \tilde{\tau}\dm = \sum_i \int_{\tilde{V}_i} \tilde{\tau}\dm=
\sum_i \lim_{\e_i\ra 0} \int_{\substack{|z_i|=\e_i\leq |x_i|\\ y_i\in V_i}} \o \dm
+ \int_{\substack{y_i\in\pa V_i\\ \e_i\leq |x_i|, |z_i|}}\o\dm
\end{align*}
The first of the integrals on the rights is easily
seen to be bounded by the Segre class around the
singularities. The latter integral ``should'' cancel
according to, 
$$\int_{\substack{y_i\in\pa V_i\\ \e_i\leq |x_i|, |z_i|}}\o\dm
=\sum_j \int_{\substack{y_i\in\pa V_i\cap\pa V_j\\ 
\e_i\leq |x_i|, |z_i|}}\o\dm$$ 
with the sum being taken over edges of the triangulation,
so that $(\pa V_i)\cap V_j$ and $(\pa V_j)\cap V_i$
have the opposite orientation. One cannot, however,
a priori arrange that the domains $|x_i|\geq \e_i$
and $|x_j|\geq \e_j$ are in fact the same. Plainly
this is the difficulty of interchanging 
integration with summation for the $\pab\rho_\a\o_\a$
that we've already noted in different clothes.
It is also, and more helpfully, a problem of
the holonomy around $Y$ of the induced foliation
in the centre manifold.

Now while the appearance of an essentially
topological obstruction such as the holonomy
may be rather encouraging, one should observe:
\begin{itemize}
\item The coordinate $x_\a$ may be supposed
to define an algebraic hypersurface, in fact 
an exceptional divisor after a blow up. The
existence domain of \cite{mp1}, however, for
such coordinates is only in sectors of the
argument of $x_\a$ up to $\pi/p$, and this
is shown to be best possible.
\item Sectors of width $\pi/p$ are uselessly
small, {\it i.e.} the centre manifold is 
hopelessly non-unique  in such sectors, and
irrespectively of how small the neighbourhood even
the local behaviour of leaves is utterly out
of control as the argument goes from
$\pi/p -\e$ to $\pi/p +\e$.
\item In particular, the holonomy about the
induced foliation in the centre manifold has
nothing more than a formal sense.
\end{itemize}

In terms of co-equalisers here is another
description of the same phenomenon. Given
that we're only describing the difficulty we'll
be a bit informal, but cf. \cite{MeOnMeasures}
for a precise discussion of invariant measures
and co-equalisers. In any case by \cite{mp1},
and the current \S 2, one may reduce to the
case where $U\bsh E$- $E$ the exceptional
divisor resulting from blowing up- is covered
(this will almost certainly employ sectors
as above) by finitely many opens $W_a$ such 
that on each $W_a$ the foliation is a 
fibration $W_a\ra T_a$ whose fibres do
not adhere in $E$, whence, a priori, not 
in the singular locus either. This implies
that on every fibre $\pi_a^{-1}(t)$,
$\o$ is integrable, while $\dm$ descends
to a measure $\dm(t_a)$. Consequently:
$$\int_{W_a} \tilde{\tau}\dm=
\int_{T_a}  \dm(t_a)
(\pi_a)_* (\tilde{\tau}) =
\int_{T_a}  \dm(t_a)
(\pi_a)_* (\o|_{\pa W_a})$$
At which point taking $T=\coprod_a T_a$
one may show that the function,
$$w= \coprod_a (\pi_a)_* (\o|_{\pa W_a}) \in C^0 (T)$$
is a co-equaliser, {\it i.e.} there is a Borel
function $f$ on the arrows of $(s,t):R\rras T$
such that,
$$w= (s_* -t_*)(f)$$
Now in principle, and by definition, the value
of an invariant measure on co-equalisers should
be zero, but this is only true if $s_* f$ and
$t_* f$ are integrable. Again this is a manifestation
that $\o$ may not be integrable on every face
of a triangle, even though it is so leaf by
leaf, and the sum over the faces is integrable.
Of course if $Y$ were $\bp^1$, and the singularities
really as simple as this toy example, $f$ would
be integrable, but such a case is not only
trivial, but, by many ways, may be shown not 
to occur.

As a final caveat to trying to do the problem
locally, one should bear in mind the simple
sub-case where the centre manifold is everywhere
convergent along $Y$. Whence, locally, there
are absolutely no obstructions on an invariant
measure supported in the same, and no reason
whatsoever for any unbounded function of $x_\a$
to be integrable, let alone $x_a^{-p}$. Of course,
there are global obstructions for the existence
of such an invariant  measure, which in a certain
sense is what we'll use, but, as we've said, by
\cite{mp1} II.1.5, we start from the situation
where the centre manifold simply fails to exist
on a sector which is sufficiently large to be
useful.

\subsection*{I.3 Isolated Residues}

Supposing resolution of singularities (a theorem in
dimension 3, \cite{mp2}), one may show that one
never needs to tackle an isolated residue which
is more complicated than,

\noindent{\bf I.3.1 Set Up} Let $U$ be a neighbourhood
of an isolated foliation singularity, such that on
completion in the singularity a (convergent) generator
of the foliation has a Jordan decomposition,
$$\pa=\pa_S + \pa_N,\,\,\, \pa_S=\sum_{i=1}^n \lb_i \hat{z_i}
\frac{\pa}{\pa \hat{z_i}},\,\,\,\, \lb_1\hdots\lb_n\neq 0$$
and the remaining (formal) coordinate $\hat{x}$
may be supposed such that $\hat{x}=0$ is not
only convergent, but the exceptional divisor
after appropriate blowing up, and $\hat{x}^{-p-1}\pa \hat{x}$
is invertible, even $1+\gmi$, for some $p\in\bn$.

Before proceeding, let us note,

\noindent{\bf I.3.2 Caveat} Already in dimension 3,
even starting from a foliated projective 3-fold $X$,
the resolution $\pi:\cX\ra X$ takes place in the
2-category of champs de Deligne-Mumford. This means
that the above $U$ will only be an \'etale neighbourhood
of $\cX$, or equivalently for some finite group $G$
(actually $\bz/2$ in dimension 3) the champ
classifiant $[U/G]$ is open in $\cX$. Manifestly,
however, this is irrelevant to the residue 
calculation which may be safely performed in $U$.

Now take coordinates $z_i, x$ such that the Jordan
decomposition holds modulo $\gmi^N$ (in fact even
modulo $\cO(-NE)$ by way of blowing up) for some
large $N$ to be chosen, and suppose that we have
an invariant measure with no Segre class around
the exceptional divisor, or, even, just the point, then:
$$\oint_{|x|=\e, |z_i|\leq \d} \frac{dx}{x}\dm
+ \sum_i \oint_{\substack{|z_i|=\d, |z_j|\leq \d\\ \mxle}}
(\frac{\pa x}{x}\frac{z_i}{\pa z_i} )\frac{dz_i}{z_i}\dm =0$$
Next observe that on the face $|z_i|=\d$,
$$|\pa z_i| \gg \d -\e^N C$$
for some constant $C$ provided $\pa_N$ vanishes to
first order, which can always be achieved by
blowing up. Whence $\e=\d$ is a good
choice, albeit we make:

\noindent{\bf I.3.3 Remark} The important thing here is
to have the modulus of $(\pa z_i)^{-1}z_i$ bounded
above, and $\d=\e^\a$ any $\a\in\br_{>0}$ will
work provided $N$ is chosen sufficiently large.
In the context of Segre classes, here Lelong
numbers, this amounts to the same with weights,
but the vanishing of the Segre class without
weights implies the vanishing of that with weights.
Alternatively blow ups in weighted centres may be
resolved by a sequence of modifications in smooth
centres. Irrespectively, this kind of estimate is
extremely robust, so for ease of notation we will
often simply take $\d=\e$ even though we occasionally
use $\d=\e^\a$ for appropriate $\a$ in \S IV.3/IV.4,
and do so without comment.

This said, we therefore have that,
$$\oint_{|x|=\e, |z_i|\leq \e}\frac{dx}{x}\dm\ll \e^p \mathrm{s}_{Z,\dm} (\e)$$
where the symbol $\mathrm{s}_{Z,\dm} (\e)$ is the part of
the Segre class/Lelong number arising from the faces $|z_i|=\e$,
at radius $\e$, which is plainly the same thing as
the Segre class/Lelong number up to $o(\e )\ra 0$.
The residue of any 1-form $\o$ that we have to
calculate has the property,
$$\o\dm= f_i \frac{dz_i}{z_i}\dm,\,\, 
\o\dm = g \frac{dx}{x^{p+1}}\dm =0$$
on the $|z_i|=\e$, respectively $|x|=\e$, faces
with $|f_i|, |g|$ bounded, and so:

\noindent{\bf I.3.4 Fact} Let things be as in I.3.1, $I$
the ideal of the singularity, and suppose (as we
may after further blowing up) that $\pa_N$ is zero
to first order, {\it i.e.} $\pa \in \gmi^2\Der(\cO)$,
then,
$$\RES:\Ext_U^{n-1}(I,\kxf ) \iso \Ext_U^n(\cO/I,\kxf)$$
vanishes on $\dm$.

\subsection*{I.4 Logarithmic holonomy}

Like \S 1.2 the current section is in-essential,
but it explains the holonomy content of the final
section, \S IV.4. It is, therefore, useful reading
either pre or post the same. It also shows how
super attracting dynamics occurs for foliations.
Our set up is as follows:

\noindent{\bf I.4.1 Set Up} Let $Y$ be a reduced
connected 1 dimensional analytic space whose
components $Y_i$ are smooth, and whose singularities
are at worst nodes. Around $Y$ we have an analytic
surface $S$ containing $Y$ as a tubular neighbourhood,
and $\wh{S}$ will be the completion of $S$ in $Y$.
The surface $S$ is foliated by way of $S\ra [S/\cF]$,
all the components of $Y$ are supposed invariant,
so, certainly, there are singularities at the
intersections of components. Not just these
singularities, however, but, for convenience, any other
singularities on $Y$ will be supposed reduced in the
sense of Seidenberg, albeit this hypothesis is really
only necessary at the intersections of components.

Now take a component $C$ of $Y$ with $C^*$ the
complement of the foliation singularities on
the same. For $T$ a 1-dimensional germ of an
analytic space we have the holonomy representation,
and its completion,
$$h:\pi_1 (C^*)\ra \Aut (T),\,\, \hat{h} :\pi_1 (C^*)\ra \Aut (\hat{T})$$
Consider the exponential, $\exp:\gh\ra T:\tau\mpo \exp(\tau)=t$
from the germ of the former around $-\infty$. 
Understanding by $\ghi$ germs of diffeomorphisms
preserving $-\infty$ we can define liftings of
the above representations as follows: in the 
first place at the singularities on $Y$ that
aren't nodes there is a further branch of the
foliation, and this may be supposed convergent,
so, say, an invariant divisor $D$ on $S$. This
may fail at nodes, but remains true formally,
and, since the latter is our real interest we'll
suppose it analytically for ease of exposition.
Consequently for $Y_C$ the part of $Y$ omitting
$C$, with $S_C$ a tubular neighbourhood of $C$,
$C^*$ is homotopic to $S_C\bsh (D+Y_C)$, and
for $V\ra C^*$ the cover afforded by the holonomy
representation inducing $\Sigma$ over $S_C\bsh (D+Y_C)$
we have a diagram,
$$\begin{CD}
V\sbs\Sigma @>>t> T \\
@VVV @. \\
C^* \sbs S_C\bsh (D+Y_C)
\end{CD}
$$
with horizontal map an \'etale covering of the pair.
This implies that the fundamental group of
$S_C\bsh (D+Y)$ is canonically an extension,
$$\begin{CD}
1@>>>\tm@>>>\pi_1(S_C\bsh (D+Y))@>>> \pi_1(C^*)@>>> 1
\end{CD}$$
or alternatively, if we make no hypothesis on the
singularities around $C$, and take $C^{\times}$
to be a suitable neighbourhood of $C^*$ punctured
in the same, then we have,
 $$\begin{CD}
1@>>>\tm@>>>\pi_1(C^{\times})@>>> \pi_1(C^*)@>>> 1
\end{CD}$$
In either case there are liftings,
$$\log h: \pi_1(C^{\times})= \pi_1(S_C\bsh (D+Y))
\longra \ghi$$

Formally the situation is more problematic than
is desirable, since the ``fundamental group of
a punctured formal disc'' has no sense except in
a pro-finite way, even though one has a perfectly
good non-profinite theory of fundamental groups
of analytic formal schemes. Nevertheless, in the
first instance, elements of $\hat{h}$ admit 
logarithms of the form,
$$\log \hat{h} = \tau + \lb + f(e^{\tau})$$
for $f$ a formal series, $\lb\in\bc$, and these
can be composed according to the obvious rules,
so, what is true, is that we get an extension,
 $$\begin{CD}
1@>>>\tm@>>>\wh{\Log}(C^{\times})@>>> \wh{\Hol}(C^*)@>>> 1
\end{CD}$$
for $\wh{\Hol}$ the image of $\hat{h}$. In the 
situation being discussed where $S$ is
convergent this subtlety is a bit irrelevant,
but if only a formal surface were given (which
would have been the case had we intended to
really use the present discussion) it's rather
critical since $\pi_1(S_C\bsh (D+Y_C))$ would have
its obvious sense, but $\pi_1(S_C\bsh (D+Y))$
would be rather problematic.

In any case, these definitions have so far
done nothing except replace the action of
the holonomy group $\Hol(C^*)$ on $T$ by that
of $\Log(C^{\times})$ on $\gh$, and the germ of
the champ classifiant $[S/\cF]$ at the
geometric point determined by $C^*$, {\it i.e}
$[T/\Hol(C^*)]$, by the groupoid $[\gh/\Log(C^{\times})]$
which is equivalent to $[T^{\times}/\Hol(C^*)]$,
$T^{\times}$ being $T$ punctured in the origin.
Now consider passing from $C$ to $D$ at some
singularity, say:
$$\pa = x \dx -\lb y \dy,\,\,\Re(\lb) > 0,\,\, D=(y=0),\, C=(x=0)$$
and everything convergent to fix ideas, with transversals
$T_C$, $T_D$ to $C$ and $D$ respectively at $y$, respectively
$x$, equal to 1. In a minor abuse of notation, put:
$$t=xy^{1/\lb},\,\, \tau=\log t,\,\, 
s=yx^{\lb},\,\, \s=\log s$$
and observe that we have inclusions,
$$
\begin{CD}
\Log(D^{\times}) @<j<< \tm^2 @>i>>
\Log(C^{\times}) 
\end{CD}
$$
where for generators $a,b\in\tm$,
\begin{align*}
i(a)&=\tau+a & j(a)(\s)&=\s+\lb a \\
i(b)&=\tau+b/\lb & j(a)(\s)&=\s+b
\end{align*}
or equivalently one conjugates the
representation in $\tau$ to that in
$\s$ by way of $\s=\lb\tau$, and the
condition $\Re(\lb)>0$ ensures that
neighbourhoods of $-\infty$ in $\tau$
go to the same in $\s$. As such we may
extend through the singularity in the
obvious Seifert-Van Kampen-esque way, 
{\it i.e.}
$$\Log(C^{\times}\cup D^{\times}):=
\Log(C^{\times})*_{\tm^2}
\Log(D^{\times})$$
which of course is perfectly consistent
with,
$$\pi_1(C^{\times}\cup D^{\times}):=
\pi_1(C^{\times})*_{\tm^2}
\pi_1(D^{\times})$$
and the representation $h$ extends, or,
similarly at the level of $\wh{\Log}$
alone if our context were purely formal.
It is furthermore true that this
correctly reflects the dynamics of
the foliation, {\it i.e.} the composition
$\gh_C\ra T_C\ra [X/\cF]$ is \'etale,
and we get a diagram,
$$
\begin{CD}
H_C @<<< R \\
@VVV @VVV  \\
[C^{\times}\cup D^{\times})/\cF ] @<<< H_C
\end{CD}
$$
such that the germ of the \'etale groupoid
$R\rras \gh_C$ around $-\infty$ is,
$$\Log(C^{\times}\cup D^{\times})\times\gh_C
\rras \gh_C$$
with the above action. Plainly this 
doesn't make sense if $\Re(\lb)\leq 0$, or
the singularity is a node. Equally plainly
the hypothesis that the singularity is
analytically linearisable isn't actually
necessary and just the given form to 1st
order will do, provided we take the actual
conjugation between $\gh_C$ and $\gh_D$
given by the foliation. The formal situation
is both easier, and more difficult, {\it i.e.}
for $\lb$ irrational the singularity is
formally linearisable and everything is
as above, otherwise in the notation post
III.4.1(c),
$$s^l=t^k\exp (\sum_{m=1}^{\infty} t^{krm} R_m(\log t))$$
for $R_m$ polynomial, and we can continue
to give a formal sense to maps preserving
$-\infty$, and how to compose them. 

Next we form the dual graph $G$ of $Y$,
or more correctly the dual graph of the
sub-curve were 
all the singularities at intersections
of components is as above, {\it i.e.}
a vertex for each component and an edge
whenever the singularity has the above
form at the intersection of components.
In particular for each edge $e$ there
are mutually inverse automorphisms,
$\phi_{+-}(e)$, $\phi_{-+}(e)$, in
$\ghi$ such that,
$$\tau_-=\phi_{-+}(\tau_+),\,\, 
\phi_{+-}(e):v_-longleftrightarrow v_+: \phi_{-+}(e),
\,\, \tau_+=\phi_{+-}(\tau_-)$$
describes the conjugation from a
logarithmic transversal at one 
vertex to the other. A loop
$\g$ in $G$ is simply a series
of directed edges bringing us
back to the same vertex, so, in
a minor abuse of notation, we
have a representation,
$$\theta:\pi_1(G)\longra\ghi:
\g\mpo\prod_{e\in\g}\phi(e)$$
Let us abuse notation a bit 
further to put things together,
{\it i.e.} suppose $Y$, $G$
connected and let $Y^{\bullet}$
be the union over components, $C$,
of the $C^*$ together with the
singularities corresponding to
the edges of $G$, and $Y^{\times}$
a neighbourhood of $Y^{\bullet}$
punctured in the same, then we
have a diagram of representations,
$$
\begin{CD}
\pi_1 (Y^{\times}) @>{\log h}>> \ghi \\
@VVV @| \\
\pi_1 (G) @>\theta>> \ghi
\end{CD}
$$
and this action accurately reflects
the holonomy groupoid defining the
foliation $S\ra [S/\cF ]$, {\it i.e.}
the germ of the induced \'etale
groupoid $R\rras\gh$ about $-\infty$
arising from the restriction to $Y^{\times}$
is the action of the image $\Log(Y^{\times})$
of $\log h$ on $\gh$. Similarly in the
formal case, albeit purely in terms of
a subgroup $\wh{\Log}(Y^{\times})$ of
automorphisms. By way of clarity let
us offer,

\noindent{\bf I.4.2 Example} Take for
$Y$ an elliptic gorenstein foliation
singularity, {\it cf.} \cite{canmod}, with
more than one components, with say 
$-n_i=Y_i^2$ the self intersection of
the components. The components $Y_1,\hdots ,Y_d$,
each isomorphic to $\bp^1$,
are ordered as a cycle: $Y_i$ meeting
only $Y_{i-1}$, $Y_{i+1}$ at $\circ_i$,
respectively $\infty_i$- $d+1=1$- and
for $y_i$ a coordinate along $Y_i$
at $\circ_i$ with $x_i$ the normal the
foliation has a generator $\pa_i$ around the
affine line through $\circ_i$,
$$\pa_i = y_i\frac{\pa}{\pa y_i} 
-\lb_i  x_i\frac{\pa}{\pa x_i}$$
where the above are related by,
$$x_{i+1}=y_i^{-1},\,\, y_{i+1}=y_i^{n_i}x_i,\,\,
\lb_{i+1}=(n_i-\lb_i)^{-1},\,\, \pa_{i+1}=\lb_i\pa_i$$
As such, for $t_i=x_iy_i^{\lb_i}$ we have a cycle of transformations
of half planes,
$$
\begin{CD}
\gh_{i-1},\, \tau_{i-1}=\log t_{i-1}
@>{\lb_{i-1}}>>
\gh_i,\, \tau_i=\log t_i
@>{\lb_i}>>
\gh_{i+1},\, \tau_{i+1}=\log t_{i+1}
\end{CD}
$$
while each $\pi_1(Y_i^*)$ is (canonically)
$\tm^2$. The topology being as simple as
possible the cycle $\bz=\pi_1(G)$ acts on
$\tm^2$ and $\pi_1(Y^{\times})$ is a 
semi-direct product,
$$\begin{CD}
1@>>>\tm^2@>>>\pi_1(Y^{\times})@>>>\bz@>>>1
\end{CD}
$$
the so called ``isotropy group of the cusp'',
whose representation in $\ghi$ is faithful,
and in terms of generators $a$, $b$ of $\tm^2$
is given by the matrices,
$$
\left[ \begin{array}{cc}
1 & a \\
0 & 1 \end{array} \right],\,\,\,\,
\left[ \begin{array}{cc}
1 & \lb_1 b \\
0 & 1 \end{array} \right],\,\,\,\,
\left[ \begin{array}{cc}
\lb_1\hdots\lb_d & 0 \\
0 & 1 \end{array} \right]
$$
and $\lb_1\hdots\lb_d$ is an automorphism
of the image of $\tm$ in this representation.
Indeed, each $\lb_i$ is a quadratic
irrationality, {\it e.g} we have the
continued fraction,
$$
\lb_1=
\cfrac{1}{n_d-\cfrac{1}{n_{d-1}-\cfrac{1}{\hdots-\cfrac{1}{n_2-\lb_1}}}}
$$  
In particular, to see that the representation
is faithful one observes, \cite{vandergeer} II.2.1,
that while the implied basis $a$, $b\lb_1$ of the
kernel to the cycle is not canonical, there is
an infinite set of such choices, which is canonical,
and on which $\bz$ acts faithfully by way of
$\lb_1\hdots\lb_d$. Furthermore, identifying,
as above, the representation of $\pi_1(G)$,
whence of $\pi_1(Y^{\times})$, with a cyclic
subgroup of $\br_{>0}^{\times}\sbs \bg_m$ this
is also the representation afforded by the
log canonical bundle $K_{S/\cF}$ which, while
trivial on each $Y_i$ is globally not so. Indeed
dually to the $\pa_i$ one has 1-forms,
$$\o_i = \lb_i \frac{dy_i}{y_i} + \frac{dx_i}{x_i},\,\,\,\,
\o_i=\lb_{i+1}\o_{i+1}$$
As such, elliptic gorenstein is a slight misnomer,
{\it i.e.} although $K_S + Y$ is trivial on the
cycle and descends to a bundle on the contraction,
since,
$$K_S = \kf + K_{S/\cF} $$
and the latter term isn't even trivial as a bundle
on $Y$, let alone $S$, by way of the above, the
canonical $\kf$ is never $\bq$-Gorenstein on the
contraction despite being torsion on a neighbourhood
of each component $Y_i$.

While straightforward this example is rather instructive.
Specifically: the champ classifiant $[Y^{\times}/\cF ]$
has the logarithmic holonomy not just as a germ, but
it is $[\gh/\Gamma]$ for $\Gamma$ a finite sub-group
of $\Aut(\gh )$ (not just $\ghi$) with, say, generators
as above. As such, it certainly admits at least one
invariant measure, {\it i.e.} the Poincar\'e metric
on $\gh$. This, in turn, admits the potential,
$\log|\log|t_1|^2|$ which, while invariant by
the action of the $\tm^2$ is not so by the cycle
of the graph which acts in a super attracting way, 
$$\log|t_1|^2\mpo (\lb_1\hdots \lb_d) \log|t_1|^2$$
In particular the assertion, valid for the usual
holonomy, that, at least formally, it must take
values in $S^1$ in the presence of an invariant
measure is false for the logarithmic since in the
logarithmic variable the cycle $\pi_1 (G)$ even
has, already at first order, 
the faithful representation in 
$\bg_m$,
$$\g\mpo\prod_{e\in\g} \lb(e)$$
for (directed) multipliers $\lb(e)$ at the edges
as above. Let us, therefore, conclude by way of,

\noindent{\bf I.4.3 Remark/Summary} The problem
of the non-triviality of the representation of
$\pi_1(G)$ appears in a closely related form in
the final section \S IV.4. The difficulty that it
poses is not dealt with by way of the relevant
(formal) logarithmic holonomy, but by a more
algebraic alternative. Ideally it would have been
treated as an extension of the basic construction
III.1.1 to the logarithmic case (notice, by the
way, there being no contradiction between III.1.1
and any of the above since logarithmic co-equalisers
must be calculated logarithmically, {\it i.e.}
with subsets of half planes). Unfortunately,
unlike the current ideal case which is akin to
the measure being supported in the centre
manifold, \S III-IV are properly ``almost holonomy"
rather than holonomy. In particular there is
always an error associated with the multiplicities
$p(v)$ at the vertices of $G$ as encountered in
\S IV.4, and this error is potentially very non-uniform
from vertex to vertex. As a result, I couldn't
get what should have been the right estimate,
{\it i.e.} up to the error, cf. \S IV.4, $\e^{p(v)}o(\e)$,
the measure is radially symmetric. Such an 
estimate would have rendered the global residue
calculation purely topological, and achieved
wholly by the method of ``almost holonomy",
as opposed to the more algebraic considerations
which are employed to treat the potential 
difficulty of many singularities of type III.4.1 (b)/(c).

\newpage

\section*{II. Particular Singularities}

\subsection*{II.1 Generically Transverse}

We begin with an investigation of certain
singularities where the relation between
the formal and analytic structure is
particularly good. In the first instance
so much so that there is no need to
restrict the dimension, {\it i.e.}

\noindent{\bf II.1.1 Set Up} In a
neighbourhood $U$ of a $2+n+m$ dimensional
complex space, $m,n\in\bn\cup\{0\}$, about
a singular point of the foliation in the
completion $\hat{U}$ in the singular locus
the foliation admits a formal generator of
the form,
$$z\dz + \dfrac{y_1^{q_1}\hdots y_n^{q_n}x^{p+1}}{p(1+\nu(y,w)x^{p})}\dx$$
where $x,y_i,z,w_j$ are formal coordinates,
$p,q_i\in\bn$, albeit $1\leq i\leq n$, $1\leq j\leq m$,
so the $y_i$, respectively $w_j$, and whence
the $q_i$, are suppressed should $m$ or $n$
be zero.

Under these hypothesis, and possibly after
blowing up in the singular locus, a 
sufficiently small open neighbourhood of the
singularity admits for any value of the
arguments of $x$ or $y_i$ an open neighbourhood
of the form, $V=S\times S_1\times\hdots\times S_n\times\D\times\D^m$,
{\it i.e.} discs for each of the coordinates
$z,w_j$, and sectors $x\in S$, $y_i\in S_i$
such that $\xi=x^{-p}y_1^{-q_1}\hdots y_n^{-q_n}$
has an argument up to $2\pi$ and a branch
within $\pi/2$ of $-1$ such that: on the same
we may find an analytic generator of the
foliation together with a conjugation to
the normal form II.1.1, \cite{mp1} VI.3.3(c).
Given this we may easily integrate the foliation.
In the first place by way of a conformal change
of variable in $S$ we may suppose $\nu=0$, and
so obtain invariant functions,
$$s=z\exp(\xi),\,\, y_i,\,\, w_j\,\,\, 1\leq i\leq n,\,\, 1\leq j\leq m$$
By way of notation, let us denote the invariant
functions $y_i, w_j$ by  a vector
$\bft$ of such. As a result if $Z_\z$ is
the transversal $z=\z\neq 0$, then $\bft$
fibres $Z_\z$ over $T$, say, with fibre embeddable
on the omission of the exceptional divisor
$E$, $x=0$ or $y_i=0$, in $\bc$ by $\xi$, {\it i.e.}
$$
\begin{CD}
Z^*_\z:=Z_\z\bsh E @>{\xi\times\bft}>> \bc\times T\\
@VV{\bft}V @.\\
T
\end{CD}
$$
where  all the maps are open, and the horizontal
one an embedding. On $\bc$ we have the action of 
$\tm$ by translations, and the groupoid,
$$R\rras Z^*_\z$$
induced by the foliation is the pull-back of the
groupoid given by the action of $\tm$ on $\bc\times T$.
As such while far from essential in the immediate
context let us reduce the question of the existence
of invariant measures to purely fibrewise considerations
in $\bft$, to wit:

\noindent{\bf II.1.2 Lemma} Let $\pi:X\ra B$ be a map
of metric spaces, and $\dm$ a measure on $X$ of finite
mass, then there is a measure $\dn$ on $B$ and a family
of measures $b\mpo d\lb_b$ with support in $X_b=\pi^{-1}(b)$ such
that for any bounded continuous function $f$ on $X$,
$b\mpo d\lb_b$ is bounded and,
$$\int_X f\dm=\int_B\dn\int_{X_b} fd\lb_b$$
\noindent{\bf proof} (for want of a reference) Define
$\dn$ by the obvious formula,
$$\dn(g)=\int_X (\pi^* g) \dm$$
then for any bounded continuous function $f$,
$$g\mpo \int_X (\pi^*g) f\dm$$
is absolutely continuous with respect to $\dn$ and
admits a $L_{\infty}(\dn)$ derivative, $D(f)$. By \cite{feder}
2.9.8 there is a formula for this, {\it viz}: let $b\in B$,
$B_\e (b)$ a ball of radius $\e$ about it and $\nu_b(\e)$
its $\nu$ measure, then,
$$D(f)(b)=\lim_{\e\ra 0} \dfrac{1}{\nu_b(\e)}\int_{\pi^{-1}(B_\e(b))} f\dm$$
which, here, is in fact defined and bounded not just
$\nu$-a.e. but everywhere. $\Box$

We can, and will, therefore, think of $\dm$ as $\dn(b)d\lb_b$,
and, of course, modulo the obvious change of notation $b\mpo\bft$,
we have:

\noindent{\bf II.1.3 Corollary/Definition} For a $R\rras Z^*_\z$
invariant measure $\dm$ described fibrewise by $d\nu(\bft)\dm_\bft$
as per II.1.2, $\dm$ is invariant iff the fibres $\dm_\bft$ are
invariant for $\nu$ almost all $\bft$.

There are, of course, lots of measures on $\bc$ invariant
by translations in $\tm$, so the only way to exclude the 
existence of $\dm$ is by considerations of mass, and this
is plainly the case. Indeed each point in each fibre admits
a neighbourhood, all translations of which are disjoint,
and there are infinitely many such. Whence every $\dm_\bft$
is zero, and since every leaf off the centre manifold $z=0$
meets $Z_\z$ we conclude,

\noindent{\bf II.1.4 Fact} Let things be as in II.1.1,
modulo blowing up in the singular locus to ensure that
$x=0$ or $y_i=0$ is the exceptional divisor and an 
analytic conjugation on the domain $V$ to the normal
form, then an invariant measure for the foliation on
$V$ lies either in the exceptional divisor or the
centre manifold.

In the particular case that $m=n=0$, one may appeal
to \cite{siu} to conclude that the centre manifold
exists in all of $U$. In general, however, such
reasoning is not valid, and we must appeal to the
Segre class around the exceptional divisor $x=0$.
The definition of Segre class is only relevant to
the part of the measure off the divisor, so,
henceforth this precision may be omitted. For any
sub-variety, and in some generality, \cite{uu} \S III,
it is defined as soon as the sub-variety is compact,
but, post factum may be localised. As such, supposing
$U$ a neighbourhood of some $X$ in which the singular
locus is compact, and letting $f=0$ be an equation
for the divisor $x=0$ in all of $U$ with $E$ the
exceptional divisor, the Segre class $\mathrm{s}^X_{E,\dm}$
dominates:
$$\mathrm{s}^U_{E,\dm}:=\lim_{\e\ra 0} \mathrm{s}^U_{E,\dm} (\e) := \oint_{|f|=\e} \frac{df}{f}\dm$$
Now, the centre manifold in $V$ embeds in $\bc\times T$
by way of $f\times\bft$, and $\dm$ descends to a
measure $\dm(\bft)$ in $\bft$- not to be confused
with $\dm_\bft$- whence:
$$\mathrm{s}^U_{E,\dm}=\dfrac{|S|}{2\pi}\int_T \dm(\bft)$$
where $|S|$ is the aperture of the sector $S$.
Consequently:

\noindent{\bf II.1.5 Corollary} Suppose $U$ is a
neighbourhood in a complex space or champ de
Deligne-Mumford analytique in which the Segre
class of an invariant measure around the hypothesised
compact singular locus in the ambient space is
zero, then modulo appropriate blowing up, an
invariant measure for the foliation $U\ra [U/\cF ]$
of implied type found in II.1.1 is zero (remember,
everything in the exceptional divisor is ignored).

One should observe that in terms of our goal of
relating residues to Segre classes, the hypothesis
of zero Segre class is unavoidable, {\it i.e.}

\noindent{\bf II.1.6 Remark} Already in dimension 3,
even supposing the centre manifold perfectly analytic one
might encounter a residue of the form,
$$\lim_{\e\ra 0} \int_{|xy|=\e} (1+\nu x^p) \dfrac{dx}{x^{p+1}y^q}\dm$$
for say $\nu\in\bc^{\times}$, which in the above
notation would become,
$$\lim_{\e\ra 0} \nu\int_{T,|y|=\e} y^{-q}\dm(\bft)$$
which admits no a priori reduction in terms of the
Segre class about $y=0$ or $x=0$.

\subsection*{II.2 Linear and unbounded}

Amongst all singularities in the centre
manifold those that are linearisable are
rather generic. Nevertheless, the behaviour
in the ambient 3-fold is dependent on the
eigenvalue of the linearisation. A case
where this behaviour is as simple as
possible is given by,

\noindent{\bf II.2.1(a) Set Up} In a
neighbourhood $U$ of 
foliated 3-fold of a point in the
singular locus after completion in the same,
the foliation admits a formal generator of
the form,
$$z\dz + \dfrac{x^{p}}{p(1+\nu x^{p})}\left(x\dx+\lb y\dy\right)$$
where $x,y,z$ are formal coordinates,
$p \in\bn$, $\lb,\nu\in\bc$, $\lb\neq 0$,
$\Re(\lb)\geq 0$. Again, after blowing up
if necessary, we may suppose that $x=0$
defines the exceptional divisor $E$.

Here by \cite{mp1} VI.4.3 the relation
between the formal and analytic structure
is straightforward, {\it viz:} 
there is an open neighbourhood
of the form, $V=S\times \D^2$,
with $x$ in a sector $S$ of aperture
up to $2\pi/p$ and $x\mpo \xi=x^{-p}$
branching within $\pi/2$ of the negative
real axis, while $y,z$ vary in discs, such
that over $V$ one may achieve an analytic
conjugation to the normal form. Again,
after a conformal change
of variable in $S$ we may suppose $\nu=0$, and
so obtain invariant functions,
$$s=z\exp(\xi),\,\, t=y\xi^{\lb/p}$$
As before, for fixed $\z$ we introduce
the transversal
$Z_\z:z=\z$, and its complement.
$Z^*_\z$ by the exceptional divisor,
which we fibre by way of $t$ over $T$
to obtain,
$$
\begin{CD}
Z^*_\z:=Z_\z\bsh E @>{\xi\times t}>> \bc\times T\\
@VV{t}V @.\\
T
\end{CD}
$$
where, again, all maps are open, and the horizontal
one an embedding. Once again 
the groupoid,
$$R\rras Z^*_\z$$
induced by the action of the foliation in $V$ on $Z^*_\z$
is that obtained by pull-back of the
groupoid in which $\tm$ acts on $\bc\times T$
by translations on the first factor.

Now we can argue as per II.1.4. Every leaf outside
the centre manifold and the exceptional divisor
meets $Z^*_\z$. The structure of the fibres is
straightforward, \cite{mp1} \S III.1, and every 
point in every fibre over $t$ not only has an
infinite orbit under $\tm$ but admits a neighbourhood
all of whose orbits are disjoint. Whence by II.1.2
we obtain,

\noindent{\bf II.2.2 Fact} Let the set up be as per II.2.1.(a),
then on the domain $V=S\times\D^2$,
any invariant measure for the foliation on
$V$ lies either in the exceptional divisor, $x=0$, or the
centre manifold, $z=0$.

Our next task is to consider the structure inside
the centre manifold. This is trickier than one
might think since apart from the fact that we
can only conjugate to the normal form for $x$ in
a sector, the invariant manifold $y=0$ need not
exist in all of $U$ but only in a sector of width
up to $3\pi/p$, \cite{mp1} \S V.2 \& \S VI.4. Nevertheless,
as might be expected, we assert:

\noindent{\bf II.2.3 Claim} If $\lb\notin\br_{+}$ then in
addition the support of any invariant $\dm$ off the
exceptional divisor is the curve $z=y=0$, whence
the latter is convergent in all of $U$ by \cite{siu}.

\noindent{\bf proof} Cover $U$ (more correctly $U\bsh E$)
by domains $V_\a=S_\a\times\D_\a^2$ where we can 
conjugate to normal form, and denote by $s_\a, t_\a$
the invariant functions previously constructed.
Plainly for any leaf $\ell$ in the centre manifold
$W_\a$ of $V_\a$, $s_\a(\ell)=0$.

For convenience choose the $V_\a$ such that $\a$
goes from 1 to $n+1$, $V_\a$ meets only $V_{\a-1}$
and $V_{\a +1}$, the overlaps are simply connected,
and we identify 1 with $n+1$. A priori on $V_\a\cap V_{\a+1}$,
$t_\a$ is a function of $t_{\a+1}$ and $s_{\a+1}$.
If, however, we consider $t_\a$ restricted to the
centre manifold $W_{\a+1}$ then it is a function
of $t_\a$ alone. Better still if $\Re(\lb)>0$,
then $t_{\a+1}$ maps $W_{\a+1}\bsh E$ onto $\bc$,
and every leaf in $W_{\a+1}$ meets $V_\a$, so
$t_\a|_{W_{\a+1}}$ is an entire function 
$h_{\a,\a+1}$ of $t_{\a+1}$. The evident
composition $h_{1,2}h_{2,3}\hdots h_{n,n+1}$
leads to an entire function $h$ of $t_1$ which
generates (a possibly complicated, and, a priori
no better than flat) groupoid,
$$H\rras\bc$$
The fact that the $W_\a$ need not patch implies that
this groupoid may have little to do with the 
action of the foliation in $U$. However, by II.2.2,
we may not only descend a hypothesised invariant
measure to $\bc$ by way of $t_1$, but it must also
be left invariant by $H$.

Now by \cite{mp1}, one may suppose that the formal
and analytic conjugations coincide modulo arbitrary,
but fixed, powers of the exceptional divisor, {\it i.e.}
not quite a full asymptotic expansion since
op. cit. is a bit lazy on this score. 
In terms of the $t_\a$ this translates into
asymptotics at $\infty$, so, in the first place
each $h_{\a,\a+1}$ sends $\infty$ to itself,
so, in fact, all of these, and whence $h$
is a polynomial. Furthermore $h$ is \'etale
at $\infty$, whence it's the first term in
the asymptotics, {\it i.e.}
$H$ is the action of
$\tm$ by $\exp(\lb\bullet)$. Finally, we
can cover by transversals $T_i$ in $V_1$
of finite mass, whose image under $t_1$
is surjective, whence the hypothesised
invariant measure is one of finite mass on $\bc$
invariant by the given action of $\tm$, which,
except for the $\d$-function at the origin, is
nonsense.

It remains to treat the case of $\lb\in\br(1)$,
which is a little easier since there is
actual holonomy about $x=0$ in the sectors
themselves.
As a function of the orientation, {\it i.e.}
$\Im(\lb) >0$, or $\Im(\lb) < 0$, 
its either contracting or expanding,
so, modulo a change of
orientation, without loss of generality contracting.
Consequently, for a sector $S\ni x$ transversal
to the said curve, the measure is invariant
by a groupoid,
$$
H  \, \textstyle\build\genfrac{}{}{0pt}{}{\longra}{\longra}_{\s}^{\tau} \, S
$$
generated by a map $h:S\ra S$ with Taylor series,
$$h(x)=e^{-2\pi|\lb|} x(1 + O(x^N)$$
for $N\in\bn$, which we could take
to be infinite, but we'll suppose
$N$ finite to indicate the robustness
of such a situation.
Indeed for $\e$
small we find a co-equaliser ({\it i.e.}
$\s_*-\tau_*$ of compactly supported) of 
the form,
$$\un_{|x|>q\e}-\un_{|x|>\e}$$
for $q$ as close to $e^{-2\pi|\lb|}$
as we please- {\it i.e.} the closer
the smaller the disc. Whence there
is no measure in the annular sector $q\e> x > \e$,
$x\in S$
for $\e$ small, but, otherwise, arbitrary,
so, the transversal has no measure off
the exceptional divisor. In this case 
every $T_\a$
is a disc, which, unlike the previous
case may be shrunk, so, without loss
of generality such transverse sectors
cover each $T_\a$ minus the origin.

Of course, in either case, the $\d$-function
at the origin, equivalently the curve $y_\a=0$
is only defined in $V_\a$, so we use \cite{siu}
to conclude. $\Box$ 

In the case that $\lb\in\br_+\bsh\bq_+$, a
so called reduced singularity, the above
discussion implies that any invariant measure
is a radially symmetric invariant measure 
on $\bc$ pulled back by $t_1$. Nevertheless,
such a measure could still be inconvenient,
so we again appeal to the Segre class, {\it viz:}

\noindent{\bf II.2.4 Fact} Suppose the $U$ 
of II.2.1.(a) is actually a
neighbourhood in a complex space or champ de
Deligne-Mumford analytique in which the Segre
class of an invariant measure around the hypothesised
compact singular locus in the ambient space is
zero, then modulo appropriate blowing up, an
invariant measure for the foliation $U\ra [U/\cF ]$
of the implied type  is zero.

\noindent{\bf proof} For $\lb\notin\br_+$ this 
follows a fortiori from II.2.3. Otherwise,
with notations as per II.1.5,
$$\mathrm{s}^U_{E,\dm}\geq \lim_{\e\ra 0} \int_{t\in T} \dm(t)  \oint_{V_t, |f|=\e} \frac{df}{f}$$
The condition $|f|=\e$ is as near $|\xi|=\e^{-p}$
as makes no difference for a sufficiently
high a priori approximation to the normal
form, so the integrand in $t$ (really $t_1$)
is an increasing function of $\e$ which 
gives the aperture of $S$ over $2\pi$ if,
$\log|t| < {\mathrm{const}.} -\log\e$, and
close to zero otherwise. As such, the
dominated convergence theorem applies,
so (shrinking $T$ a little if $\lb\in\br(1)$),
$$\mathrm{s}^U_{E,\dm}\geq \dfrac{|S|}{2\pi}\int_T\dm(t)$$
where, as before, we've profited from the
support being in the centre manifold to
descend the measure $\dm$ to $T$. $\Box$

We require to extend this rather satisfactory discussion to

\noindent{\bf II.2.1.(b) Set Up}
Everything as per II.2.1.(a) but now
the normal form of a
generator is,
$$(p+q\lb) z\dz + \dfrac{x^{p}}{1+\nu x^{p}}\left(x\dx+\lb y\dy\right)$$
where $x,y,z$ are formal coordinates,
$p,q \in\bn$, $\lb,\nu\in\bc$, $\lb, p+q\lb\neq 0$,
$\Re(p+q\lb), \Re(p/\lb +q) \geq 0$,
with suitable a priori blowing up
so that $x=0$, $y=0$ are exceptional
divisors, $E_1$, $E_2$, and $E=E_1+E_2$.

Here the best existence domains are 
obtained by passing to logarithmic
coordinates, or, better still,
viewing $\exp\times\exp:\gh\times\gh\ra\D^2$
as an \'etale neighbourhood of $\D^2\bsh E$,
for $\D^2\ni (x,y)$. In particular, the
most symmetric analogue of the previous
variable $t$ is to introduce,
$$\tau:=\dfrac{1}{p+q\lb}\log y -\dfrac{\lb}{p+q\lb}\log x$$
then fibre $\gh\times\gh$ over $T$ by $\tau$,
which can subsequently be embedded in
$\bc\times T$ by way of,
$$
\begin{CD}
\gh\times\gh @>>{(p\log x + q\log y)\times \tau}> 
\bc\times T \\
@V{\tau}VV @. \\
T
\end{CD}
$$
at which point one may describe the
existence domain $V$ for a conjugation 
to the normal form as the pull-back 
to $\gh\times\gh$ of a domain in which
$\log\xi:=-(p\log x + q\log y)$ varies
in a strip of width up to $2\pi$
resulting in a branch in the $\xi$
plane within $\pi/2$ of the real 
axis, and, of course, multiplied by
a disc in the variable $z$, cf. \cite{mp1} \S III.2.

Consequently, again after a conformal mapping
in $\xi$, we now have invariant functions,
$$s=z\exp(\tau),\,\, \tau$$
As such, we again look to the transversal,
$Z^*_\z:=Z_\z\bsh E$, factor it as, 
$$
\begin{CD}
Z^*_\z:=Z_\z\bsh E @>{\xi\times \tau}>> \bc\times T\\
@VV{t}V @.\\
T
\end{CD}
$$
and observe, once more, that every leaf outside
the centre manifold and the exceptional divisor
meets $Z^*_\z $, so the fibres $\dm_{\tau}$ of
a hypothesised invariant measure are, again,
invariant by the groupoid induced from translations
by $\tm$ in the $\xi$-variable.

Whence, we are again reduced to studying the 
possible structure inside the centre manifold.
Observe that the hypothesis exclude the
possibility that $\lb\in\br$, and the analogue
of $\lb\in\br(1)$ encountered in II.2.1.(a) is
$p+q\lb\in\br(1)$ or $p/\lb +q \in\br(1)$,
albeit at most one of these can occur. Now we
argue as in the proof of II.2.3, descending the
measure to a measure $\dm(\tau)$ which is
invariant by a groupoid,
$$H\rras T$$
generated by two  mappings, $A$, $B$, - {\it i.e.}
as per op. cit. but coverings indexed for
sectors in $x$ and $y$- entire  if neither
$p+q\lb\notin\br(1)$ nor $p/\lb +q \notin\br(1)$,
while the domain of $\tau$ and the said
mappings are right, respectively left,
half planes otherwise.
In either case, as before, we know by
\cite{mp1} the form of the mappings
modulo some large, but fixed, power
of the exceptional divisor, {\it i.e.}
for $a,b$ generators of $\tm^2$ they
have the form,
\begin{align*}
A:&T\longra T  :\tau\mpo \tau + \dfrac{a}{p+q\lb} + O((x^py^q)^N) \\
B: &T\longra T :\tau\mpo \tau + \dfrac{b\lb}{p+q\lb} + O((x^py^q)^N)
\end{align*}
Unfortunately this error doesn't quite
translate into an expansion at $\infty$
in the $\tau$ variable. As such, put,
$Y=(p+q\lb)\tau$, then for $\Re(Y)\ra\-\infty$
in strips with $\Im(Y)$ bounded, the
induced groupoid has asymptotically the
same generators as that induced by
translations in $\tm$ and $\tm\lb$.
Whence taking the strip to have width
in the imaginary direction at least $2\pi$,
and $\Re(Y)< -R$ for some sufficiently
large (determined by $\Im(\lb)\neq 0$)
$R$ each point in such a strip has an
open neighbourhood with an orbit enjoying
infinitely many connected components,
each with mass bounded below if the
measure were supported at the point in
question. The strip has, however, finite
mass, so this is nonsense. Arguing
similarly in the $x$-variable, and possibly
shrinking the initial neighbourhood $U$
if necessary,
we cover all of $T$ by points with orbits
meeting such strips, and deduce:

\noindent{\bf II.2.5 Fact} Suppose we are in
the situation II.2.1.(b) then, without any
hypothesis on the Segre class, any invariant
measure must be supported in the exceptional
divisor.

\subsection*{II.3 Exceptional nodes}

A saddle node singularity inside the
exceptional divisor is an obvious
source of concern, and, again, like
linearisable singularities the
relation between formal and
analytic theory is not always
straightforward. A particular
case where it is not too complicated
is when the centre manifold of the
node is an exceptional divisor,{\it i.e.}

\noindent{\bf II.3.1 Set Up} In a
neighbourhood $U$ of 
foliated 3-fold of a point in the
singular locus after 
appropriate blowing up, we have a
singularity such that on
completion in the point (not,
as before, completion in the singular locus)
we find a formal generator of
the form,
$$qz\dz + \dfrac{x^{p}y^q}{1 + y^q(R(x)+\lb x^{p+r})}\left(y\dy+
\dfrac{x^r}{1+\nu x^r}(qx\dx-py\dy)\right)$$
where $x,y,z$ are formal coordinates,
$p,q,r \in\bn$, $\lb,\nu\in\bc$, 
deg$R\leq r-1$, and the exceptional
divisor $E$ has 2 components $E_1:x=0$,
$E_2:y=0$, and we may even assume that
the normal form is defined after 
completion in the former, albeit 
not the latter.

The analytic theory has a number
of complications. In the first place
it is convenient to have all of $\lb,\nu,R$
equal to zero, a so called {\it monomialised}
form. This may be achieved,
\cite{mp1}, VI.2.2.(g),
once and for all, {\it i.e.} 
depending only on the said parameters,
in a sub-domain of $S\times\Sigma\times\D$ where
$x$ varies in a sector of aperture
up to $2\pi/r$, $z$ in the disc $\D$,
and $y$, or better $\log y$, in a
so called spiralling domain: $\log y$
in a sector of width up to (which
means strictly less than) $\pi$ 
wholly contained in the left half plane.
Having effected this change of variables
on some coordinates affording a sufficiently
large approximation to the normal form,
we restrict our attention to the domain
away from the negative real-axis in
$w=x^{-r}$, {\it i.e.}
 $|\Im (w)| > R$, or $\Re(w)> R$, $R$ large-
otherwise leaves may be bounded,
and invariant measures may exists.
Now the further restrictions we
must make to achieve a conjugation
to the monomialised form are,
for $\xi$, 
again equal to $(x^py^q)^{-1}$, we
must,
necessarily, restrict it
to a sector $S_\xi$ (or, better $\log\xi$
is restricted to a strip) of apertures
up to $2\pi$ branched
within $\pi/2$ of the negative real axis,
in addition the parameter $w$ can vary
through at most $3\pi/2$, or, equivalently,
away from the negative real axis in
the same, from above it we can go to
within the lower imaginary axis, while
from below to within the upper imaginary
axis. This yields, \cite{mp1} VI.4.5, a domain 
$V$ such
that after the formal to analytic conjugation
of the variables,
$$s=z\exp(\xi),\,\,\,\, \tau=w-\log\xi$$
are invariant functions. Now we can argue
exactly as before: restrict to a
transversal, fibre in $\tau$
observe, cf. \cite{mp1} \S IV.4, that
the condition of being away from the
negative real axis in $w$ is exactly
what guarantees that the fibres in $\tau$
embedded by $\xi$ in $\bc$ contain,
apart from the branch, full neighbourhoods
of $\infty$, and the discrete action
by translation of $\tm$ on the same
leaves the measure invariant, apply
II.1.2, and conclude,

\noindent{\bf II.3.2 Fact} Suppose we are in
a domain $V$ as above, of which the
critical feature is being bounded away
from the negative real axis in $x^{-r}$,
then the support of any invariant measure
is contained in the centre manifold
or the exceptional divisor.

Now we must aim to study the situation
inside the centre manifold, which risks
being a little complicated since we have
as many of these as we have determinations
of $\log\xi$. As such fix analytic 
functions $x,y$ in $U$ defining the
exceptional divisor, and for $\log(x^py^q)^{-1}$
in a strip $\a$ with branching as described
above, denote by $V_\a$ the corresponding
existence domain where the conjugation
to monomial/normal form is valid, with
$Z_\a$ its centre manifold and $\tau_\a$
as above. The domain $T_\a$ of $\tau_\a$
is all of $\bc$, and as before if $\a\cap\b\neq 
\emptyset$ then every leaf in $Z_\a$
meets $V_\b$, and conversely. Whence
we may write,
$$\tau_\a|_{Z_\b}=h_{\a\b}(\tau_\b)$$
for $h_{\a\b}$ entire. Identifying this
function is easy since an a priori
choice of approximation to the normal
form translates into a bound at $\infty$
with knowledge of a finite number of
terms in the Taylor expansion. Whence,
say for convenience, all strips of the
same width $|S|$, then there are $\lb_\a\in\br(1)$
such that,
$$h_{\a\b}=\tau_\b + (\lb_\b-\lb_\a) + O(R^{-N}),$$
$N$ as large as we like, so $h_{\a\b}$ is
translation by $\lb_\b-\lb_\a$. Observe that
the intersection $W_\a$ of $V_\a$ with
$Z_\a$ is embedded via:
$$(w,  \log(x^py^q)^{-1}): W_\a\hookrightarrow
S\times L$$
where $S\sbs\bc$ is defined by way
of our sectorial and modular restrictions in $w$,
with $L$ the surface of the logarithm,
and the above procedure allows us to
extend say $\tau_0$ from $W_0$, $0$
a strip, to $\tau$ on all of $S\times\bc$.
Again this map would have little or
nothing to do with the ambient foliated
3-fold were it not for the fact that
the hypothesised invariant measure
$\dm$ is necessarily supported in $Z_\a$,
while the existence domain $V_\a$ is 
naturally a sub-domain of $\D\times S\times L$,
and the pull-back of the measure to the
latter is invariant by translations 
under $\tm$ on the last factor. As
such, pulling back $\tau$ by the
projection $\D\times S\times L\ra S\times L$,
it follows that $\dm$ not only descends
to a measure $\dm(\tau)$ on $T_0$,
but it is also invariant by translations
in $\tm$. Finally, the image under $\tau$
of a transversal, say $S\times 0$ has
finite mass, and each point an open
neighbourhood with an infinite disjoint
orbit, so the said image has zero mass, while
every point in $T_0$ has a point of this
image in its orbit, so, indeed the measure
on $T_0$ is zero, and we deduce:

\noindent{\bf II.3.3 Fact} Suppose we
are in the set up II.3.1, and, in a
minor notational confusion, suppose
$x$ is an analytic function in $U$
defining the exceptional divisor $E_1$,
then for $R$ sufficiently large, outwith
the domain: $|\Im (w)| <R, \Re(w)< -R$,
any invariant measure is supported in
the total exceptional divisor $E_1+E_2$.

While similar to the argument of II.2.4/5
for excluding support in the centre manifold,
the above is more complicated again, so
let us make,

\noindent{\bf II.3.4 Remark} Say for simplicity,
even if it can be excluded by blowing up,
$p=0$, and the centre manifold converges in
$U$, then normally one excludes an invariant
measure outside $x=0$, $y=0$, the so called
strong and weak branches, by examining the
holonomy around the strong branch which in
an appropriate Fatou coordinate is translation
by $\tm$. On the other hand, before even
beginning we made a spiralling restriction in
$\log y$, so it might seem that there should
be insufficient holonomy to conclude. However,
this isn't how it works. More precisely
under the above
sectorial restrictions, say $x^{-r}\in S$ one finds
a first integral $t$, of which $\tau$ is
the logarithm, and the implied \'etale 
covering of $S\times\D^{\times}$ is equally
that given by $\log y$. On this covering,
$S\times L$,
irrespective of spiralling restrictions, no
leaf has a return map on any $S\times\ell$.
On the other hand the restriction of $\tau$
to any such is Schlict, and we have a natural
action of $\tm$ on $\tau$ which leaves the
descended measure invariant. All of which
remains true even under the spiralling 
restriction on $\log y$, or, if one prefers
the image in $S\times\D^{\times}$ of a
$\tm$ orbit of leaves in $S\times L$ is
pretty much the same leaf as before, {\it
i.e} contains the full leaf in $S\times\D'$,
for $\D'$ a slightly smaller disc.

\newpage

\section*{III. Almost holonomy}

\subsection*{III.1 Local Set Up}

As before $U$ will be a small open
neighbourhood, say a polydisc, in
a 3-fold foliated by curves, about
a singular point, say the origin,
and the singular locus will be
denoted $Y$. In particular, modulo
blowing up, $Y$ is either a smooth
curve, or a plane curve with a node,
and in the completion $\wh{U}$ of
$U$ in $Y$ there is a well defined
formal centre manifold $\wh{Z}$
provided the singularity is not
a beast in the sense of \cite{uu} \S I.5,
\cite{mp1} II.2, so for our present
purposes, essentially for notational
convenience, we'll omit this case.
As such, by way of appropriate
blowing up we may further suppose that
there is an invariant simple normal
crossing divisor $E$, all of whose
components are smooth and in 1-1
correspondence with components of Y
by way of,
$$E_i\in |E| \longra E_i\cap\wh{Z}\sbs Y$$
whence there will always be at least
one component of $E$ meeting $U$, and
a local equation for it will always
be given by the variable $x$ equal 
to zero. As per the convention already
implicit in \S II if there is a further
component it will be given by the
variable $y$ equal to zero, albeit
the latter may, otherwise, just be
another variable. This convention
is followed throughout \cite{uu}
and \cite{mp1}, and the above reduction is
explained in \cite{uu} \S I.4, or
\cite{mp1} \S V.1. In any case there
is a well defined (formal) foliation
in the formal centre manifold induced
from $U$ and we'll suppose further
that our point of interest is a
singular point of this formal
foliation with $Y$ being invariant
by the same. Consequently for any
simple closed curve $\g$ in $Y$ there
is well defined formal holonomy
of the said induced foliation in
the centre manifold, or, equivalently
an inverse system of representations,
$$\rho_n :\pi_1 (\g)\longra\Aut\left(\bc[t]/t^n)\right)$$
for $n\in\bn$. Plainly this 
representation can be read off
from the normal form. A complete
list is provided in \cite{mp1} VI.1.1,
and we always have the right to choose a generator
of the foliation (excluding beasts,
and the case already encountered
in II.3.1) without any loss
of domain in the form,
$$\pa=z\dz + x^py^q u(x,y) D\,\,\mod\,\, \cO(-NE)$$
for $u$ a unit, $p,q\in\bn\cup\{0\}$, $p\neq 0$,
$N\in\bn$ as large as we like, $x,y,z$ honest
coordinates, and $D$ some 2-dimensional normal
form.

In so much as planar normal forms consist
of partial derivatives and at worst
fractions in polynomials, we can define
real co-dimension 1 subsets $\G_\e$ as
follows:

\noindent{\bf III.1.1 Basic Construction}
Obviously a simple closed curve $\g$
in $Y$ lies in a unique irreducible
component of the same, and we can solve
the equation $Dt=0$ around $\g$ in the normal
direction to the component- {\it i.e.} $t$
is a unit times $x$ or $y$ according to
the case. There will, however, most likely
be a discontinuity arising from the
holonomy of the approximating plane
foliation generated by $D$ as we turn
through $2\pi$. Consequently for $\e>0$
we have a real co-dimension 1 variety
$\G_\e$ (strictly $\G_{\g,\e}$ since
often estimates will have dependence
on $\g$) defined by $y$ (respectively
$x$) in $\g$, $|z|\leq \e$, $|t|\leq\e$,
with a possible discontinuity in the $t$
variable at $2\pi$- alternatively, if
one prefers things smoother, say tubes
$\G_{\e,s}$ around intervals $I_s$ in
$\g$ with $I_s\ra\g$. Now suppose that
we have an invariant  measure $\dm$ for
the foliation in $U$, then modulo the
precision (which will frequently be
omitted) of excluding sets of zero
Lebesgue measure be it in $\e$ or
perturbations of $\g$, $\dm|_{\G_\e}$
is well defined, and for any function
$f(t)$ we can apply Stokes to obtain:
\begin{align*}
& \int_{\G_\e}  df(t)\dm =\\
& \int_{|z|=\e,|t|\leq\e}  f(t)\dm
  + \int_{|t|=\e,|z|\leq\e} f(t)\dm 
+ \int_{\substack{|t|\leq \e,|z|\leq\e\\ \bfp }} \left(f(t^h)-f(t)\right)\dm
\end{align*}
The last integral being at the
discontinuity, and $t\mpo t^h$ being
the (convergent) holonomy of the
approximating plane foliation defined
by $D$ around $\g$. In particular one recognises the final
integrand
as the co-equaliser of an approximation
to the formal holonomy groupoid in
the centre manifold. It is, however,
an approximation to an object whose
action on $\dm$ has no sense. Nevertheless,
amongst the above terms, the integral
over $\G_\e$, and that over $|t|=\e$
can invariably be made $O(\e^N)$, for
any $N\in\bn$. For example, the choice
of $f$ the characteristic function $\un_F$
of some set is usually best, so the
formula becomes,
$$
\int_{\substack{|t|\leq \e,|z|\leq\e\\ \bfp }} \left(\un_F(t)-\un_F(t^h)\right)\dm=
\int_{t\in\pa F, |z|\leq\e}\dm + \int_{|z|=\e,t\in F} \dm$$
There will always be a plane meromorphic
1-form $\o$ with poles on $E$ such
that $\pa(\o)=1$, and $\g$ will never
be allowed to shrink arbitrarily. As
such, only the pole around the irreducible
component containing $\g$ is ever relevant,
and in a way that depends on $\g$ this
is bounded by some a priori fixed power
$|t|^{-p}$. On the other hand $\pa t$
vanishes to some high order, say, in
a minor abuse of notation, $N+p$, so
$|dt|$ wedged against $\dm$ will be
no worse than that of $|\o|$ times
$|t|^N$. Presently we'll work through
some examples in detail, but, already
we can reasonably say that the dominating
term on the right of this formula is
that on the face $|z|=\e$. This will,
however,
prove ($s$, Lebesgue a.e allowing $F$ to move in some
family $F_s$) to be bounded by $\e^p$ times the part, $ \mathrm{s}_{Z,\dm}(\e)$,
of the Segre class around the singular locus
in the normal direction to the centre
manifold, so we'll have,

\noindent{\bf III.1.2 Almost holonomy estimate} Let things
be indicatively as above with $p$ the order of
vanishing of $\pa$ restricted to the formal centre
manifold around the component of $Y$ containing
$\g$, then for a sufficiently good approximation
indexed by $N\in\bn$ to the formal foliation
in the formal centre manifold,
$$
|\int_{\substack{|t|\leq \e,|z|\leq\e\\ \bfp} } \left(\un_F(t)-\un_F(t^h)\right)\dm|
\leq \e^p \mathrm{s}_{Z,\dm}(\e) + O(\e^N)
$$

An evident variant, indeed corollary, of
the above reasoning occurs on varying
the branch point $\bfp$. Typically, in
practice for, say, loops $\g$ of the
form $|y|=r$, $r\in I$, and the branch
the interval $I$ itself viewed in $\br_+$
inside the domain of $y$ in $\bc$, or,
for that matter $\br_+\lb$ for any direction
$\lb$. In such a situation we have,

\noindent{\bf III.1.3 Variant} Notations as above, then:
\begin{align*}
& \int_{|t|\leq \e,|z|\leq\e, \bfp\in b } \dfrac{d|y|}{|y|}\left(\un_F(t)-\un_F(t^h)\right)\dm =\\
& \int_{t\in\pa F, |z|\leq \e, |y|\in I}\ \dfrac{d|y|}{|y|}\dm +
\int_{|z|=\e,t\in F, |y|\in I}  \dfrac{d|y|}{|y|}\dm
\end{align*}
which in turn could have been viewed at
Stokes applied to the region, $t\in F$,
$|y|\in I$, $|z|\leq \e$ since $d|y|$ 
vanishes on the boundary of $I$. In any
case, the important point is that 
$d|y|$ 
has a sign on the branch, so, without
loss of generality we may suppose that
it is an actual length form.

\subsection*{III.2 Linear and bounded}

We proceed to use the basic construction
III.1.1 to get estimates at the singularities
not encountered in \S II beginning with,

\noindent{\bf III.2.1(a) Set Up} In a
neighbourhood $U$ of 
foliated 3-fold of a point in the
singular locus after completion in the same,
the foliation admits a formal generator of
the form,
$$z\dz + \dfrac{x^{p}}{p(1+\nu x^{p})}\left(x\dx+\lb y\dy\right)$$
where $x,y,z$ are formal coordinates,
$p \in\bn$, $\lb,\nu\in\bc$, $\lb\notin \br$,
$\Re(\lb)< 0$, and as per III.1, after blowing up
if necessary, we may suppose that $x=0$
defines the exceptional divisor $E$.

Plainly we take an approximation modulo
a large, to be decided, power $\cO_U(-NE)$ 
of the exceptional divisor as discussed
pre III.1.1, so that, in a minor notational
confusion, $x,y,z$ are honest holomorphic
coordinates. The holonomy of the approximating
plane foliation is, therefore, always determined
by its first order part, and given by,
$$\pi_1(\g)\iso\tm\longra\bc^{\times}: a\mpo \exp(-a/\lb)$$
Even though the above map is wholly independent
of a choice of square root of $-1$, there
is an orientation convention in Stokes'
theorem, so as we follow an oriented loop
there is an implied choice of such, and
a contraction if $\Im(\lb^{-1})> 0$, expansion
otherwise, so either changing the orientation,
or replacing $\e$ by $\e\exp (-2\pi\Im(\lb^{-1}))$
if $\Im(\lb^{-1}) < 0$, there is $q\in\br_+$ of
modulus less than $1$ such that,
$$\int_{\substack{\mzle, q\e'\leq\mt\leq\e'\\ \bfp}} \dm =
\int_{\substack{\mze,\mt\leq\e'\\ y\in\g}} \dm
+ \int_{\substack{\mt=\e',\mzle\\ y\in\g}}\dm
$$
which leads to the formula,
$$
\int_{\substack{\mzle, \mt\leq\e'\\ \bfp}} \dm =
\sum_{n=0}^{\infty} \left( \int_{\substack{\mz=q^n\e,\mt\leq q^n\e'\\ y\in\g}} \dm
+ \int_{\substack{\mt=q^n\e',\mz\leq q^n\e\\ y\in\g}}\dm\right)
$$
Now, say $\g$ is the loop, $\my=Y$ for $Y$ in
a range $I=(a,b)$, $a > 0$, then for constants $c_I, C_I$
depending on $I$,
$$c_I\mx < \mt < C_I\mx$$
so for possibly different constants $c,C$ again
depending on $I$ we get an estimate,
$$
\int_{\substack{\mz,\mx\leq\e\\ \bfp}}\dm\leq
\sum_{n=0}^{\infty} \left( \int_{\substack{\mz=q^n\e,\mx\leq cq^n\e\\ y\in\g}}\mmu 
+ \int_{\substack{\mt=Cq^n\e',\mz\leq q^n\e\\ y\in\g}}\mmu\right)$$
where here, and elsewhere. $\mmu$ is the total variation 
of the sliced normal distribution with measure regularity,
cf. \cite{feder} 4.3.2, whereas on the left $\dm$ is
naturally a measure, and indeed the left hand side is
the measure of a transversal. To proceed further
let us suppose that the Segre class is well defined
around $Y$, {\it e.g.} $U$ a neighbourhood in a complex
space or champ de Deligne-Mumford analytique with
compact singular locus. Every choice of distance
function to $Y$ gives rise to different measures
depending on $\e$ that limit on $\mathrm{s}_{Y,\dm}$, but
the difference is $o(\e)$, so let us omit this
from the notation, and write:
\begin{align*}
\mathrm{s}_{Y,\dm}(\d) & \geq \oint_{\substack{\mzd,\mx\leq c\d\\ \my\leq A}} \dfrac{dz}{z}
+ \oint_{\substack{\mx= c\d,\mzld\\ \my\leq A}} \dfrac{dz}{z} \\
& \geq \oint_{\substack{\mzd,\mx\leq c\d\\ \my\leq A}} \dfrac{dz}{z}:= \mathrm{s}_{Z,\dm}(\d)
\end{align*}
where $\d>0$ is anything less than some fixed constant
outside a set of nil Lebesgue measure, $A$ is just
the upper limit of the size of our disc in the
variable $y$, and $c$ is a constant to be chosen, 
{\it i.e.} here our distance function to $Y$ is
$\max\{\mz,c^{-1}\mx\}$ in the neighbourhood $U$.
Consequently,
\begin{align*}
& \int_a^b dr\int_{\substack{\mzd, \mx\leq c\d\\ \my=r}}\dm =\\
& 2\pi\oint_{\substack{\mzd, \mx\leq c\d\\ a\leq \my\leq b}}
\Bigl\lvert \Im\{\dfrac{y\bar{z}\pab \bar{ y}}{\my\pab\bar{z}}\}\Bigr\rvert\dfrac{dz}{z}\dm
\leq C \d^p\mathrm{s}_{Z,\dm}(\d)
\end{align*}
where the implied constant depends on the
interval $(a,b)$, and, of course, the level
$N$ of approximation that we're working to.
As such one could say more precisely that $C$ is of the 
form $C_1 c^p$ for $C_1$ depending only on $A$,
and $\d\leq C_2c^{-1}$, again $C_2$ depending
only on $A$, and a further dependence on the
interval built into the choice of distance
function. In any case, this yields an estimate,
$$\int_a^b dr\sum_{n=0}^{\infty}\int_{\substack{\mzd, \mx\leq c\d\\ \my=r}}\dm 
\leq \dfrac{C\e^p}{1-q^p} \sup_{\d\leq\e} \mathrm{s}_{Z,\dm}(\d)$$
Consider, therfore, the set of $r\in(a,b)$ where
the integrand on the left exceeds $B$ times the
estimate on the right. This has Lebesgue measure
at most $B^{-1}$, and so we deduce,

\noindent{\bf III.2.2 (a) Estimate} Let $\g$ of the 
basic construction III.1.1 be taken of the form
$\my=r$ for $r$ in the interval $I=(a,b)$, then
there is a constant $C$ and limiting Segre classes
$\mathrm{s}_{Z,\dm}(\d)$ depending on $I$ such that for
every $\e$, and $r$ in a set (depending on $\e$)
of measure at least $(b-a) - B^{-1}$ we have,
$$\sum_{n=0}^{\infty}\int_{\substack{\mzd, \mx\leq c\d\\ \my=r}}\dm 
\leq CB\e^p \sup_{\d\leq\e} \mathrm{s}_{Z,\dm}(\d)$$

Next for $\e<\d$ say, and $\my\in I$ as before
with $m\in\bn$ to be chosen, consider,
$$\int_0^\d \dfrac{d\e}{\e^{m+1}}
\int_a^b\sum_{n=0}^{\infty} \int_{\substack{\mt=Cq^n\e,\mz\leq q^n\e\\ \my=r}} \dm
\leq \sum_{n=0}^{\infty} \int_{\substack{\mt\leq Cq^n\d,\mz\leq \d\\ \my\in I}}
\Bigl\lvert\dfrac{dtd\my}{t^{m+1}}\Bigr\rvert\dm$$
As a function of the level $N$ of approximation
modulo powers of $(x^N)$ the exceptional divisor,
the integrand on the right is at most,
$$C\mx^{N-2(p+1)-m} d x\otimes d\bar{x}\dm$$
again for $C$ depending on the interval $(a,b)$.
As such, if $J(\e)$ is the (necessarily non-negative)  value of the
integrand on the left with respect to $\e$
it satisfies a bound,
$$\int_0^{\d} \dfrac{J(\e)}{\e^m}\dfrac{d\e}{\e} \leq \dfrac{C\d^n}{1-q^n}$$
provided $n=N-2(p+1)-m > 0$. Choosing $N$
appropriately large, implies,
that the set of $\e$ where $J(\e)> \e^m$ has 
finite measure with respect to $\e^{-1}d\e$
and we obtain,

\noindent{\bf III.2.2 (b) Estimate} Again let
$\g$ be of the form $\my=r$ for $r$ in $(a,b)=I$,
with $C$ (possibly bigger) and $\szmd$ as in
III.2.2 (a), together with an exceptional set
of $\e$ (again depending on $I$) of finite
$\lem$, then for non-exceptional $\e$ and $r$ in
a set (again depending on $\e$) of measure at
least $(b-a)-2B^{-1}$ the estimate III.2.2 (a)
holds along with the estimate,
$$\sum_{n=0}^\infty \int_{\substack{\mt=Cq^n\e, \mz\leq q^n\e\\\my=r}} \dm \leq CB\e^N$$
where, in what will be a repeated notational
confusion, $N\in\bn$ is just some a priori
choice of an integer as large as we please,
which may be the order of approximation along
$E$ or some finite shift of the same by
irrelevant and very fixed constants,
albeit that other implied constants, sets of
exceptional measure, {\it etc.}, may depend on $N$.

Putting all of this together implies that 
on the transversal $T_\bfp$ at $\bfp$ we have,

\noindent{\bf III.2.2 (c) Conclusion} Taking $N$
sufficiently large, with quantifies
as above in III.2.2 (b) then the measure
of the transversal 
$T_\bfp: \mxle, \mzle, y=\bfp$ satisfies,
$$\int_{T_\bfp} \dm\leq B\e^p\max\{\sup_{\d\leq\e}\szmd,\e\}$$

At which point, one should recognise that this
is a variant on the principles that yielded
the much stronger II.2.3. The results of \cite{mp1},
which are nevertheless optimal, cf. op. cit. \S III.3,
are nowhere close to permitting II.2.3 in this case.
One particular (and by no means the worst) problem
is the divisor $y=0$ which, formally, may be taken
invariant, cannot be supposed such, whence, in
our immediate context, the emphasis on the interval
$I$ which is bounded away from 0. Nevertheless the
formal invariant curve $y=z=0$ might converge, and
it could then support an invariant measure. On the
other hand, the above conclusion obviously cannot have
any relevance to such a measure, so, once again, the
Segre class must intervene. To this end, suppose the
Segre class around the exceptional divisor is zero,
then we have the identity,
\begin{align*}
\oint_{\substack{\mxe\\ \mzle\\ \my\leq r}} \dfrac{dx}{x}\dm
\, + \, \oint_{\substack{\mze\\ \mxle\\ \my\leq r}} \Re \left(\dfrac{z\pa x}{\pa z x}\right)\dfrac{dz}{z}\dm \,
+ \, \oint_{\substack{\my=r\\ \mzle\\ \mxle}} \Re \left(\dfrac{y\pa x}{\pa y x}\right)\dfrac{dy}{y}\dm
\,=\,0
\end{align*}
for any $r\leq A$, albeit outside a set of zero
Lebesgue measure. Now the integral over the face
$\mze$ plainly admits the estimate $\e^p\szme$.
The integral over the face $\my=r$ is just the
average of the previous discussion over $\bfp\in\g$
with respect to circle measure, and, at the price
of repeating everything mutatis mutandis with an
extra integral thrown in, all of III.2.2 (a),(b),(c)
remain valid with the same quantification when 
averaged over $\bfp\in\g$. Alternatively one can
homotope between different $\bfp$ (as we'll do later)
and use the estimates as they stand. In any case,
the integral over the $\mxe$ face is strictly
positive, and we have,

\noindent{\bf III.2.3 Corollary} Suppose further
that the Segre class around the exceptional divisor
is zero, then for a constant $C$, and an exceptional
set of finite $\lem$ measure both depending on the interval
$(a,b)$ as also the classes $\szme$ we have,
$$\oint_{\substack{\mxe, \mzle\\ \my\leq a}} \dfrac{dx}{x}\dm\leq \e^p\szme$$
Alternatively for $r\in (a,b)$ one has the same statement
with the more convoluted quantification encountered in
III.2.2.

The zero Segre class hypothesis, so $\szme\ra 0$ as
$\e\ra 0$ permits an inessential, but convenient 
improvement of the quantification encountered
in III.2.2 (a) to something more similar to III.2.2 (b).
Specifically, put,
$$\Theta_Z(\e):=\sup_{\d\leq\e} \szmd$$
The latter, understood globally for added convenience,
decreases to zero, so the measure,
$$d\theta_Z(\e) = \dfrac{\Theta'_Z}{\Theta_Z}d\e$$
has infinite mass, while,$\Theta_Z d\theta$ has finite
mass. Consequently we can on dividing by $\e^p$, take
a supremum over $\e$ of the integrand, then integrate
against $d\theta_Z (\e)$ to obtain,

\noindent{\bf III.2.2 (a)/(c) bis. Alternative}
Let $g, I$ be as before, then there is a function
$o(\e)$ going to zero a $\e\ra 0$ such that for
$r$ in a set (depending on $\e$) of measure at
least $(b-a) - B^{-1}$ the estimated quantity 
be it III.2.2 (a) or (c) is bounded by,
$$B\e^p o(\e)$$

This is a somewhat more convenient way to proceed
when addressing,

\noindent{\bf III.2.1 (b) Set Up} Everything as before
in III.2.1 (a), except that now $\lb\in\br$. Furthermore
if $\lb\in\bq$ we insist that its big height is
sufficiently large compared to $p$. About $(p+1)$
will do, but one can achieve anything by blowing up
which certainly improves the situation in the formal
centre manifold since at other singularities in its
proper transform one actually has 2 rather than
just 1 component of the exceptional divisor.

In this situation there is absolutely no obstruction
to the existence of a non-trivial invariant measure, 
and we proceed by computing meromorphic residues
around the exceptional divisor. Specifically for
$i> 0$, $r\leq A$, one has,
$$\Bigl\lvert \int_{\substack{\mxe,\mzle\\ \my\leq r}} dx^{-i}\dm\Bigr\rvert \leq
\e^{-i} \int_{\substack{\mxe, \mze\\ \my\leq A}} \mmu 
+ \e^{-i} \int_{\substack{\mxe, \mzle\\ \my=r}} \mmu $$
The estimates of the terms on the right is very
much similar to III.2.2 (a), respectively (b).
Indeed for $\chi$ in say $(e^{-1}, 1)$ consider,
$$
\int_{e^{-1}}^1 \dfrac{\chi}{\chi} \int_{\substack{\mx=\chi\e,\mze\\ \my\leq A}}\mmu
= \int_{\substack{e^{-1}\e\leq\mx\leq 1 \\\mze, \my\leq A}}
\Bigl\lvert \dfrac{dx}{x}\Bigr\rvert\dm \leq C\e^p \szme$$
for a constant depending only on $A$, so that
for say $J$ the range of $\e$,
$$\int_J d\theta_Z(\e) \sup_{\d\leq\e}\left(
\d^{-p} \int_{e^{-1}}^1 \dfrac{\chi}{\chi} 
\int_{\substack{\mx=\chi\d, \mzd\\ \my\leq A}} \mmu\right) \leq C$$
 
For a possibly different constant $C$, and we
put ourselves in the situation of zero Segre
class, equivalently $d\theta_Z(\e)$ has infinite
measure to deduce for an $o(\e)$ as per
III.2.2 (a)/(c) bis, but actually no dependence
on $I$ (III.1.1 giving no information here)
albeit certainly depending on $A$,
$$\int_{e^{-1}}^1 \dfrac{\chi}{\chi} 
\int_{\substack{\mx=\chi\d, \mzd\\ \my\leq A}} \mmu
\leq \e^po(\e)$$
From which we obtain,

\noindent{\bf III.2.4 (a) Estimate} For $o(\e)$
as above and every $B>0$ there is subset of $(-1,0)$
(depending on $\e$) of measure
at least $1-B^{-1}$ such that for $\log\chi$
belonging to the same we have the estimate,
$$\int_{\substack{\mx=\chi\d, \mzd\\ \my\leq A}} \mmu\leq B\e^po(\e)$$

The other term, like III.2.2 (b), is much more
robust. Indeed, introducing once more our
interval $I$, bounded away from 0, in which
$\my$ varies, for $m>0$ to be chosen, we have
for $\e\in (0,\d)$,
$$\int_0^\d\dfrac{d\e}{\e^{m+1}}\int_a^b dr
\int_{\substack{\mxe,\mzle\\ \my=r}}\mmu\leq
\int_{\substack{\mz,\mx\leq\d\\ \my\in I}} \Bigl\lvert \dfrac{d\mx d\my}{\mx^{m+1}}\Bigr\rvert\dm$$
Furthermore under the hypothesis III.2.1 (b)
we have for some large $N$, possibly not arbitrary
if $\lb\in\bq$, but the big height condition
says that its large enough,
$$\Bigl\lvert \dfrac{d\mx d\my}{\mx^{m+1}}\Bigr\rvert\dm
\leq C\mx^N dx\otimes d\bar{x}\dm$$
for $C$ depending only on $I$, and $\lb\in\br$
is used in an essential way. Integrating against
$\chi$ one therefore obtains,

\noindent{\bf III.2.4 (b) Estimate} For $\e$
outside a set with finite $\lem$ measure,
there is a constant $C$, depending on $I$,
such that for every $B>0$ and $(r,\log\chi)$
in $I\times (-1,0)$ outside a set (depending
on $\e$) of measure
$B^{-1}$ we have,
$$\int_{\substack{\mxe\\ \mzle,\my=r}}\mmu\leq CB\e^{p+1}$$
or indeed $\e^{N}$ rather than $\e^{p+1}$,
$N$ only limited by the big height proviso
in the case $\lb\in\bq$, but all constants,
exceptional sets, {\it etc.} depending on it.

As such combining III.2.4 (a), (b) we obtain,

\noindent{\bf III.2.4 (c)} Suppose things are
as per the set up III.2.1 (b), and that the
Segre class of our invariant measure around
the singularities is zero, then taking $B$
sufficiently large so that III.2.4 (a), (b)
hold for every $\e>0$, and $(r,\log\chi )$
varying in a set (depending on $\e$) which
is as close to full as we please,
$$\lim_{\e\ra 0} \int_{\substack{\mx=\e\chi^{-1}\\ \mzle,\my\leq r}}
\dfrac{dx}{x^j}\dm=0,\,\,\, 1\leq j\leq p+1$$

\subsection*{III.3 Nodes}

Next on our to do list is
to consider a node in the
centre manifold, {\it i.e.},

\noindent{\bf III.3.1(a) Set Up} In a
neighbourhood $U$ of 
foliated 3-fold of a point in the
singular locus after completion in the same,
the foliation admits a formal generator of
the form,
$$(p+r)z\dz + x^{p}\left(R(x) y\dy + \dfrac{x^{r+1}}{1+\lb x^{p+r}}\dx \right)$$
where $x,y,z$ are formal coordinates,
$p,r \in\bn$, $\lb\in\bc$, deg$R\leq r$,
$R(0)\neq 0$, and as per III.1, after blowing up
if necessary, we may suppose that $x=0$
defines the exceptional divisor $E$.

In addition to immediately applying the previous
conventions to take a generator $\pa$ convergent
in analytic functions $x,y,z$ which agrees
with the normal form modulo $\cO(-N E)$ for
some $N\gg 0$, define functions of a complex
variable $x$ by way of,
$$\z'(x) = -\dfrac{R(x)}{x^{r+1}}(1+\lb x^{p+r}),\,\,\,\,
\xi'(x) = -(p+r) \dfrac{1+\lb x^{p+r}}{x^{p+r+1}}$$
where it is convenient to view $\z,\xi$ as 
functions of $X=x^{-1}$, {\it i.e.} in
neighbourhoods of infinity. As such,
$\xi (X)$ is asymptotically $X^{p+r}$,
and $\z(X)$, $c_0 X^r$ for some constant
$c_0\neq 0$ (equals $R(0)/r$) whose value
is important. Were the normal form to be
convergent, we would then have first integrals,
$$s=z\exp(\xi),\,\,\,\, t=y\exp(\z)$$
Unlike the cases encountered in \S II.1-II.3,
the existence domain for a conjugation to the
normal form are complicated, \cite{mp1},
\S IV.2-IV.3, VI.4.5, nevertheless for $\Re(\z)>0$
they share many of the salient features
already encountered, and 
by way of an illustrative, albeit logically
irrelevant to the proof of the main lemma, complement
we may deduce,

\noindent{\bf III.3.2 Fact} Provided imaginary
rays in the plane of $\z$ do not go to imaginary
rays in the plane of $\xi$, {\it i.e.} 
no $j$ satsfies $\Re(c_0 j^r)=\Re(j^{p+r})=0$, 
there are  open sectors strictly (at both ends) containing
$\Re(\z) >0$ such that any 
invariant measure with
zero Segre class must be supported in 
the exceptional divisor.
Otherwise, idem, but for any open sectors
strictly contained  in $\Re(\z) >0$.

\noindent{\bf proof} This is basically an
appendix to \cite{mp1} \S IV, and one should
probably have a copy to hand, the notation
being the same except that $\zeta$ here is
$z$ in op. cit. This said, one is basically trying
to bring the field into the form $-\frac{\pa}{\pa\xi}$
in $\xi, s,t$ variables, which results in a
fibration,
$$
\begin{CD}
U @>{b\times \xi}>> B\times\bc \\
@VV{b=(s,t)}V @. \\
B\sbs\bc^2
\end{CD}
$$
with, as ever, all maps open and the horizontal
arrow an embedding. Whence, modulo connectedness
issues, the leaves may be identified with the
fibres $U_b$. Supposing for the sake of argument $y,z$
varying in unit discs, the fibres have the form,
$$\log |s| < \Re (\xi),\,\,\, \log\mt < \Re (\z )$$
Now for $\Re(\z) > 0$ the leaves are, on the whole,
unbounded, but some cases are more unbounded
than others, {\it viz:} if in the plane of $\xi$
the restriction imposed by $t$ is also open
to $\Re(\xi)\ra +\infty$ then we'll say that
we're in a thick sector, otherwise, it will be
called thin. Observe furthermore that for zero
Segre class the technique affording I.3.4
implies that $\mx^{-(r+2)}dx\otimes d\bar{x}$
is absolutely integrable. Indeed,
\begin{align*}
\oint_{\substack{\mxe\\ \my\leq A\\ \mz\leq B}} \dfrac{dx}{x}\dm 
\,+ \, \oint_{\substack{\my= A\\ \mxle\\ \mz\leq B}} \Re\{\dfrac{\pa x y}{x\pa y}\}\dfrac{dy}{y}\dm \,
+\, \oint_{\substack{\mz=B\\ \mxle\\ \my\leq A}} \Re\{\dfrac{\pa x z}{x\pa z}\}\dfrac{dz}{z}\dm 
\, =\, 0,\,\, \e>0
\end{align*}
the latter 2 terms are bounded by the measure
on the appropriate face in an $\e$ neighbourhood
of the exceptional divisor times $\e^r$ and $\e^{p+r}$
respectively, so the initial term is at worst $\e^r o(\e)$,
and applying Stokes yields the assertion.
It follows that the set of leaves where the
integral of the density $\mx^{-(r+2)}dx\otimes d\bar{x}$
is infinite has zero measure. This is certainly
true of leaves in thick sectors with $\Re(\z) >0$, but
there may be a doubt about it for $p\gg r$ in
thin sectors. As such we need to add some
comment to \cite{mp1} around the imaginary
axis in $\z$ when the imaginary axis is thick.
In this situation, op. cit. IV.3.2 (a)-(d)
apply without change to construct a sector
around the imaginary axis where the existence
domain can be taken to be a bi-disc in $(y,z)$
and a sector in $x$, {\it i.e.} thick sectors
in $\Re(\z) > 0$ bordering on the imaginary
axis in $\z$ continue to admit a conjugation
into the half space $\Re(\z)<0$ of the analytic
field to its normal form without taking a logarithm
in $y$, albeit perhaps on quite a small sector
if $p\gg r$.

This implies that the only possibility for the
measure to be supported in $\Re(\z)>0$ are
thin sectors which meet an imaginary line,
since, otherwise, the leaf may be analytically
continued from a thin region to a thick one
where there is no measure. Now, there are cases
to consider. The generic one is that the imaginary
axis in $\z$ is strictly thin, {\it i.e.} $\Re(\xi)$
goes to plus or minus $\infty$ along it with the
asymptotics of some nearby non-imaginary ray.
Here one could solve the centre manifold
problem of op. cit. \S IV with an existence
domain sectorial in $x$ and a bi-disc in
$(y,z)$ but maybe not the conjugation to
normal form in a way better than op. cit.
IV.2.5.  Fortunately the difference is slight.
In the case where it can be done, {\it i.e.}
the imaginary axis has $\Re(\xi)\ra +\infty$,
the angle between the $s$ and $t$ level curves
never goes to zero, so, in fact, 
$\mx^{-2}\cd$, and not just $\mx^{-(r+2)}\cd$,
has infinite measure even in leaves in leaves
in sectors around the imaginary axis. In the
other case the only leaves which can't be
analytically continued into the neighbouring
thick sector must be in a bounded region for
the variable $X$, actually an annulus of bounded 
moduli with inner boundary the boundary for
the $X$ variable itself, so, shrinking $x$
as necessary, all leaves in this sector will
continue to a thick sector where there
is zero measure.

This leaves the possibility that purely
imaginary in $\z$ goes to purely
imaginary in $\xi$, equivalently the
imaginary axis is critical in the
notation of op. cit., so for $j$ the
point of $S^1$ in $X$-space we have
$c_0 j^r$ and $j^{p+r}$ purely imaginary.
Obviously rare, but not impossible.
If this situation doesn't occur we could have taken
the asserted sector of zero measure to include
not just $\Re(\z) > 0$, but also extend
around both imaginary axis. When it does
occur, 
one
must just accept a small loss, move slightly
(sectorial sense) into the domain of $\Re(\z)>0$,
and argue as above in the cases where the
imaginary axis is strictly thin. $\Box$

This is somewhat different as to how
we proceeded  in
\S II, so let us offer:

\noindent{\bf III.3.3 Remark} To the extent that
one is prepared to invoke \cite{mp1} one
could reasonably pursue a somewhat stronger 
proposition without hypothesis on the
Segre class \`a la II.3.2/3. This can be
done with some variation on the above. The
problem in either case is that while elementary
the existence domains in the region $\Re(\z)>0$
depend on $p,r$ and even $R(x)$ in a rather
fastidious manner, and may be well short of
width $2\pi$ in the $\xi$-variable. As such,
op. cit. takes the point of view that there
is a finite division of the $\z$-plane with
a conjugation to the normal form in every region,
albeit possibly with further restriction on $\log y$
if $\Re(\z)<0$. Consequently to compute the
holonomy of a transversal $z=$constant, one
has to patch around to find the infinite
orbit which would eliminate the possibility
of an invariant measure. 
Whence, the best way to prove this
stronger proposition is to employ an
approximate holonomy estimate by way
of the evident variant of III.1.1 for
the holonomy around a loop in the
curve $x=y=0$, albeit this would repeat
much of \cite{mp1} \S IV.3 in another guise. 
In what follows, however,
although we will not use III.3.2, we
will employ the fact encountered in
the proof that, in the presence of
zero Segre class,
$\mx^{-(r+2)}\cd$ is 
absolutely integrable against $\dm$.

Now, in applying the basic construction III.1.1
around loops in the singular locus $x=z=0$,
with, as before $\my=\rho$ for $\rho$ in
some interval $I$ bounded away from 0, it's
helpful to have another normal form at
our disposition, {\it viz:}

\noindent{\bf III.3.1 (b) Alternative}
As per III.3.1 (a) but with the
normal form,
$$z\dz + x^p(\tr (x) + 
\tilde{\lb} x^{p+r})\left(y\dy + \dfrac{x^{r+1}}{r(1+\nu x^r)}\dx\right)$$
for $\tilde{\lb}, \nu\in\bc$ and deg$\tr\leq r-1$.

Each normal form has its own role, and III.3.1(b)
is more convenient for the almost holonomy
around $\g$ when $\my=\rho$, as opposed
to the actual, or even approximate, holonomy
around $\mz=\rho$ which underlies III.3.2.
In particular we now suppose that a convergent
generator $\pa$ in functions $x,y,z$ approximates
III.3.1 (b) to some large order modulo $\cO(-N E)$,
$N\in\bn$ to be specified. Denoting by $t$ the
resulting approximately invariant function
envisaged by III.1.1, and $\g$ the loop, or
its canonical image in $\tm$, for any function
of $x$ we may express it as a function of $t$,
and vice versa, by way of $x(t)=x|_T$, for $T$
our plane transversal. As such for $\z$ (basically
the same function as before in different clothes)
equal to $x^{-r} + \nu \log x^{-r}$ the approximate
holonomy is given by,
$$h^*\z = \z +\g$$
for any determination of $\log x$ in a sector
of width $2\pi/r$ in the argument of $x$.

Observe that 
a convenient norm to work with is,
$$\nx^{-r}=\min\{ |\Re(x^{-r})|,\, |\Im(x^{-r}|\}$$
In particular $R=\e^{-r}$, modulo the minor notational
confusion with III.3.1 (a), is often easier to 
work with, while, in addition, $t(x)=x(1+ O(x^r))$.
Certainly there is a branching issue in $\z$ but
in practice it will pose no problem, while if
we suppose III.3.2, it
may be taken in the domain of null measure.
Either way,
$|\Re(\z)|\geq R$ and $\Im (\z)|\geq R$
is often the most convenient domain to work with.
There may also an orientation issue as to whether
we have to worry about there being measure in the
upper or lower half planes in $\z$, so, say the
former, albeit this implies that we'll be doing
III.3.1 with loops oriented in the opposite direction
to that implied by the choice of an imaginary part
function.

All of which said, we have the following co-equaliser
estimate,
$$\int_{\substack{R> \Im(\z) > R+ 2\pi\\ \mzle, y=\bfp}} \dm \leq
\int_{\substack{ \Im(\z(t))=R\\ \mzle, \my=\rho}} \mmu 
+ \int_{\substack{\mze, \my=\rho\\ \Im(\z) < R}} \mmu $$
Evidently, we've already encountered how to estimate
the terms in question and choosing the order of
the approximation to the normal form to be
sufficiently high, we have the estimates,

\noindent{\bf III.3.4 (a) Estimate} For a constant
$C$ depending on the interval $I=(a,b)$ in which
$\my$ varies, and $R$ outside a set of finite
$R^{-1}dR$ measure, 
\begin{align*}
\int_a^b d\rho\biggl(\sum_{n=0}^{\infty}  
& \int_{\substack{|\imzt|=R+n\\ \mz\leq (R+n)^{-1/r}\\\my=\rho}}   \mmu 
+ \sum_{m=-\lceil R\rceil}^{\lceil R\rceil} 
\Bigl[ \int_{\substack{\imzt=m\\ \mz\leq R^{-1/r},\my=\rho}}  \mmu  \\
+ &\int_{\substack{|\rezt|= R, |\imzt|\geq m\\ \mz\leq R^{-1/r}, \my=\rho}}
 \mmu \Bigr]\biggr) \leq C R^{-(2+ r)}
\end{align*}
where the implied possibility to extend the
approximate holonomy estimate over strips
with $\imz$ between $-R$ and $R$ results
from the fact that $\rezt$ is holonomy
invariant. This is akin to III.2.2 (b),
whereas the analogue of III.2.2 (a) is,

\noindent{\bf III.3.4 (b) Estimate} For a
constant $C$ depending only on the maximum
value $A$ of $\my$, we have,
$$\int_a^b d\rho \int_{\substack{\nx\leq R^{-1/r}, \mz=R^{-1/r}\\ \my =\rho}} \mmu
\leq CR^{-p/r} \sup_{R_0>R} \szm (R_0^{-1/r})$$

Before applying these estimates we need
a consequence of the fact encountered
in the proof of III.3.2 that 
$\mx^{-(r+2)}\cd$ is
absolutely integrable against $\dm$
whenever the Segre class is zero.
More precisely for any $R<B<\infty$,
\begin{align*}
& \int_{\substack{\my=\rho\\ R\leq \imx\leq B\\ |\rex|\leq R, \mz\leq R^{-1/r}}} d\rex \dm
\,+\, \int_{\substack{\mz=R^{-1/r}\\ R\leq \imx\leq B\\ |\rex|\leq R, \my\leq\rho}} d\rex \dm 
\,=\, \\
& \int_{\substack{\imx=B \\ |\rex|\leq R\\ \my\leq\rho, \mz\leq R^{-1/r}}} d\rex \dm
\,-\, \int_{\substack{\imx=R \\ |\rex|\leq R\\ \my\leq\rho, \mz\leq R^{-1/r}}} d\rex \dm
\end{align*}
The integral, $I(B)$, say, over the $\imz=B$ face may be
integrated against $B^{-1}dB$, and by the
aforesaid absolutely integrability of
$\mx^{-(r+2)}\cd$ the result is finite,
so the liminf of $I(B)$ as $B\ra\infty$
must be zero. Whence, 
\begin{align*}
& \int_{\substack{|\rex|\leq  R, \imxr\\ \my\leq \rho,\mz\leq\re}} d\rex  \dm\,=\, \\
-4\pi & \oint_{\substack{\my=\rho, \mz\leq\re\\ |\rex|\leq R, |\imx|\geq R}} \Im\left(\dfrac{x^{-(r+1)}\pa x y}{\pa y}\right)\dfrac{dy}{y}
 \dm \,\, - \\
4\pi &  \oint_{\substack{\mz=\re, \my\leq\rho\\|\rex|\leq R, |\imx|\geq R  }}
 \Im\left(\dfrac{x^{-(r+1)}\pa x z}{\pa z}\right)\dfrac{dz}{z}\dm
\end{align*}
At which point, observe that there is a considerable
tangency along the foliation between the curves $\my$=const.,
and $\rex$=const., {\it i.e.} for $\nx\leq\re$, and $\my$
confined to $I$ there is a constant $C$ depending on $I$
such that,
$$\oint_{\substack{\my=\rho\\ \nx,\mz\leq\re}}
\biggl\lvert\Im\left(\dfrac{x^{-(r+1)}\pa x y}{\pa y}\right)\biggr\rvert\dm\,\leq\,
C\oint_{\substack{\my=\rho\\ \nx,\mz\leq\re}} |x^r|\dfrac{dy}{y}\dm$$
and, of course, a sufficiently high order of approximation
to the normal form. Similarly, for such an approximation,
and another $C$ depending only on the same,
$$ \oint_{\substack{\mz=\re\\ \nx\leq\re,\my=\rho}}
\biggl\lvert\Im\left(\dfrac{x^{-(r+1)}\pa x z}{\pa z}\right)\biggr\rvert\dm\,\leq\,
CR^{-p/r}\szm(\re)$$
The estimates III.3.4 (a)-(b) do not depend on the
point $\bfp$ in $\g$ where we take the discontinuity,
so summing the inequality arising from the basic
construction III.3.1 with quantification implied
by III.3.4 (a)-(b), and averaging over $\bfp$, we have,
$$\int_a^b d\rho\oint_{\my=\rho,\nx,\mz\leq\re}
\mx^r\dfrac{dy}{y}\dm \leq
CR^{-p/r}\max\{\sup_{R_0>R} \szm (R_0^{-1/r}), R^{-1}\}$$
Arguing similarly for $\rex=-R$, and
combining all of which, we obtain,

\noindent{\bf III.3.5 Estimate} There is a constant $C$
depending on the interval $I$ in which $\my$ varies
such that for $R$ outside a set of finite $R^{-1}dR$
measure, for any $B>0$ and $\rho$ excluded from a
set (depending on $R$) of measure at most $B^{-1}$ in $I$,
\begin{align*}
 \int_{\substack{|\rex|\leq  R\\ |\imx|=R\\ \my=\rho,\mz\leq\re}}
|d\rex|\dm 
 \leq BCR^{-p/r}\max\{\sup_{R_0>R} \szm (R_0^{-1/r}), R^{-1}\}
\end{align*}

The integrand above being positive for the given
choice of orientation, we have an evident 
improvement in the quantification in $\rho$
at our disposition. We now require a similar
shape  of estimate around the boundary parallel
to the imaginary axis, to wit, for $\chi$ to be chosen:
\begin{align*}
& \int_{\substack{|\imx|\leq R,\rex=-\chi R\\ \my\leq \rho,\mz\leq\re}} d\imx\dm =
 \int_{\substack{\imxr,\rex=-\chi R\\ \my\leq \rho,\mz\leq\re}} \imx\dm\,\, +\\
& \int_{\substack{|\imx|\leq R,\rex=-\chi R\\ \my= \rho,\mz\leq\re}} \imx\dm\,\, +
 \int_{\substack{|\imx|\leq R,\rex=-\chi R\\ \my\leq \rho,\mz=\re}} \imx\dm 
\end{align*}
To deal with the first of the integrals on the right
we need to use III.3.5, {\it i.e.}
\begin{align*}
&\int_{e^{-1}}^1\dfrac{d\chi}{\chi}\int_{\substack{\imxr\\ \rex=-\chi R\\ \my\leq \rho,\mz\leq\re}} \mmu =
\int_{\substack{\imxr \\ R\leq \rex\leq- e^{-1} R\\ \my\leq \rho,\mz\leq\re}} 
\Bigl\lvert\dfrac{d\rex}{\rex}\Bigr\rvert\dm  \\
& \leq \dfrac{e}{ R}  \int_{\substack{\imxr,|\rex|\leq R\\ \my\leq \rho,\mz\leq\re}} 
d\rex\dm
\end{align*}
Integrating both sides against $d\rho$, $\rho\in I$ is a
little more convenient than applying III.3.5 directly,
and leads to:

\noindent{\bf III.3.6 (a) Estimate} Quantifiers
as per III.3.5 but now with $\log\chi\times\rho$
varying in $(-1,0)\times (a,b)$ excluded from 
a set depending on R of measure at most $B^{-1}$,
then:
\begin{align*}
 R\int_{\substack{\imxr,\rex=-\chi R\\ \my\leq \rho,\mz\leq\re}} \mmu \leq
BCe
R^{-p/r}\max\{\sup_{R_0>R} \szm (R_0^{-1/r}), R^{-1}\}
\end{align*}

Of the two remaining terms the easier is that over
the face $\mz=\re$, which integrates against $\chi^{-1}d\chi$
to give,
\begin{align*}
 \int_{e^{-1}}^1 \dfrac{d\chi}{\chi}
\int_{\substack{|\imx|\leq R\\ \rex=-\chi R\\ \my\leq \rho,\mz=\re}} \dm & =
 \int_{\substack{|\imx|\leq R\\ - R\leq \rex\leq- e^{-1} R\\ \my\leq \rho,\mz=\re}} 
\Bigl\lvert\dfrac{d\rex}{\rex}\dm \Bigr\rvert  \\
& \leq \dfrac{eC}{ R^{1+p/r}}\szm (\re)
\end{align*}
for a constant $C$ depending only on the
order of approximation, so:

\noindent{\bf III.3.6 (b) Estimate}
Quantifiers as above, {\it i.e.} nothing
to worry about except the order of 
approximation, then for $\log\chi\in (-1,0)$
excluded from a set (depending on R) of
measure at most $B^{-1}$ for any $B$,
\begin{align*}
\Bigl\lvert \int_{\substack{|\imx|\leq R,\rex=-\chi R\\ \my\leq \rho,\mz=\re}} \imx \dm\Bigr\rvert\leq
\dfrac{eCB}{ R^{p/r}}\szm (\re)
\end{align*}

The final term will again employ the tangency
along the foliation between $\my$=const. and
$\rex$=const., {\it viz:}
\begin{align*}
\int_a^b\dfrac{d\rho}{\rho}
 \int_{e^{-1}}^1 \dfrac{d\chi}{\chi}
\int_{\substack{|\imx|\leq R\\ \rex=-\chi R\\ \my= \rho,\mz\leq\re}} \mmu =
 \int_{\substack{|\imx|\leq R\\ - R\leq \rex\leq - e^{-1} R\\ \my\in I,\mz\leq\re}} 
\Bigl\lvert\dfrac{d\my}{\my}\dfrac{d\rex}{\rex}\Bigr\rvert \dm 
\end{align*}
While the integrand on the right may be
bounded by,
$$C \Bigl\lvert\dfrac{\Im(\nu x^r}{\rex}\Bigr\rvert\dfrac{dy\otimes d\bar{y}}{\my^2}$$
for a constant $C$ depending on $I$ and the order
of approximation, so that:
\begin{align*}
\int_a^b & \dfrac{d\rho}{\rho} 
 \int_{e^{-1}}^1 \dfrac{d\chi}{\chi}
\Bigl\lvert \int_{\substack{|\imx|\leq R\\ \rex=-\chi R\\ \my= \rho,\mz\leq\re}} \imx\dm\Bigr\rvert \\
\leq & \dfrac{eC\nu}{ R}  \int_{\substack{|\imx|\leq R\\ - R\leq \rex\leq - e^{-1} R\\ \my\in I,\mz\leq\re}} \dfrac{dy\otimes d\bar{y}}{\my^2}\dm
\end{align*}
As before to estimate the term on the right one
again uses co-equalisers, about $2R$ of them,
of the form $\imx>m$, $m$ between $-\lceil R \rceil$ and $\lceil R \rceil$,
and $\rex\leq - e^{-1} R$, or, better, idem
for $\imz$ and $\rez$ which are as near to the
same thing as makes no difference, and whence:

\noindent{\bf III.3.6 (c) Estimate}
For quantifiers as per III.3.6 (a),
\begin{align*}
\int_a^b & \dfrac{d\rho}{\rho}
 \int_{e^{-1}}^1 \dfrac{d\chi}{\chi}
\Bigl\lvert \int_{\substack{|\imx|\leq R\\ \rex=-\chi R\\ \my= \rho,\mz\leq\re}} 
\imx\dm\Bigr\rvert\\ 
\leq & \dfrac{BeC\nu}{ R^{p/r}}\max\{\sup_{R>R_0} \szm (\re), R^{-1}\}
\end{align*}

In order to lighten the notation a little, for
every $\chi$ between $e^{-1}$ and $1$, put,
$$\nx^{-r}_{\chi}=\min\{ ( \chi)^{-1} |\Re(x^{-r})|, |\Im(x^{-r}|\}$$
and profit from the positivity of the integrands
around the real and imaginary boundaries to
improve the quantification in $\rho$ by way
of replacing $I$ by an interval $(a', b')\sbs (0,a)$ to obtain
on repeating the above for $\rex=R$,

\noindent{\bf III.3.7 Estimate} There is a constant $C$
depending on the interval $I$, such that in the 
presence of zero Segre class, for $R$ outside a
set of finite $R^{-1}dR$ measure, $\rho\in I$, and
given $B$ a set (depending on R) of $\log\chi\in (-1,0)$
of measure $B^{-1}$ outwith which,
\begin{align*}
\int_{\substack{\nx_\chi=\re\\ \mz\leq\re\\ \my\leq\rho}}
& \bigl(|d\imx| + |d\rex\bigr|)\dm 
\leq\\ 
& \dfrac{BC}{R^{p/r}}\max\bigl(\sup_{R_0 > R} \szm (R_0^{-1/r}), R^{-1}\bigr)
\end{align*}

Plainly  this
computes any residue $dx^{-i}$, $1\leq i\leq p+r$,
{\it i.e.}

\noindent{\bf III.3.8 Corollary} Everything
as above in III.3.7, then for $\chi$ understood
to depend on $R$ in the way implied by III.3.7,
and $1\leq i\leq p+r+1$,
$$
\lim_{R\ra\infty} \, \biggl\lvert
\int_{\substack{\nx_\chi =\re\\ \mz\leq\re, \my\leq\rho}}
\dfrac{dx}{x^i}\dm\biggr\rvert \, =\, 0
$$

It also has the further important corollary
for globalisation,

\noindent{\bf III.3.9 Corollary}
At the risk of a certain notational
confusion for $\Vert \, \Vert$ any
norm in the $x$-variable, $\my$
limited to $I$ as above, $o(\e)\ra 0$
as $\e\ra 0$ a small function 
depending on the limiting Segre
classes (with limit supposed zero)
and the norm,
$$
\oint_{\substack{\my=\rho\\ \nx,\mz\leq\e}}
\dfrac{dy}{y}\dm\, = \, \e^po(\e)
$$

\noindent{\bf proof} Under the 
zero Segre class hypothesis, and,
of course $\dm$ without support
in the exceptional divisor, we
have as per the proof of III.3.2,
\begin{align*}
\int_{\substack{\my=\rho\\ \nx, \mz\leq\e}}
\o\dm +
\int_{\substack{\nx=\e\\ \my\leq\rho, \mz\leq\e}}
\o\dm +
\int_{\substack{\mze\\ \my\leq \rho, \nx\leq\e}}
\o\dm\, = \, 0
\end{align*}
For $\o$ the 1-form,
$$(1+\nu x^r)\dfrac{dx}{x^{r+1}}$$
The last integral obviously has
modulus $\e^po(\e)$, and the middle
one too by III.3.7 provided that
$\Vert\,\Vert$ is taken to be $\Vert\,\Vert_\chi$.
Meanwhile, on the face $\my=\rho$ identified
with $\g\in\tm$,
$$\dfrac{dy}{y}\dm\,=\, -(1+ O(x^N))\o\dm$$
so dividing out the first integral by $-\g$
it is the asserted quantity up to an
irrelevant error. The asserted quantity is,
however, the average measure of an $\e$-ball
in  a translate, so,  the
estimate for some norm implies it for any. $\Box$
 
Even were we to make use of III.3.2,
the above complication would not be
lessened because of the possibility
that, under very particular conditions,
imaginary $\z$ lines could
go
to imaginary $\xi$ lines, so, we have to
worry about (a probably imaginary) measure
around the boundary $|\imx|=R$. 
In the case of
exceptional nodes we have no such
problem, and things are a lot easier, {\it i.e.}

\noindent{\bf III.3.1 (c) Set Up}
As per III.1.1 (a), but as per II.3.1
a formal generator in the completion
at the point is given by,
$$qz\dz + \dfrac{x^{p}y^q}{1 + y^q(R(x)+\lb x^{p+r})}\left(y\dy+
\dfrac{x^r}{1+\nu x^r}(qx\dx-py\dy)\right)$$
where $x,y,z$ are formal coordinates,
$p,q,r \in\bn$, $\lb,\nu\in\bc$, 
deg$R\leq r-1$, and the exceptional
divisor $E$ has 2 components $E_1:x=0$,
$E_2:y=0$, and we may even assume that
the normal form is defined after 
completion in the former, albeit 
not the latter.

In our immediate context, take a
sufficiently good approximation
modulo a large power of $E_1$
in convergent variables $x,y,z$,
and consider for $\chi\leq 1$,
\begin{align*}
&\int_{\substack{\rex=-R\chi\\ |\imx|\leq R\\ \my\leq\rho, \mz\leq\re}}
d\imx\dm\,=\\
&\int_{\substack{\rex=-R\chi\\ |\imx|\leq R\\ \my=\rho, \mz\leq\re}}\imx\dm \, +\,
\int_{\substack{\rex=-R\chi\\ |\imx|\leq R\\ \my\leq\rho, \mz=\re}}\imx\dm
\end{align*}
which is much simpler than before
thanks to II.3.3, and we can argue
as previously, in fact with quite a
lot of simplification, to conclude:

\noindent{\bf III.3.10 Fact}
As per III.3.9,
$\Vert \, \Vert$ any
norm in the $x$-variable, 
and supposing zero Segre
class, with
$\my$
limited to $I$ as above, 
there is a function
$o(\e)\ra 0$
as $\e\ra 0$ 
such that,
$$
\oint_{\substack{\my=\rho\\ \nx,\mz\leq\e}}
\dfrac{dy}{y}\dm\, = \, \e^po(\e)
$$

At the risk of a certain notational
confusion we need a similar estimate
on faces $\mx=\rho$, for $\mx$ confined
to an interval $I$ bounded away from
zero. The quantity in question is
simply the average measure of
transversals at $\bfq$ in the loop
$\g:\mx=\rho$, and we already know
that this is zero for $\bfq$ outside
of a strip around the negative real
axis in the $x^{-r}$ plane by II.3.3.
Evidently for $\bfp,\bfq\in\g$ such
that the measure at $\bfp$ is zero
we'd like to employ a variant of
the basic construction III.1.1 by
way of an approximation of the form:

\noindent{\bf III.1.1 bis Variant} 
Define a co-dimension 1 
real manifold $\G_\e$ by,
$$x\in(\bfp,\bfq),\,\, \mzle, \mtle$$
where $t$ is a function of $x,y$
which is an approximately invariant
function defining $E_2$.

Now, while generically a fairly trivial 
variant, here we should be cautious
because III.3.1 (c) is not valid
modulo arbitrary powers of $E_2$.
The reason for this problem,
\cite{mp1} VI.2.1 (g) \& VI.2.2 (g)
is because the operators,
$$\un-\a x^{r+1}\dx$$
for $\a$ a constant, cannot be inverted in convergent power
series.  They can, however, be inverted
in domains of argument up to $3\pi/r$
in the $x$-variable. To achieve the
normal form one has to invert several
of these with $\a$ going through finitely
many positive and negative integers, which
implies a certain incompatibility in
the branching, so that ultimately
one achieves 

\noindent{\bf III.3.11 Fact}
For $N\in\bn$, $x$
confined to a sector $S_N$ of aperture
up to $2\pi/r$ branched within $\pi/2$
of the negative real axis, and $y,z$
varying in a bi-disc $\D^2_N$, the normal
form III.3.1 (b) can be realised modulo
$\cO(-N E_2)$.

Plainly given II.3.3 this will be
adequate for our current considerations
and we apply it as follows: take $N\in\bn$
to be chosen, restrict ourselves to
$U_N=S_N\times\D^2_N$, take $t$ to
be our approximation (basically $y\exp(x^{-r})$)
such that $\pa t\in (y^N)$, confine $\mx$
to an interval $I_N$ bounded away from
$0$ with $\g_{\bfp,\bfq}$ a simple
path in $\g$ joining $\bfp,\bfq$, then:
\begin{align*}
\int_{\substack{T_\bfp\\ \mt,\mz\leq\e}} \dm\,\leq\,
\int_{\substack{\mze,\mtle\\ x\in\g_{\bfp,\bfq}}} \mmu
\, +\, \int_{\substack{\mzle,\mte\\ x\in\g_{\bfp,\bfq}}} \mmu
\end{align*}
We estimate the integrals on the right, or, better,
their averages in the problematic region exactly
as in III.2.2 (a)/(b) respectively, and obtain:

\noindent{\bf III.3.12 Fact} There is an interval
$I$ bounded away from zero such that in the
presence of zero Segre class, there is a function
$o(\e)\ra 0$ as $\e\ra 0$ for which if $\rho\in I$
is excluded from a set of measure $B^{-1}$, then,
$$\oint_{\substack{\mx=\rho\\ \my,\mz\leq\e}}
\dfrac{dx}{x}\dm\,=\,B \e^q o(\e)$$

Here the question of zero Segre class is very
much a hypothesis of convenience, which for
$\xi=(x^py^q)^{-1}$ we may apply to compute
residues as follows,
\begin{align*}
0= & \oint_{\substack{|\xi|=\e^{-1}, \mzle\\ \mx\leq X \my\leq Y}} \dfrac{d\xi}{\xi}\dm \, 
+\, \oint_{\substack{|\xi|\leq\e^{-1}, \mze\\ \mx\leq X \my\leq Y}} 
\Re\Bigl(\dfrac{\pa\xi z}{\xi\pa z}\Bigr)\dfrac{dz}{z}\dm \\
+ & \oint_{\substack{|x|=X, \mzle\\ \my^q\leq\e X^{-p} }}
\Re\Bigl(\dfrac{\pa\xi x}{\xi\pa x}\Bigr)\dfrac{dx}{x}\dm
\, +\,\oint_{\substack{\my=Y, \mzle\\ \mx^p\leq \e Y^{-q}}} 
\Re\Bigl(\dfrac{\pa\xi y}{\xi\pa y}\Bigr)\dfrac{dy}{y}\dm
\end{align*}
where the coordinates are chosen so that $x=0$, $y=0$
are the invariant exceptional divisors, the
approximation is as good as we like modulo $E_1$,
and we use \cite{mp1} V.1.9 to guarantee $\pa z=z$
modulo $(x,y)^N$ for $N\in\bn$ as large as we
please, while $X, Y$ are just things in the
intervals appearing in III.3.9/12. In particular,
by the above, the final two integrals are
$\e o(\e)$, while $\pa z =z$ mod $\xi^N$,
as ever $N$ large,  so the term on the $\mze$
face goes like $\e$ times the Segre class, whence:

\noindent{\bf III.3.13 Corollary} Let things be as in
III.3.1 (c) with $X, Y$ arising from the intervals
of III.3.9 and III.3.12 respectively, then for $o(\e)\ra 0$
as $\e\ra 0$,
$$ 
\oint_{\substack{|\xi|=\e^{-1}, \mzle\\ \mx\leq X \my\leq Y}} \dfrac{d\xi}{\xi}\dm \, =\,\e o(\e)
$$
and (locally at least) all residues of interest are zero.

\subsection*{III.4 Other linear(ish) singularities}

Let us gather together the remaining possibilities
for the singularities that we have to deal with
starting with the most straightforward,

\noindent{\bf III.4.1(a) Set Up} In a
neighbourhood $U$ of 
foliated 3-fold of a point in the
singular locus after completion in the same,
the foliation admits a formal generator of
the form,
$$z\dz + \dfrac{x^{p}y^q}{1+\nu x^py^q}\left(x\dx + \lb y\dy \right)$$
where $x,y,z$ are formal coordinates,
$p,q \in\bn$, $\lb\in\bc\bsh\br $, 
($\Re(p+q\lb)$, $\Re(p/\lb +q)$ not
both negative if we avail ourselves
of II.2.1 (b)), and sufficient a priori
 blowing up having been performed to
guarantee that
$x=0$, $y=0$ are local equations for
exceptional divisors $E_1$, and $E_2$
respectively, with $E= E_1 + E_2$.

Here with exactly the same proof 
(and, in fact, somewhat easier since
the convergence of $y=0$ is given)
III.2.2 (c) or III.2.2 (c) bis. holds.
The residue calculation is, therefore,
particularly straightforward since, for,
say $\mx\leq X$, $\my\leq Y$ in the
the domain of the appropriate coordinates,
and again $\xi=(x^py^q)^{-1}$,
\begin{align*}
0= & \oint_{\substack{|\xi|=\e^{-1}, \mzle\\ \mx\leq X \my\leq Y}} \dfrac{d\xi}{\xi}\dm \,
+\, \oint_{\substack{|\xi|\leq\e^{-1}, \mze\\ \mx\leq X \my\leq Y}} 
\Re\Bigl(\dfrac{\pa\xi z}{\xi\pa z}\Bigr)\dfrac{dz}{z}\dm \\
+ & \oint_{\substack{|x|=X, \mzle\\ \my^q\leq\e X^{-p} }}
\Re\Bigl(\dfrac{\pa\xi x}{\xi\pa x}\Bigr)\dfrac{dx}{x}\dm
\, +\,\oint_{\substack{\my=Y, \mzle\\ \mx^p\leq \e Y^{-q}}} 
\Re\Bigl(\dfrac{\pa\xi y}{\xi\pa y}\Bigr)\dfrac{dy}{y}\dm
\end{align*}
under the hypothesis of zero Segre class.
Again, the second term is $\e\szme$,
the other two are governed by III.2.2 (c), so:

\noindent{\bf III.4.2 Fact} Shrinking $X$, $Y$ a
little to simplify the quantification, then
for $\e$ outwith a set of finite $\lem$ measure,
and zero Segre class around $E$,
$$ 
\oint_{\substack{|\xi|=\e^{-1}, \mzle\\ \mx\leq X \my\leq Y}} \dfrac{d\xi}{\xi}\dm \, =\,\e o(\e)
$$
Whence, any residue that we need to calculate is
not just zero locally in a natural way, but the
measure admits best possible estimates of the
form III.2.2 (c) for its mass on transversals.

The next cases are quite different, {\it i.e.}

\noindent{\bf III.4.1(b) Set Up} 
Exactly as per III.4.1 (a) with $\lb\in\br_{<0}\bsh\bq$,
or, $k,l\in\bn$ relatively prime,
$$z\dz + \dfrac{x^iy^j(x^{k}y^l)^n}{1+ x^iy^j\nu( x^ky^l)}\left(lx\dx -k y\dy \right)$$
where $x,y,z$ are formal coordinates,
$n,i,j \in\bz_{\geq 0}$, 
$\nu$ a formal function in one variable.
In the latter case, when choosing an
approximation, it is, therefore, to be
understood that $\nu$ is truncated to
some appropriately high order.

Plainly the approximately invariant
function $t=yx^{-\lb}$ may have 
holonomy, but $\mt$ does not- one
could use $x^ky^l$ in the rational
case, but that only leads to a
simplification globally if $i=j=0$.
In any case since $\lb\in\br_{<0}$,
the boundary $\max\{\mz,\mt\}=\e$
is a perfectly good way to
calculate the residue. Further,
for $N\in\bn$ as large as we like,
we have $\pa t$ at worst of the
form $t^{N+1}\pa x/x$, while any
1-form $\o$ of which we have to
compute the residue is no worse
than $t^{-m}\pa x/x$ for some fixed
$m$. More precisely, locally the
residue has the form,
$$\int_{\mze,\mtle}\o\dm\, +\, \int_{\mte,\mzle}\o\dm$$
where the first term will never be worse
than the Segre class, while for 
$m\in\bn$  large, albeit smaller than say $N/2$,
\begin{align*}
\int_0^\e\dfrac{d\d}{\d^{m+1}}
\biggl\lvert \int_{\substack{\mtd\\ \mzld}}\o\dm\biggr\rvert
\,\leq\,
\int_{\substack{\mtle\\ \mzle}} \mt^{-m} \bigl\lvert \dfrac{dt}{t}\o\bigr\rvert \dm
\, \leq\, o(\e)
\end{align*} 
and whence,

\noindent{\bf III.4.3 Fact} For $\o$ a
1-form of which we require to calculate
the residue, for $\e$ outwith a set of
finite $\lem$ measure,
$$\biggl\lvert \int_{\max\{\mt,\mz\}}\o\dm\biggr\rvert
\,\leq\, \max\{\szme,\e\}$$
One should, however, bear in mind that this
is a local estimate, which lacks the good
bound on transversals occurring in III.4.1 (a),
so, patching it to the global procedure is
a little more delicate.

This leaves us with the general rational
case to tackle, {\it i.e.}

\noindent{\bf III.4.1(c) Set Up} 
Exactly as per III.4.1 (b) 
but with normal form,
\begin{align*}
z\dz + \dfrac{x^iy^j(x^{k}y^l)^n}{1+ x^iy^j (R(x^iy^j)+\lb( x^ky^l)^{n+r})}
\biggl( &  \bigl(lx\dx -k y\dy\bigr) \\ 
+ & \dfrac{(x^ky^l)^r}{1+\nu (x^ky^l)^r}
\bigl(jx\dx-iy\dy\bigr)
\biggr)
\end{align*}
where $x,y,z$ are formal coordinates,
$n,i,j \in\bz_{\geq 0}$, 
$\lb,\nu\in\bc$, $R$ of degree at most
$r-1$, $r\in\bn$, and one replaces the
fields $jx\dx-iy\dy$ by, say $x\dx$ if
$i=j=0$. Alternatively, should this 
occur, one may have recourse to the
normal form,
$$z\dz + (x^{k}y^l)^n
\left(R(x^ky^l)(lx\dx -k y\dy) + \dfrac{(x^ky^l)^r}{1+\nu (x^ky^l)^{n+r}}
x\dx
\right)$$
where, now, $R$ is of degree at most $r$, 
$R(0)\neq 0$, and all cases, $p=i+kn$,
$q=j+ln$ are in $\bn$.

Taking a convergent approximation to
an order to be decided, we have for
the approximating foliation, invariant
formal functions of the form,
\begin{align*}
t=xy^{l/k} \exp\bigl(\sum_{m=1}^\infty (x^ky^l)^{mr}P_m(\log y)\bigr),\,\,\,
s=yx^{k/l} \exp\bigl(\sum_{m=1}^\infty (x^ky^l)^{mr}Q_m(\log x)\bigr)
\end{align*}
for $P_m, Q_m$ polynomials of degree at most
$m$, and,
$$P_1(\log y) = k^{-2}(k\tilde{j} -l \tilde{i})\log y,\,\,
Q_1(\log x) = -l^{-2}(k\tilde{j} -l \tilde{i})\log x$$
where $(\tilde{i},\tilde{j})=(i,j)$ if one of $i,j$ in
III.4.1 (c) is non-zero, and $(1,0)$ otherwise.
As such, say $P_1(\log y) =p_1/k\log y$, and
$Q_1(\log x)=q_1/l\log x$ to lighten the notation,
and, observe, that these have formal holonomy
of the form,
\begin{align*}
& t\mpo \exp(\frac{l\g}{k})t\bigl(1+ \frac{p_1\g}{k} t^{kr} + O(t^{2kr})\bigr) \\
& s\mpo \exp(\frac{k\g}{l})s\bigl(1+ \frac{q_1\g}{l} s^{kr} + O(s^{2lr})\bigr) 
\end{align*}
for $\g$ a loop in $x=0$, respectively $y=0$,
identified to its canonical image in $\tm$.
Truncating $t$ and $s$ to an appropriately
high order, we therefore have approximately
invariant functions of the form,
\begin{align*}
& t^k = x^ky^l \bigl( 1+ p_1 (x^ky^l)^r\log y + O((x^ky^l)^{2r}\log^2 y)\bigr)\\
& s^l = x^ky^l \bigl( 1+ q_1 (x^ky^l)^r\log x + O((x^ky^l)^{2r}\log^2 x)\bigr)
\end{align*}
and an appropriate function to use in constructing a
boundary around the singularity is $\max\{\mz,\mw\}$,
where,
$$\mw:= \max\{\mt^k, \ms^l\}$$
Clearly $w$ is discontinuous, so an
appropriate estimate for the mass
will come from the almost holonomy
estimate.
To this end, observe that the approximating
holonomy has the form,
$$t^{-kr} \mpo t^{-kr} -rp_1\g + o(t)$$
whence for $F(\e)$ the region $\mt^k<\e$,
the sub-region in which, $\Re(t^{-kr}/p_1\g)\leq 0$
maps into itself, while that where
$\Re(t^{-kr}/p_1\g)\geq 0$ maps out of itself.
Consequently, we need both regions, say,
$F_-(\e)$ and $F_+(\e)$ to estimate the
discontinuity, {\it i.e.} their signs are
opposite in the basic construction, albeit:
$$\bigl\lvert \un_{F(\e)}(t) - \un_{F(\e)}(t^h)\bigr\rvert
\,\leq\,
\bigl\lvert  \un_{F_+(\e)}(t) - \un_{F_+(\e)}(t)\bigr\rvert
\,+\, \bigl\lvert  \un_{F_-(\e)}(t) - \un_{F_-(\e)}(t) \bigr\rvert
$$
Notice also, that $|t|$ doesn't change much in
the region of discontinuity, more precisely,
$$\dfrac{\bigl\lvert \un_{F(\e)}(t) - \un_{F(\e)}(t^h)\bigr\rvert}{\mt^{p}}
\,\leq\,
\dfrac{(1+O(\e^r))}{\e^{p/k}}\bigl\lvert \un_{F(\e)}(t) - \un_{F(\e)}(t^h)\bigr\rvert$$
With this in mind, say
$\log y$, $\log x$, are
branched along $b_y$, $b_x$ in $\br_+$ 
respectively, for $\my\leq Y$, respectively
$\mx\leq X$, in the spirit of
the variation III.1.3, so that $d\my$
may be taken as a length form along the
branch. Now divide through by $\e^{p/k}$, and integrate both sides of the 
basic construction III.1.1 against
$\my^{-1-q+pl/k}d\my$ to obtain,
\begin{align*}
& \int_{\substack{\mz,\mt^k\leq\e\\ y\in b_y,\my\leq Y}}
\dfrac{\bigl\lvert \un_{F(\e)}(t) - \un_{F(\e)}(t^h)\bigr\rvert}{\mt^{p}\my^{1+q-pl/k}}
d\my\dm \,\leq \, 2'\int_{\substack{\mze,\mt^k\leq\e\\ \my\leq Y}}
\Bigl\lvert\dfrac{d\my}{\mt^{p}\my^{1+q-pl/k}}\Bigl\lvert\dm\, +\\
&  \e^{-p/k}\Biggl\lvert  \int_{\substack{t\in\pa F_+(\e)\\ \mzle, \my\leq Y}} 
\dfrac{d\my}{\mt^{p/k}\my^{1+q-pl/k}}\dm\Biggr\rvert
\, +\,
\e^{-p/k}\Biggl\lvert\int_{\substack{t\in\pa F_-(\e)\\ \mzle, \my\leq Y}} 
\dfrac{d\my}{\mt^{p/k}\my^{1+q-pl/k}}\dm\Biggr\rvert
\end{align*}
For $2'$ anything bigger than 2 provided $\e$ is sufficiently small.
The choice of integrand was cooked up to be commensurate
with $|x^py^q|^{-1}d\log\my$, so the first integral on
the left is certainly no worse than $\szme$. As for the
other two, $\pa t$ is at worst $|t|^N\pa y/y$ for $N$
as large as we like, so arguing in the previous way,
these terms are as irrelevant as we wish to make them,
and whence,

\noindent{\bf III.4.4 Fact} For $\e$ outside a set of
finite $\lem$ measure we have the following estimate
for the discontinuity in $t$ (and similarly for $s$),
\begin{align*}
\int_{\substack{y\in b_y\\ \my\leq Y}}\dfrac{d\my}{\my^{q+1}}
\int_{\mz,\mt^k\leq\e} \mx^{-p}
\bigl\lvert \un_{\mt^k\leq\e}(t) - \un_{\mt^k\leq\e}(t^h)\bigr\rvert
\,\leq\,
2'\max\{\szme,\e\}
\end{align*}

Consequently we have estimated the difference
between the current case, III.4.1 (c), and
the previous one III.4.1 (b). By way of
notation, to emphasise this call the discontinuity
$w(\e)$ and interpret the integral over $w(\e)$
as an integral against an elementary function in $t$-
{\it i.e.} $\un_{w(\e)}$ is a difference
$\un_{w_+(\e)}-\un_{w_-(e)}$ of 
characteristic functions, where $w_*(\e)\sbs F_*(\e)$.
As such for $\o$ a 1-form of which we must
calculate the residue, this will be the
limit in $\e$ of,
$$
\int_{\substack{\mze\\ |w|\leq\e}}\o\dm \, +\,
\int_{\substack{|w|=\e\\ \mzle}}\o\dm \, +\,
\int_{\substack{w(\e)\\ |w|,\mz\leq\e}}\o\dm
$$
with the implicit restrictions $\mx\leq X$, respectively
$\my\leq Y$. As such, we have an extra term
not appearing in III.4.1 (b) that we've
seen how to bound, whence for $3'>3$,

\noindent{\bf III.4.5 Fact}
Everything as per III.4.3 then,
$$
\Biggl\lvert\int_{\max\{\mz, |w|\}=\e}\o\dm \, +\,
\int_{\substack{w(\e)\\ \max\{\mz, |w|\}\leq\e}}\o\dm
\Biggr\rvert\,\leq\,
3'\max\{\szme,\e\}
$$
Furthermore, while 
unlike III.4.3
there are estimates on the
mass of transversals- which we've basically
seen and we'll quantify in III.4.7-
one should be aware that they're
well short of those found in III.4.2,
so, again, there is need for caution
when going from local to global.

Before proceeding to the said mass
estimates, let us note that a small,
and necessary, variation is possible,
{\it viz:},

\noindent{\bf III.4.6 Remark} Globally
one cannot quite use the normal form,
and one must make a coordinate change
of the form,
$$x\mpo \phi x (1+ a(x^ky^l)),\,\,\,\,
y\mpo \theta y (1+ b(x^ky^l))$$
for $\phi,\theta$ constants and $a,b$
functions of a single variable. Such
changes do not quite preserve the normal
form, since the nilpotent part of the
plane field in III.4.1 (c) becomes,
$$
(x^ky^l)^r\left(\wa(x^ky^l)x\dx\, +\, \wb(x^ky^l)y\dy\right)
$$
again for $\wa, \wb$ functions of a single
variable. On the other hand this is an
irrelevant perturbation, {\it e.g} the
formal power series for $t$ has the 
form,
$$xy^{l/k}\exp\bigl( 
\sum_{m=r}^\infty (x^k y^l)^m P_m (\log y)\bigr)$$
where $P_m$ is a polynomial of
degree at most $\lfloor m/r\rfloor$,
and up to homotheties $x\mpo\phi x$,
$y\mpo\theta y$ the leading term
is exactly as before. In particular
such coordinate changes, together
with the implied change in $w$,
change absolutely nothing whether
in III.4.1 (c), or, easier, the
rational case of III.4.1 (b).

Finally let us give the mass bound
proper to the set up III.4.1 (c).
Specifically to keep ourselves in
the notation of \S III.3, let us
write the approximating holonomy as,
$$\z:=t^{-kr}\mpo\z-\g + o(\frac{1}{|\z|})$$
then we can bound strips,
$R< \imz< R+2\pi$, or
$C< \imz< C+2\pi$, should $\rez>R, C<R$
in the desired way, {\it i.e.}
$R^{-p/kr}$ times the Segre class.
Plainly this is inadequate to get
the same bound for the whole
transversal, but arguing as per
III.3.4 it does give,

\noindent{\bf III.4.7 Fact}
Quantifiers as per III.2.2 (c),
so, in particular $\my$ in some
interval $I$, then transversals
$T_\bfp$ satisfy the bound,
$$
\int_{\substack{\mz,\mt\leq\e\\ \bfp}}
\mt^{kr}\dm
\,\leq\,
B\e^p\max\{\sup_{\d\leq\e} \szmd,\e\}
$$
and similarly for the transversals
$\mz, |s|\leq\e$. Indeed, one could
achieve better quantification, but
as above it also applies to 
the rational case of
III.2.1(b), where one could use
this bound to remove the big
height condition.

\subsection*{III.5 The Beast}

While in a certain sense straightforward,
there is a risk of complication as one
passes along the singular locus and 
the number of eigenvalues jumps from
1 to 2. This occurs at,

\noindent{\bf III.5.1 Set Up} In a
neighbourhood $U$ of the singular
locus of a foliated 3-fold, after
sufficient blowing up we have 
exceptional divisors $E_1, E_2$
such that after completion in the
former, but not the latter, the
foliation has a formal generator,
$$z\dz\, +\, x^py\dy$$
for $x,y,z$ formal coordinates, 
$E_1$ given by $x=0$, $E_2$,
$y=0$, and $E=E_1+E_2$.

The point here is that even the
centre manifold need not be 
defined around the completion
in $E_2$, cf. \cite{uu} \S I.5,
\cite{mp1} II.2. In any case,
take convergent coordinates $x,y,z$
such that $x=0$, $y=0$ define
$E_1$, and $E_2$ respectively,
and the singular locus is
the $x$ and $y$ axis, with
generically 2 eigenvalues on
the former versus 1 on the latter.
We can, therefore, suppose that
there are functions $a,b,c,d$
and some large $N\in\bn$ such that,
$$
\pa=(z+ax^Ny)\dz + x^p(1+bx^N)\dy
+ (cy+dz)x^{N+1}\dx
$$
is a convergent generator. In
particular blowing up in the
2-eigenvalue component, we have
$\pa x\in (x^Ny)$ and otherwise
the singularity becomes log-flat
({\it i.e.} trivial from the 
point of view of computing 
residues). We can also improve the
order of vanishing of $\pa z-z$
along $E_2$ albeit at the price
of a loss of domain in $x$. 
Regardless, for $\o$ a 1-form of
which we must calculate the
residue, {\it i.e.} $\pa(\o)$
a function, against an invariant
measure without support in $E$,
and zero Segre class around the same,
consider the following strategy
for $\e,\d >0$,
$$
\int_{\substack{\max\{\mz,\mx\}=\e\\ \my\geq\d}}
\o\dm\,+\,
\int_{\substack{\un_{E_2}=\d\\ \max\{\mz,\mx\}\geq\e}}
\o\dm
$$
where $\un_{E_2}$ is a distance function
to $E_2$ coinciding with $\my$ locally.
The former of these two integrals divides
up as,
$$
\int_{\substack{\mze,\mxle\\ \my\geq\d}}
\o\dm\,+\,
\int_{\substack{\mxe,\mzle\\ \my\geq\d}}
\o\dm
$$
Of which, the first integrand is plainly
absolutely integrable even as $\d\ra 0$,
and the integral itself is bounded by
the Segre class. As to the second, we
may write,
$$\o=f \dfrac{dy}{x^py}$$
for $f$ some function, so that for
some $m$ to b chosen,
\begin{align*}
\int_0^\e \dfrac{dr}{r^{m+1}}
\Bigl\lvert \int_{\substack{\mx=r,\mz\leq r\\ \my\geq\d}}
\o\dm\Bigr\rvert
\,\leq\,
\int_{\substack{\mz,\mx\leq\e\\ \my\geq\d}}
\Bigl\lvert f\dfrac{d\mx}{\mx^{m+1}} \dfrac{dy}{x^py}  \Bigr\rvert\dm  
\end{align*}
while by the above, after blowing up,
$\pa x\in (x^{N+1}y)$, so up to a
constant the above right hand side
is at worst,
$$\e^{N-(m+p)} \int_{\substack{\mz,\mx\leq\e\\ \my\geq\d}}
dy\otimes d\bar{y}\,\dm
$$
which does not depend on $\d$. Whence,
keeping $\e>0$ we can let $\d\ra 0$
first. By the hypothesis on the
Segre class, there is no
residue around $E_2$ for $\e>0$, so
our residue is, in fact,
$$\lim_{\e \ra 0}\,\, \int_{\max\{\mz,\mx\}=\e}
\o\dm
$$
with no restriction on $\my$, and whence
zero for zero Segre class. Consequently,

\noindent{\bf III.5.2 Remark/Summary} 
Not withstanding the rather complicated
formal structure around $E_2$ the beast
poses, in the presence of zero Segre
class, it is no more or less of a problem
than any other point in the singular
locus where the induced foliation is
smooth. If, however, we wished to get 
quantification of a given residue in
terms of a possibly non-vanishing Segre
class, it would be a different story.
Fortunately, this isn't our situation,
and the only further point to remark
is that for globalisation, we'll have
to apply the above not for $x$ but for
$h(x)=\phi x (1+ O(x))$ some convergent
automorphism, but, plainly this is a
non-difficulty. As such, in practice,
the beast poses no problem, and it will
be passed over without comment, {\it i.e.}
it is to be understood that in collapsing
the boundary to the singular locus the
above strategy of first letting $\d\ra 0$
for $\e>0$ fixed, then taking $\e\ra 0$
is being observed.

\newpage

\section*{Globalisation}

\subsection*{IV.1 Warm up case}

As the title suggests we first do an
example in order to understand where
the difficulties lie, to wit:

\noindent{\bf IV.1.1 Set Up} Let
$U$ be a foliated 3-dimensional
tubular neighbourhood of a smooth
compact curve $Y$, such that:
\begin{itemize}
\item[(a)] $Y$ is the singular locus of the foliation.
\item[(b)] The foliation has canonical singularities
along $Y$, but only 1-eigenvalue at each point.
\item[(c)] In the formal centre manifold
$\wh{Z}$ obtained on completing in $Y$,
the induced foliation is smooth and $Y$
is invariant.
\item[(d)] A priori blowing up has been
performed so that for some smooth 
connected exceptional divisor $E$,
$Y=E\cap\wh{Z}$ and a hypothesised
invariant measure $\dm$ has neither
support on $E$ nor Segre class around the
same.
\end{itemize}

As it happens, one can show under these
hypothesis that $Y$ is an elliptic curve, 
but we'll eschew this so as to illustrate
some general features. Around every 
point in $Y$ one may (usual conventions)
find formal coordinates $x,y,z$ such
that the foliation is given by,
$$z\dz + x^p\dy$$
where formal means completion in $Y$,
or, indeed $E$, so the function $y$
is actually convergent, and, of course,
$p$ is fixed. Now there are two problems,
\begin{itemize}
\item[(i)] The existence domain for
a conjugation from analytic to such
a normal form while a bi-disc in
$(y,z)$ may be a sector of aperture
no more than $\pi/p$ in $x$, and this
is too small to be useful, \cite{mp1}
II.1.5.
\item[(ii)] Even if one had a conjugation
on a larger sector, there would still
be the problem of how these coordinates
patch as one moves whether in  
$Y$ or the argument of $x$, with the
latter being rather bad.
\end{itemize}
Let us therefore circumvent these difficulties
by an appropriate approximation procedure,
beginning with the formal centre manifold.
On an open cover $\coprod_\a U_\a$ $\ra U$, we have
functions $z_\a$ such that $z_\a=0$ is as
good an approximation, say modulo $\cO(-NE)$
to the centre manifold as we please, and so,
$$z_\a=g_{\a\b}z_\b\,\, \mod\, \cO(-NE)$$
where  $g_{\a\b}\in\H^1(U_N,\bg_m)$ and 
$U_N$ is defined via the exact sequence,
$$
\begin{CD}
0@>>>\cO_U(-NE) @>>>\cO_U@>>>\cO_{U_N}@>>>0
\end{CD}
$$
so for $\cA_*$ smooth functions we may 
profit from the preparation theorem 
to conclude that,
 $$
\begin{CD}
0@>>>\cA_U(-NE) @>>>\cA_U@>>>\cA\otimes_{\cO}\cO_{U_N}=:\cA_{U_N}@>>>0
\end{CD}
$$
is exact.
Furthermore,
$$
\begin{CD}
|g_{\a\b}|^2\in\H^1(U_N,\cA^{\times}_{U_N,+})@<{\sim}<{\exp}< \H^1(U_N,\cA_{U_N,\br})
\end{CD}
$$
where the subscripts $+$, $\br$ indicate positive real
valued, and real valued respectively. The above arrow
being an isomorphism, we may find real valued functions
$\phi_\a$ such that,
$$e^{-\phi_\a}|z_\a|^2-e^{-\phi_\b}|z_\b|^2
\,\in\,
\G(U_{\a\b}, \cA_U(-NE))$$
whence for some $\psi_\a$ vanishing
to order $N$ along $E$, we have:
$$
\nzs\,:=\, e^{-\phi_\a}|z_\a|^2 + \psi_\a\,=\,
 e^{-\phi_\b}|z_\b|^2 + \psi_\b \, 
\in \G (U, \cA_U)
$$
where, notation not withstanding, $\nzs$
might take negative values, albeit we only
care about its positive values so one could
take a max with zero if one prefers. For
convenience, we can arrange local generators
$\pa_\a$ of the foliation so that,
$$\pa_a (z_\a) \, =\, z_\a\,\,\, \mod\, \cO(-NE)$$
and whence for some $f_\a$ vanishing to 
order $N$ along $E$,
$$
\pa_\a\nzs\,=\, (1-\pa_\a\phi_\a)\nzs + e^{-\phi_\a} 2\Re (z_a\bar{f}_\a)
+ (\pa_\a \psi_\a -(1-\pa_\a\phi_\a)\psi_\a)
$$
Consequently if $\nzs=\e^2$ and the distance
to $E$ is also $\e$ and $N\gg p$, we have:
$$\dfrac{\pa_\a\nzs}{\nzs} \, =\, (1+O(\e^p))$$
and similarly for $\pab_\a$.

The discussion of approximation to the
exceptional divisor is a little more
complicated. In the first place as per
\S III.1 we have the holonomy representation
of the formal foliation in $\wh{Z}$ to
order $N$, {\it i.e.}
$$h:\pi_1(Y)\longra \Aut \bigl(
\dfrac{\bc[t]}{t^N} \bigr)$$
and for $g\geq 0$ the genus of $Y$
we choose a basis of simple closed
curves, $\s_1,\tau_1;\hdots;\s_g,\tau_g$
in the usual way, so that,
$$
\s_i.\tau_j=\d_{ij},\,\,\,\,
\s_i.\s_j=\tau_i.\tau_j=0
$$
Even though $Y$ is compact, its
technically useful (and eventually
necessary in the presence of
singularities of the induced
foliation in $\wh{Z}$, so the
current notation is temporary)
to view the $\s_i, \tau_i$ as
homology classes, with a dual
basis $\s^\vee_i(=\tau_i)$,
$\tau^\vee_i(=\s_i)$ in co-homology,
since the holonomy around the 
homology class gives rise to
a discontinuity in the 
(approximately) invariant function
describing $\cO_{\wh{Z}}(-E)$ along
the co-homology class. More 
precisely there is a formal function
$t$ (even on all of $\wt{Y}\times\Spf\bc[[t]]$,
$\wt{Y}$ the universal cover) such that
the co-equalisers satisfy,
$$t|_{\tau^+_+}=h(\s_i)(t|_{\tau^-_{i}})$$
where $\tau^+_i$, $\tau^-_i$ are the
respective sides of $\tau_i$, similarly
for $\s_i$ and $\tau_i$ interchanged,
and is continuous otherwise. A priori
this only holds in $\wh{Z}$ but blowing
up $N$-times in $Y$ it becomes true in
the completion $\wh{U}$ of $U$ in $Y$
modulo $\cO(-NE)$. Now for 
$\tu$ the universal cover of $U$,
view $\mt^2$ as a smooth section
of $\cA_{\tu_N}(E\otimes\bar{E})$,
then it may be lifted to a
global section of
$\cA_{\tu}(E\otimes\bar{E})$
which we denote by the same letter,
so that for $\pa_\a$ any local
generator of the foliation,
$$\pa_\a \mt^2\, \in \,\cA_\tu(-NE)+ \cA_\tu(-N\bar{E})$$
while on restricting to the fundamental
domain implicit in our basis for
the homology we have the same
defined on $U$ with a discontinuity
such that,
$$\mt^2|_{\tau^+_i}\,=\, |h(\s_i)(t|_{\tau^-_i})|^2
\,\,\,\mod \,\, \cO(-NE)$$
and, in addition, the discontinuity
not withstanding, $\mt^2$ is commensurable
in $U$ to an actual continuous distance
function $\un_E$ to the exceptional
divisor, {\it i.e.} for constants $c,C$
one has,
$$c\,\leq\, \dfrac{\mt^2}{\un^2_E}\,\leq\, C$$

We can now form a distance function
$\max\{\nzs,\mt^2\}$ which will be
discontinuous along real hypersurfaces
(denoted by the same letter) which
cut $Y$ in $\s_i$, respectively $\tau_i$.
As post III.4.4 this discontinuity will
be a difference of set functions, so we
put,
$$\un_{\tau_i(\e)}\, := \un_\e (t^{\s_i}) -\un_\e (t)$$
for $\un_\e$ the characteristic function
of $\mt^2\leq\e$, and understand by 
the integral ``over'' $\tau_i(\e)$
the integral against this elementary
function, and similarly for the $\s_i$,
and $\tau_i$ interchanged. Consequently
for $\o$ any 1-form of which we must
calculate the residue, it will be the
limit in $\e\ra 0$ of,
\begin{align*}
\int_{\substack{\nz=\e\\ \mtle}} \o\dm \, +\, 
\int_{\substack{\mte\\ \nz\leq\e}} \o\dm \, +\, 
\sum_{i=1}^g \biggl(
\int_{\substack{\s_i(\e)\\ \nz\leq \e}} \o\dm \, +\, 
\int_{\substack{\tau_i(\e)\\ \nz\leq \e}} \o\dm \biggr)
\end{align*}
Much of which will be handled exactly as
in the local cases, {\it i.e.} the first
term is plainly bounded by the Segre,
the second will prove negligible outside
a set of finite $\lem$ measure, and
we'll bound the final ones using the
almost holonomy estimate. The only really
new feature is that the holonomy group 
may be non-commutative.

To this end we have to consider the
structure of the representation
beginning with its first order
behaviour, {\it i.e.} the character:
$$
\begin{CD}
\chi:\pi_1(Y)@>>h> 
\Aut \bigl(
\dfrac{\bc[t]}{t^N} \bigr)
@>>>\bg_m
\end{CD}
$$
and, rather more importantly, $|\chi|$
in $\br^{\times}_+$ which we previously
encountered in \S III.2. The situation
where $|\chi|\neq 1$ is particularly
advantageous. Indeed $|\chi|(\s_i)\neq 1$,
and $\e$ sufficiently small imply that
the co-equaliser $\un_\tie$ is actually
the characteristic function of a set
(or maybe its negative). More precisely,
say, without loss of generality,
$\mch(\si)>1$ then $\un_\tie > 0$,
and dominates the characteristic
function of the annulus,
$$\{\e\,>\, \mt\,>\, q\e\}$$
for any $q>\mch (\si)^{-1}$, again
for $\e$ depending on $q$ sufficiently
small. This implies exactly as in
III.2.2 that for $\bfp$ a point of
intersection between $\si$ and its
dual (more correctly between a family
of perturbations of $\si$ which play
the role of the interval $I$ of op. cit)
we can use the basic construction
III.3.1  to estimate the measure
not just of $\tie$, but the whole
transversal $T_\bfp(\e)$ given by
$\mtle$ at $\bfp$. Better still
modifying the basic construction
of III.1.1 according to the variant
III.1.1 bis (post III.3.10) we can
for $\bfq$ in a given bounded smooth path 
$\rho$ get estimates for the mass
of all transversals through small
perturbations of it as follows.
Firstly join some small perturbations
$\bfp'\in P$ of $\bfp$ to small perturbations
$\bfq'\in Q$ of $\bfq\in\rho$ then for quantifiers
as per III.3.12 we'll find many
many variants of the basic construction
with mass on the $\mt=\e$ and
$\nz=\e$ faces bounded in the
fashion that we've become accustomed to.
Whence we get bounds on transversals
for $\bfq'\in Q$ on a set of measure
that is as close to full as we please.
Now perturb $\rho$ to $\rho'$ in a way that is
parametrised by $Q$. At this point
the bounds that we get on the
$\mt=\e$ and
$\nz=\e$ faces in applying the
basic construction to $\rho'$ 
are absolute and may be applied
equally to the variant III.1.1 bis
at all (actually Lebesgue almost all)
points of the perturbed path where
the anticipated bounds are achieved,
so, in summary:

\noindent{\bf IV.1.2 Fact} Suppose
$\mch\neq 1$ and
the order of approximation
sufficiently large, then for
$\e$ outside a set of finite $\lem$
measure, given a
smooth bounded path $\rho\sbs Y$
({\it e.g.} $\si$ or $\ti$) 
with $I$ a small interval in
the normal direction parametrising
perturbations $\rho'\in I$ there
is a constant $C$ such that for
$\rho'\in I$ outside a set (depending
on $\e$) of measure at most $B^{-1}$,
and $\bfq\in\rho'$ outside a set of
null Lebesgue measure, we have:
$$\int_{T_\bfq(\e)}\dm \,\leq\, \e^p BC\max\{\szme,\e\}$$ 
for $T_\bfq$ the transversal 
$\mt,\nz\leq \e$ at $\bfq$.

Plainly, this is way more than we
need and yields a bound on the
total variation of $\dm$ on
every $\tie$ or $\sie$ which is
much more than adequate. Consequently,
we can simplify the notation a
little by supposing that $\chi$
takes values in $S^1$.

At this point we must be more carefull
since as we've seen in \S III.3-4,
co-equalisers of maps tangent to the
identity do not afford an estimate
such as IV.1.2, while maps conjugate
to irrational rotations are only conjugate
to such. Consequently we distinguish two 
cases, beginning with the possibility that
the linear/first order holonomy is torsion.
We therefore make,

\noindent{\bf IV.1.3 Supposition} According to this case,
at the negligible cost of passing to a finite
\'etale cover, the formal holonomy of the
induced foliation in the formal centre
manifold is everywhere tangent to the identity.

Under this hypothesis, we're in the situation
already encountered in III.4.1 (c), {\it i.e}
the co-equaliser $\un_{\tie}$ contains both
positive and negative set functions, but the
co-equalisers in appropriate half planes
$\tie_+$, $\tie_-$ may be supposed strictly
positive, respectively negative, and,
$$\bigl\lvert \un_{\tie}(t) \bigr\rvert\,\leq\,
\bigl\lvert   \un_{\tie_+}(t^{\si}) - \un_{\tie_+}(t) \bigr\rvert
\, +\, \bigl\lvert   \un_{\tie_-}(t^{\si}) - \un_{\tie_-}(t) \bigr\rvert
$$
Now, while overkill from the point of
view of computing residues, we can argue
exactly as above to conclude,

\noindent{\bf IV.1.4 Fact} Suppose we're in the
situation of IV.1.3, and otherwise quantifiers
as per IV.1.2 albeit we only look at small normal
perturbations $\ti'\in I$ of the given $\ti$ (or for that
matter $\si$), then,
$$\int_{T_\bfq(\e)\cap\tie}\mmu \,\leq\,\e^p BC\max\{\szme,\e\}$$ 
where we write $\mmu$ to emphasise that the
left hand side is the total variation of
$\un_\tie$ over the transversal.

Again after passing to an \'etale cover,
it therefore remains to discuss,

\noindent{\bf IV.1.5 Remaining Case} $\chi$
is $S^1$ valued and the formal holonomy of
the induced foliation in the formal centre
manifold is non-trivial (otherwise there is
nothing to do) extension of the image of $\chi$
by a group of automorphisms tangent to the
identity.

The estimate will be somewhat akin to IV.1.2,
in a IV.1.3'sh sort of way. We normalise the
representation so that $r$ is the smallest 
integer for which we find a non-trivial map,
say,
$$t^a = t (1+\g t^r + O(t^{r+1}))^{-1/r}$$
while, as automorphisms of,
$\bc[t]/t^r$ everything commutes and is
identified to an element of $S^1$. This
implies that the total variation of any 
possible discontinuity
is strictly (understood in terms of
the relevant sets rather than their 
measure) less than that of the mass of
the annuli,
$$A(\e) := \{ \e^{-r} \,\leq \mt^{-r}\leq \e^{-r} + c\}$$
for some absolute constant $c$- basically
the maximum modulus of the $(r+1)$th
Taylor coefficient of a generator of
the representation, normalised as above.
By hypothesis, we also have irrational
rotations,
$$t^b = \lb t (1 + O(t^r))$$
for $\lb\in S^1$ non-torsion. Consequently,
$$bab^{-1}:t\mpo t(1 + \lb^{-r}\g t^r + O(t^{r+1}))$$
The supporting hyperplanes of the regions
$\tie_-$, $\tie_+$ encountered previously
for the action of $h(\si)$ may first be
considered for the action of the element
$a$, say, $a(\e)_-$, and $a(\e)_+$. The supporting
hyperplanes of these are perpendicular to
$\g$. Equally we have the same for the
conjugate element $bab^{-1}$, say,
$\tilde{a}(\e)_-$, and $\tilde{a}(\e)_+$,
but perpendicular to $\lb^{-r}\g$, so
together these cover an annulus,
$$ \e^{-r} \,\leq\, \mt^{-r}\,\leq \e^{-r} + c' $$
where $c' > 0$, but in general less than $c$.
However this is hardly a problem, since we have
every right to iterate the basic construction
III.3.1 to estimate co-equalisers- indeed, already 
neither $a$ nor its conjugate need have been simple.
Whence dividing the loops up into compositions 
of smooth (to avoid technical issues about
slicing) simple ones we estimate the co-equalisers
as anticipated, reason as per IV.1.2, and obtain,

\noindent{\bf IV.1.6 Fact} Replacing the hypothesis
$|\chi|\neq 1$ by that of IV.1.5, and otherwise
everything exactly as per IV.1.2, we have:
$$\int_{A_\bfq(\e)}\dm \,\leq\, \e^p BC\max\{\szme,\e\}$$ 
for $ A_\bfq(\e)$ any annulus of the form,
$$ \e^{-r} \,\leq\, \mt^{-r}\leq\, \e^{-r} + c $$
inside $T_q(\e)$ with $c$ independent of $\e$, and,
implicitly $C$ depending linearly on $c$.

At which point we may conclude,

\noindent{\bf IV.1.7 Fact} Suppose the set up
IV.1.1, let $\iy$ be the ideal of the singular
locus, and $\dm$ an invariant measure with
support outside the exceptional divisor and
zero Segre class around the same, then:
$$\RES:\Ext^2_U (\iy,\kxf)\longra\Ext^3(\cO/\iy,\kxf)$$
is zero on $\dm$.

\noindent{\bf proof} As already observed the
residue of a 1-form $\o$ that must be 
calculated is the limit as $\e\ra 0$ of,
\begin{align*}
\int_{\substack{\nz=\e\\ \mtle}} \o\dm \, +\, 
\int_{\substack{\mte\\ \nz\leq\e}} \o\dm \, +\, 
\sum_{i=1}^g \biggl(
\int_{\substack{\s_i(\e)\\ \nz\leq \e}} \o\dm \, +\, 
\int_{\substack{\tau_i(\e)\\ \nz\leq \e}} \o\dm \biggr)
\end{align*}
where there is no obligation to keep the
$\tie$ or $\sie$ fixed, {\it i.e.} it is
sufficient that we collapse down to $Y$
as $\e\ra 0$. Consequently, 
observing that we may safely replace
$U$ by any finite \'etale cover,
then perturbing
the $\sie$, $\tie$ in $\e$ as appropriate,
and appealing to IV.1.2, IV.1.4, or
IV.1.6 according to the circumstances,
we may handle all of these discontinuities.
The leading term is trivially bounded by
the Segre class, and the approximately invariant
term is easy, {\it i.e.} for local projections
$\pi_\a:U_\a\ra Y$, $\o$ may, via a partition
of unity $\rho_\a$, be written,
$$\sum_\a \rho_\a f_\a x_\a^{-p}\pi_\a^*\o_\a$$
where $f_\a$ is a function on $U_\a$, and
$\o_\a$ a holomorphic 1-form on $U_\a\cap Y$.
As ever, $x_\a$ is an equation for the exceptional
divisor, and for any local generator $\pa_\a$ of the
foliation $\pa_\a t$ is at worst 
$x_\a^{N-p}\o_\a$ for $N$ as large as we like.
Consequently, this term is trivially small
outside of a set of finite $\lem$ measure. $\Box$

\subsection*{IV.2 Warmer Case}

We treat the increase in difficulty by 
increments, and, whence, address:

\noindent{\bf IV.2.1 Set Up} Let
$U$ be a foliated 3-dimensional
tubular neighbourhood of a connected
compact curve $Y$ each component
of which is smooth, the singularities
are at worst plane nodes, and furthermore:
\begin{itemize}
\item[(a)] The curve $Y$ is the singular locus of the foliation.
\item[(b)] The foliation has canonical singularities
along $Y$, but only 1-eigenvalue at each point.
\item[(c)] In the formal centre manifold
$\wh{Z}$ obtained on completing in $Y$,
the components $Y_i$ invariant by the
induced foliation are disconnected.
\item[(d)] A priori blowing up has been
performed so that for some connected
exceptional divisor $E$, with smooth
components $E_k$, 
$E_k\mpo E_k\cap\wh{Z}$ 
is a 1-1 correspondence between 
components of $E$ and components of $Y$.
\item[(e)] The singularities that occur
are no worse than those of II.1.1, 
II.2.1 (a), III.2.1 (a), III.2.1 (b),
or III.3.1 (a).
\item[(f)]  A hypothesised
invariant measure $\dm$ has neither
support on $E$ nor Segre class around the
same.
\end{itemize}

Here one should observe that IV.1.2 (c) basically
implies IV.1.2 (e). The exceptions are rather
particular things such as II.3.1 with $q=0$,
or the big height condition of III.2.1 (b) 
not being satisfied. As such we're not quite
treating the general case implied by IV.1.2 (c),
but it's perfectly sufficient for our immediate
goal of understanding how to treat singularities.
In this context, the occurrence of any of,
II.1.1, II.2.1 (a), III.2.1 (a), or III.3.1 (a)
put us more or less immediately in the good
case IV.1.2, so, this is really only a discussion
about III.2.1 (b).

As such, to begin with, every $Y_i$ of IV.2.1 (c)
is either all of $Y$, or it meets a component
$Y'_j$ where the induced foliation is
everywhere transverse, so by II.1.5 the measure
can have no support in a neighbourhood of
$Y'_j$. Now all we need is some notation,
{\it i.e.} for,
$$z_i\frac{\pa}{\pa z_i} + x^{p_i}\frac{\pa}{\pa y_i}$$
the normal form at a general point of $Y_i$,
call $p_i=p(Y_i)$ the {\it multiplicity} along
$Y_i$, then IV.1.2 applies with the following
quantifiers,

\noindent{\bf IV.2.2 Fact} Suppose
that some $Y_i$ of IV.2.1 (c)
is not all of $Y$, and
the order of approximation
sufficiently large, then for
$\e$ outside a set of finite $\lem$
measure, given $U\sbs Y$ bounded
away from the singularities of
the induced foliation in $\wh{Z}$,
and a
smooth bounded path $\rho\sbs U$ 
with $I$ a small interval in
the normal direction parametrising
perturbations $\rho'\in I$ there
is a constant $C$ 
(so, a priori depending on $U$)
such that for
$\rho'\in I$ outside a set (depending
on $\e$) of measure at most $B^{-1}$,
and $\bfq\in\rho'$ outside a set of
null Lebesgue measure, we have:
$$\int_{T_\bfq(\e)}\dm \,\leq\, \e^{p_i} BC\max\{\szme,\e\}$$ 
for $T_\bfq$ the transversal 
$\mt,\nz\leq \e$ at $\bfq$.

We need some more notation,
say, $Y_i^*$ the complement
of $Y_i$ by small neighbourhoods
(to be specified) of the singularities
of the induced foliation in $\wh{Z}$,
and $Y'_k$ the components where
this induced foliation is 
everywhere transverse. If
$Y_i^*$ is all of $Y_i$ then
we're in the situation of 
IV.1.1, so, without loss
of generality $Y_i^*$ is 
non-compact. Consequently,
as we have warned, we change
notation, and for each $i$,
take $\s_{ij}$ to be a basis
of the homology, with $\tau_{ij}$
a dual basis in co-homology
formed as follows: for the
contribution from the genus
everything is a loop as per
\S IV.1, for each puncture
$\s_{ij}$ is a simple loop
around the same, and simple
paths from puncture to puncture
for the $\tau_{ij}$, albeit
for convenience we take a
slit from the puncture to one
one of the closed loops
arising from the genus should
there be only one puncture.
According to the various cases
specified by IV.2.1 (e) there
is, by IV.2.2, nothing else to
do as soon as $Y_i\neq Y$. 
Indeed on $Y_i^*$ we argue as
in \S IV.1, but in the easy
way when IV.1.2 applies, the
local computations around the
punctures have all been done
in \S II, III, and while there
could be a discontinuity between
the local coordinates employed
in these computations, and 
the global coordinate on $Y_i^*$
this is comfortably dealt with
by IV.2.2. Consequently, we
can lighten the notation by
supposing that $Y=Y_i$,
drop the $i$ from the notation,
and let $\sj$, $\tj$ be 
basis in homology and its dual
as above. The relevant holonomy
representation is that of 
the induced foliation in the
formal centre around $Y^*$,
which again we denote by $h$,
with $\chi$ its first order
part, and whence,

\noindent{\bf IV.2.2 bis. Fact}
Suppose either
$\mch\neq 1$, or there is a singularity
of type III.3.1 (a) rather
than the hypothesis $Y_i\neq Y$,
then IV.2.2 holds with the same
quantification. 

As a result, we may argue exactly
as above to reduce ourselves to the
case where  all the singularities
have the form III.2.1 (b), and
the holonomy around $Y^*$ is
$S^1$-valued. Consider, therefore,
for a suitable 1-form $\o$ the
form of the residue calculation
with $\e\ra 0$,
\begin{align*}
\int_{\substack{\nz=\e\\ \mtle}} \o\dm \, +\, 
\int_{\substack{\mte\\ \nz\leq\e}} \o\dm \, +\, 
\sum_i &
\int_{\substack{\tau_i(\e)\\ \nz\leq \e}} \o\dm\, + \\ 
\sum_{j} &  \biggl(
\int_{\substack{\sjes\\ \nz\leq \e}} \o\dm \, +\, 
{\mathrm{Local}}_j(\e, \o\dm) \biggr)
\end{align*}
where the first line are all terms around $Y^*$, and the 
new terms in $j$, are any possible
discontinuities around the punctures, interpreted
as before in a signed way, between the local
coordinates employed in \S III, and the global
approximately invariant function $t$ on $Y^*$,
while Local$_j$ simply indicates the local
strategy of \S III for calculating the residue.

Again, we distinguish between the cases of
$\chi$ having torsion or non-torsion image
in $S^1$. In the latter case, we have sub-cases:
the formal holonomy is linearisable to
sufficiently high (basically p) order,
or, it is not. In the former case there
is no issue since $\mt$ may be supposed
continuous. In the latter case we may
normalise the representation so that for
some $\lb,\nu\in S^1$ with the latter non-torsion,
we find an irrational rotation, $a$ by
$\nu$, and any other element has the form,
$$b(t)=\lb t (1 +\b t^r + O(t^{r+1}))$$
for $r\in\bn\cup\{\infty\}$, albeit there
is a minimal $r\in\bn$ for which $\b\neq 0$.
Continuing to denote this element by $b$,
we have,
$$c\,=\,[a,b]\,= t(1 + \b\lb^{-r}(\nu^r-1) + 
O(t^{r+1}))$$.
Consequently, we're in exactly the 
situation of IV.1.6, {\it i.e.} 
$c$ and $aca^{-1}$ are elements of
the formal holonomy, combinations
of which yield co-equalisers which,
in turn, afford an estimate of the
mass of any annulus of the shape
indicated in op. cit., and this is
plainly enough to deal with any
possible discontinuities along the
co-homology classes $\ti$ in $Y^*$.

The situation around the punctures
is, however, more subtle. The easy
case is when the eigenvalue at the
puncture is irrational, {\it i.e.}
it's not a singularity of the
form III.2.1 (b) with rational
eigenvalue. By the above the 
commutator of an irrational rotation
is just the group of all rotations.
Consequently, up to scaling the
global $t$ on $Y^*$ is (more
precisely modulo the order of
approximation) the function 
$xy^{-1/\lb}$ encountered in 
III.2.1 (b) {\it et sequel.},
or the holonomy around the puncture
in the $t$ variable has the form
$b(t)$ as above, albeit with $r$
possibly larger than the minimal
one, while the local variable
$xy^{-1/\lb}$ is a function of
$t$ of the same form, {\it i.e.}
the first non-linear term in
either case occurs to the same
order. As such the estimate \`a la
IV.1.6, with quantifiers as per
IV.2.2, for the mass of annuli:
$\e^{-r}\leq\mt^{-r}\leq\e^{-r} +$ const.,
for the minimal $r$ as above is
more than adequate to estimate
the total variation over any
$\sjes$.

We thus arrive to the principle
pre-occupation of this section 
which is singularities of the
form III.2.1 (b) with rational
eigenvalue. The completion $\wh{U}$
of $U$ in the singular locus 
admits a covering $\pi:\wh{V}\ra\wh{U}$
ramified only in the formal divisors
$y=0$ for $y$ a local equation in
the normal form of a singularity
of type III.2.1 (b). In most cases,
we can suppose that 
the ramification is only at the
singularities with rational eigenvalue.
The exceptional cases being $Y\iso\bp^1$
with at most 2 rational eigenvalues
with distinct numerators when it may
also be necessary to ramify in 
some irrational ones too. In any case,
call $R$ the set of ramification points
in $Y$, and observe that the holonomy
of the induced foliation in the 
centre manifold $\wh{W}$ in $\wh{V}$ has no
torsion. On the other hand the only
way to conserve canonical singularities
after such a covering is to ramify
in what may be purely formal divisors,
so this is not a trivial operation
like its smooth counterpart implicit
in IV.1.3/5. Nevertheless the singularities
in $\wh{W}$ with rational eigenvalue
only (big height condition) contribute
holonomy to a negligibly high order, so
for a basis 
$\sit, \tit$ of homology, respectively
co-homology, of the corresponding 
covering $\pi:\ty\ra Y$ of $Y$ the
discontinuity is confined to 
co-homology classes, and homology classes
around irrational punctures, say,
$\sjst$ for some indices $j$. Furthermore,
away from the ramification $\pi$ is 
\'etale, so, it extends convergently to
a covering $\pi:\tu\ra U_R$ where $U_R$
is the complement of $U$ by polydiscs
around the ramification points $R\sbs Y$.
Consequently if $s$ is the invariant formal
function on $\wh{V}$ continuous outwith
the $\tit$ and possibly some loops
$\sjst$ at irrational punctures, then we
have a smooth lifting to $\tu$ of the
same modulo as large a power of the
exceptional divisor as we please, which,
in a minor notational confusion, we
continue to denote by the same letter.

Now denote by $d$ the degree of $\pi$,
so $\pi_* |s|$ is commensurate to the
$d$th power of the distance to the
exceptional divisor, and for $\o$
a 1-form on $U$ we can do the residue
calculation as,
\begin{align*}
& \int_{\substack{\nz=\e, |\pi_*s|\leq \e^d\\ y\in Y^*}} \o\dm \, +\, 
\int_{\substack{|\pi_*s|= \e^d, \nz\leq\e\\ y\in Y^*}} \o\dm \, +\, 
\sum_i 
\int_{\substack{\tiet\\ \nz\leq \e}} \pi^*\o\dm \\
+\,& \sum_{j}\int_{\substack{\sjest\\ \nz\leq \e}} \pi^* \o\dm\, +\ \,
\sum_k 
\int_{\substack{\skes\\ \nz\leq \e}} \o\dm \, +\, 
\, \sum_{jk}{\mathrm {Local}}_j(\e, \o\dm) 
\end{align*}
where by way of notation we distinguish 
the co-equalisers $\sjest$, $\skes$ according
as $\pi$ is un-ramified, or otherwise at
the singularity which affords the homology
class, and ``local'' again just means the
local strategy. Now all the eigenvalues at
the $\sjest$'s are around singularities with
 irrational eigenvalue, so, we
know how to estimate the total variation here.
For the same reason, or arguing as per IV.1.3/4
if all the eigenvalues were rational, we have
the expected mass estimate also over the
$\tiet$. Amongst the $\skes$'s, again the
irrational ones pose no problem, which leaves
a possible discontinuity around punctures
with a rational eigenvalue between $\pi_* s$,
and the special coordinates of \S III.2.1 (b)
to discuss. To this end let $x,y,z$ be local
coordinates for the normal form at such a
point, then modulo some large power, $N$, of the
exceptional divisor $\pi_* s$ is invariant
and continuous, so for $k|l$ the eigenvalue,
it has, up to homothety, the form,
$$ (x^ky^l)^{d/k} (1 + a (x^ky^l))^d\,\,\, \mod\,\, \cO(-NE)$$
for some polynomial $a$ in a single variable.
As such the change of local coordinates,
$$x\mpo x (1+ a(x^ky^l))$$
together with rescaling in $y$ so that we
have a disc of radius 1, ensures that, we
change nothing in III.2.1 (b) and, at the
same time, avoid any discontinuity with
$\pi_* s$. We, of course, know how to
calculate all the local terms in $j$ and $k$, 
and whence,

\noindent{\bf IV.2.3 Fact} Let things be as in
IV.2.1 with $\iy$ the ideal of the singular
locus, and $\dm$ an invariant measure with
support outside the exceptional divisor and
zero Segre class around the same, then:
$$\RES:\Ext^2_U (\iy,\kxf)\longra\Ext^3(\cO/\iy,\kxf)$$
is zero on $\dm$.

\subsection*{IV.3 Penultimate Case}

We move towards a conclusion by way of
introducing the difficulties associated
with connected but not irreducible
components of the singular locus which
are left invariant by the induced foliation
in the formal centre manifold, {\it i.e.},

\noindent{\bf IV.2.1 Set Up} Let
$U$ be a foliated 3-dimensional
tubular neighbourhood of a connected
compact curve $Y$ each component
of which is smooth, the singularities
are at worst plane nodes, and furthermore:
\begin{itemize}
\item[(a)] The curve $Y$ is the singular locus of the foliation.
\item[(b)] The foliation has not just canonical singularities,
but are as convenient as possible in the sense
of \cite{mp1} V.1.8. In addition the locus
$Y_1$ where there is exactly 1-eigenvalue
(counted with multiplicity) will be
supposed of pure dimension 1.
\item[(c)] In the formal centre manifold
$\wh{Z}$ obtained on completing in $Y_1$,
the sub-curve $Y'$ formed by components
$Y'_i$ invariant by the induced foliation
has a dual graph without cycles.
\item[(d)] A priori blowing up has been
performed so that for some connected
exceptional divisor $E$, with smooth
components $E_k$, 
$E_k\mpo E_k\cap\wh{Z}$ (strictly
speaking $\wh{E}_k\mpo \wh{E}_k\cap\wh{Z}$
when the number of eigenvalues jumps 
from 1 to 2) is a 1-1 correspondence
between components of $E$ and components of $Y$.
\item[(e)] A singularity of the
induced foliation in $\wh{Z}$ at which
there is only one component of $Y$
is no worse than II.2.1 (a), III.2.1 (a)/(b) or
III.3.1 (a)- which, in fact, already follows
by (b) above.
\item[(f)]  A hypothesised
invariant measure $\dm$ has neither
support on $E$ nor Segre class around the
same.
\end{itemize}

From what we've already seen in \S IV.1/2, 
there is an increasing degree of nuisance
value associated with holonomy in the
centre manifold which is torsion to first
order, and we'll require a parenthesis to
examine this more carefully. To this end
let $C$ be a component of $Y'$, and $C^*$
the complement of $C$ by discs around any
singularities of the induced foliation in
$\wh{Z}$. We will suppose,

\noindent{\bf IV.3.2 Hypothesis} Notations
as above the formal holonomy around $C^*$,
so a fortiori around the excluded singularities,
is to first order with values in $S^1$,
while the punctures are not nodes as 
encountered in III.3.1 (a)/(b). 

Now we have various sub-cases to consider.
The easy one we have already largely seen, {\it viz:}

\noindent{\bf IV.3.3 (a) Easy Case} The
holonomy in $C^*$ contains a map conjugate
(formally) to an irrational rotation. There
may be several such, but we choose one, and
normalise the representation so that,
$t\mpo \lb t$, $\lb\in S^1$ non-torsion, for
some variable $t$ around the loop, say,
$a$, in question.

Now consider the holonomy $t\mpo b(t)$ around
any puncture, which, by hypothesis is 
conjugate to that given by one of the
normal forms III.3.1 (b), III.4.1 (b)/(c).
Now, say the conjugation is $\s b\s^{-1}$,
for $\s$ of the form,
$$t\mpo \nu t (1+ \b t^m + O(t^{m+1}))$$
and $\b\neq 0$. Should the normal form be
an irrational rotation with multiplier
$\rho$, we therefore have,
$$b(t)= \rho t ( 1+ \b (1-\rho ^m) t^m + O(t^{m+1}))$$
So, as already observed, the commutators 
$[a,h]$, $[a^2, h]$ afford maps which 
allow, as per IV.1.6, the appropriate
mass estimate on annuli, $\e^{-m}\leq \mt^{-m}\leq \e^{-m}+$ 
const. Of course $m$ may be infinite, but
that equates, up to homothety, to an identity between $t$
and the coordinates affording the normal
form around the puncture.

It may be, however, that the holonomy around
the puncture is torsion to first order, and
so the normal form of the singularity is as
per III.4.1 (c), or even just torsion in the
case of the second possibility in III.4.1 (b),
or torsion to all intents and purposes in
the rational case of III.3.1 (b). Irrespectively,
the normal form of the holonomy about the 
puncture is to leading order,
$$ n(s) =\rho s (1 + \g s^{kr} + O( s^{k(r+1)})) $$
where, now, $\rho$ is a $k$th root of unity, 
and possibly $r$ infinite in the pure torsion case,
with (at a minor risk of notational confusion)
$s$ the normalised local variable employed in
the local residue calculations of III.4.1 (b)/(c),
or its natural extrapolation to the rational
case of III.3.1 (b). As such, whenever the
above conjugation $\s$ between $t$ and $s$ has
$m\geq kr$, we find,
$$b(t) = \rho t ( 1+\g\nu^{kr} t^{kr} + O(t^{kr+1}))$$
and by way of commutators, we gain bound the
mass of transversals of an annulus,
$\e^{-kr}\leq \mt^{-kr}\leq \e^{-kr}+$const. in
an appropriate way. In the situation that $m<kr$,
we must distinguish cases. Should $k$ not divide
$m$, then we're akin to the previous irrational
case, {\it i.e.} for $b=\s^{-1}n\s$, the
commutators $[a,b]$, $[a^2,b]$ again yield a
bound on the annuli $\e^{-m}< \mt^{-m}<\e^{-m}+$ const.
Otherwise, $m=m_0k$, and, as per the end of the
previous section, we change the local coordinates
around the puncture by way of,
$$x\mpo \nu^{-1} x ( 1+ \a (x^ky^l)^{m_0}),\,\,\,\,\, y\mpo y$$
for $\s^{-1} t= \nu^{-1} t ( 1+\a t^{m_0k} + O( t^{m_0k+ 1}))$,
and $x$ the normal coordinate to $C$ in $\wh{Z}$.
This has no effect on the normal form III.4.1 (b),
a meaningless (big height condition again) effect
for III.3.1 (b), and 
a minor effect on the normal form in
the case of III.4.1 (c), whose negligibility has
already been noted in III.4.6. Consequently the
possibility that $k|m$ may be excluded. 
This equally applies, and indeed a particular
case was already encountered at the end of \S IV.2, 
to the extreme case that there is no first
order holonomy around $C^*$, {\it i.e.} we
may simply change the local variable of
III.3.1 (b), or III.4.1 (b)/(c)- in all cases
k=1- so that the local approximately invariant
function $s$ around the puncture agrees with the
global one $t$ (irrespective of its determination,
and, indeed to any order) without prejudice to
the local residue calculations of \S III. As
such, it remains to consider,

\noindent{\bf IV.3.3 (b) Fastidious Case} The first
order holonomy of $C^*$ in $\wh{Z}$ is torsion,
and non-trivial.

By hypothesis, therefore, the first order 
holonomy is cyclic of some order $d$, say,
generated by a $d$th root of unity $\lb$
on a variable $t$, normalised by way of,
$$a(t)= \lb t ( 1 + t^{nd} + O(t^{d(n+1)}))$$
for some maximal $n\in\bn\cup\{\infty\}$ amongst
all maps which are primitive $d$ torsion to first
order- so, in fact we put $n=\infty$ if we go
over our large a priori order of approximation
$N$, so, we can make the better normalisation, 
$$\lb t( 1 + A(t^{nd}))$$
for $A$ a function of a single variable.
Now take $t$ for the globally approximately
invariant function, and consider the conjugation
to the local variable $s$ of \S III by way of
$\s$ as above, so, plainly $k|d$. Again, if
$k|m$, or for that matter $m$ bigger than
our sufficiently large $N$, we can just change
the local variable $x$ so that $s$ and $t$
coincide without prejudice to the local strategies
of \S III. Thus, without loss of generality,
$k\nmid m$, and, a fortiori $\lb^m\neq 1$.
On the other hand,
$$a^{-1}(t)=\lb^{-1}t(1+ B(t^{nd}))$$
again for $B$ in a single variable.
Furthermore the holonomy, around the
puncture in the variable $s$ has the
form,
$$n(s) = \rho s (1 + f(s^k))$$
for $f$ in a single variable of sub-degree r-
the normal form is actually a function
of $s^{kr}$, but we may be in the
situation of III.4.6., so better than
the above cannot be guaranteed. This
implies that,
$$b=\s^{-1}n\s: t\mpo \rho t\{ 
1+\b (1-\rho^m) t^m + f(\nu^kt^k) + O(t^{m+1})\}$$
so, irrespectively of the relative sizes
of $m$ and $k$, the coefficient of $t^m$
does not vanish. To get a suitable mass
estimate (the possibility $kr<m<2kr$
is where the real nuisance value lies)
on transversal we consider the
commutator of $b$ with $a^{d/k}$,
$$[a^{d/k},b] =t\{ 1+\b(1-\rho^{-m})(1-\lb^{-md/k})t^m
+ g (t^k) + O (t^{m+1})\}$$
where, again, $g$ is polynomial in
a single variable, and, in addition
it vanishes to order strictly greater
than the maximum of $r$ and $dn/k$.
Consequently, by repeatedly forming
commutators with $a^{d/k}$ we eventually
obtain an element of the holonomy that
has the form,
$$c(t)=t\{ 1+\b(1-\rho^{-m})(1-\lb^{-md/k})^e t^m
 + O (t^{m+1})\}$$
for some $e\in\bn$. Now we have to
distinguish two cases. The easy one
is $\lb^{m}\neq -1$. In this case
the $(m+1)$th Taylor coefficients of
$c$, and $[a,c]$ (which is again tangent
to the identity to order $m$) are
linearly independent over $\br$, and
we can, once more, bound the mass of an
annulus, $\e^{-m}<\mt^{-m}<\e^{-m}+$ const.
in the usual way, which in turn is a
more than adequate bound for the
discontinuity between $t$ and $s$.
Otherwise, $\lb^m=\rho^m=-1$ and we
may only be able to bound the total
variation of mass between the regions,
$\mt^{-m}> \e^{-m}$ and $|c(t)|^{-m}>\e^{-m}$.
Fortunately, however, $\rho^m=-1$, so
up iterating this bound, {\it i.e.}
using $c^f$ for $f$ large (about 2 will do),
we bound the mass of the discontinuity
between $s$ and $t$. We have, therefore,
achieved the following dichotomy, 

\noindent{\bf IV.3.4 Summary} Suppose the holonomy
of the induced foliation in $\wh{Z}$ around
$C^*$ is to first order in $S^1$ and that
the multiplicity of the singularity around
$C$ is $p$. Then we may find an approximately
invariant function $t$ around $C^*$ such that
for any singularity of
the form III.4.1 (b),(c) and $\s$ a homology
class around the same either we can choose
the local coordinates of III.4.1 such that
the local invariant function $s$ of op. cit.
(actually $t$ in the notation therein) agrees
with $s$ or the total variation in the
discontinuity $\s^*(\e)$ between $t$ and $s$
around $\s$
(strictly speaking around a small perturbation $\s'$ of $\s$)
admits the bound,
$$\int_{T_\bfq\cap\s^*(\e)}\mmu \,\leq\, \e^{p} BC\max\{\szme,\e\}$$ 
for $\bfq\in\s$ and quantifiers, {\it etc.},
as per IV.2.2.

Which doesn't cover the case of III.2.1 (b),
but since $s=xy^{-1/\lb}$, $\lb\in\br_{<0}$,
and the local strategy there is to use the
boundary $\mxe$, since the invariant manifold
$y=0$ may only exist formally, we can just
rescale so that $\my=1$ 
on the boundary, and
the relevant discontinuity there,
{\it i.e.} between $\mx$ and $\mt$, is, in fact,
that between $|s|$ and $\mt$. Whence,

\noindent{\bf IV.3.4 bis. Similarly} Exactly
as per IV.3.4, but for singularities of
type III.2.1(b) and the discontinuity understood
as that between $\mt$ and the local variable $\mx$.

We can also clear up the relation between
the discontinuity along cohomology classes
$\tau$ occasioned by the dual homology class
$\s$. Again, we're in the situation of IV.3.3 (a) or (b),
with holonomy around $\s$ of the form,
$$h(\s) (t) = \rho t \{1 +\a t^r + O(t^{r+1})\}$$
with $t$ normalised as above. Whence in
the easy case 4.3.3 (a), as previously
observed in \S IV.1/2, we always get bounds
on the mass of an annulus,
$\e^{-r}<\mt^{-r}<\e^{-r}+$ const., and, idem
for IV.3.3 (b) if $\lb^r\neq \pm 1$. Should
$\lb^r=-1$, we bound the mass of a smaller
region, but as above this is the mass of the
discontinuity anyway, so we're okay. Finally
if $\lb^r=1$ then, a fortiori, $\rho^r=1$,
and again the mass is appropriately bounded,
so that without recourse to finite coverings
to kill the torsion,

\noindent{\bf IV.3.5 Variation} Everything as per IV.3.4,
with $\ti$ a dual basis in co-homology to
the homology basis $\s_i$ for $C^*$, but here
we bound
the total variation in mass over $\tie$,
or, more correctly a small perturbation of it.

Now let us organise the global residue
computation. By II.1.5 and III.5.2 there
is nothing to do except at components 
$C$ of the singular locus where (counted
with multiplicity) there is everywhere one
eigenvalue (counted with multiplicity),
the foliation is not log-flat, and $C$ is
invariant by the induced foliation in the
formal centre manifold $\wh{Z}$. Now as
per \S I.4, form the dual graph $G$ with
vertices such components $C$, and edges
intersections of the same. By the hypothesis
IV.3.1 (c), $G$ is a tree, and, we do not
prejudice the residue calculation by
supposing that it is connected. As such,
we may choose a root $R$, say, which in
turn defines a unique path to each vertex,
whence (depending on $R$) a direction on
each edge according to the sense of
increasing distance from $R$. At a
directed edge $e$ we introduce local
coordinates $x_e$, $y_e$, and a multiplier
$\nu_e$, {\it i.e.},
$$
\begin{CD}
v_-:(x_e=0)@>{\nu_e}>>v_+: (y_e=0)
\end{CD}
$$ 
according to the rule: $\nu_e=p(v_-)/p(v_+)$,
for $p(v)$ the multiplicity along the vertex
$v$ (so $q/p$ in the notation of \S III) 
whenever the singularity at $e$ is not
of the form III.4.1 (b) or (c), and otherwise
according to the form of the first order
part of the singularity in $\wh{Z}$, {\it i.e.,}
$$y_e\dfrac{\pa}{\pa y_e} - \lb_e x_e 
\dfrac{\pa}{\pa x_e}$$
$\lb_e\in\br_+$, and irrespectively of whether
its rational or irrational
we put $\nu_e=\lb_e$. 
We continue to
denote by $C^*$ the curve $C$ minus discs
about singularities in $\wh{Z}$, but we
also introduce a curve $C^{\bullet}$ where
we re-fill $C^*$ at singularities of the 
form  III.4.1 (b)/(c) or III.2.1 (b) where
the dichotomy of 
IV.3.4 \& bis implies that the
normalised approximately invariant function
$t_C$ extends  around the 
singularity by way of agreement with its
local equivalent in III.4.1 (b)/(c) or its
modulus with that of $\mx$ in III.2.1 (b).
As such, continuing in the abuse of notation
between $t_C$ and a smooth lift of the same
modulo a large unspecified power of the
exceptional divisor, we extend the abuse
to all of $C^{\bullet}$ by way of the above
rule. Finally if $\{R,C\}$ is the path from
$R$ to $C$ we introduce a multiplier $\lb_C$
according to $\lb_R=1$, and, otherwise:
$$\lb_C=\prod_{e\in\{R,C\}} \lb_e$$
together with a sub-graph $G^{\bullet}$
whose edges are intersections of 
filled curves $C^{\bullet}$. As such, by the
definition of $|w|$ in III.4.1 (c), at
such an edge, $|t_{v_+}|=|t_{v_-}|^{\nu_e}$.

The function $\mt:= |t_C|^{1/\lb_C}$ has,
therefore, on the curve defined by $G$ a discontinuity,
$\tau(\e)$,
not only at co-homology classes $\tau$ in
any $C^*$ for every vertex, but also
(bearing in mind we take a slit at the
puncture if the corresponding vertex has
only one singularity) the extension of these
up to, and around (it continues naturally
from one component to another by way of
$w(\e)$ post III.4.4)
a singularity when this has the
form III.4.1 (c) which is an edge in $G^{\bullet}$. In this latter case,
the mass bound comes from 
IV.2.2 \& bis, or
IV.3.4 
when the former does not apply, 
for the
part of the discontinuity in $C^*$ and
always
from III.4.4 close to the singularity.
Otherwise at an edge of $G$ 
not in $G^{\bullet}$
we have one,
or two loops, $\s^*$, according as the
singularity is filled in one or neither
of the vertices. Certainly there may be
a discontinuity between $|t_C|$ and the
boundary employed in the appropriate local
strategy of \S III at $\s^*(\e)$, or, indeed, any other
$\s^*$ which goes around a puncture in $C^*$
that is not filled in $C^{\bullet}$,
but in all such cases we may either 
trivially bound it by
IV.2.2 \& bis, 
and the evident extension of the same
to the singularities III.3.1 (c) and
III.4.1 (a),
or, when necessary by the more delicate IV.3.4. 

Now let's look at the form of the residue
calculation for an appropriate 1-form $\o$,
and $\e\ra 0$, {\it i.e},
\begin{align*}
&\sum_C \biggl( \int_{\substack{\nz=\e, |t_C|\leq \e^{\lb_C}\\{\mathrm{nhd.}} 
\,\mathrm{of}\, C^{*}}}\o\dm
\, +\,  \int_{\substack{\nz\leq\e, |t_C|= \e^{\lb_C}\\ \mathrm{nhd} \,
\mathrm{of}\, C^{*}}}\o\dm
\biggr) \, +\, \sum_i \int_{\tie} \o\dm \\
& \sum_j \int_{\sjes} \o\dm \, +\,
\sum_k {\mathrm{Local}}_k (\e, \o\dm)
\end{align*} 
where, by way of notation, all local strategies
at singularities,
irrespective of whether it's at a singularity 
which is an edge of $G^{\bullet}$ or not, have
been lumped together even if its more correct
to think of such edges globally since at them
the local and global strategies coincide. The
discontinuities $\tie$, $\sjes$ have been defined
above. and observed to be bounded in an appropriate
way, while the first two integrals are our friends
the Segre class, and something negligible. Whence:

\noindent{\bf IV.3.6 Fact} Let things be as in
IV.3.1 with $\iy$ the ideal of the singular
locus, and $\dm$ an invariant measure with
support outside the exceptional divisor and
zero Segre class around the same, then:
$$\RES:\Ext^2_U (\iy,\kxf)\longra\Ext^3(\cO/\iy,\kxf)$$
is zero on $\dm$.

\subsection*{IV.4 General Case}

We now bring our calculation to a
conclusion, {\it viz:}

\noindent{\bf IV.4.1 Set Up} Everything as
per IV.3.1, but without the no-cycles hypothesis
IV.3.1 (c), and suppose (pro tempore) that,
in addition, the neighbourhood $U$ is
projective, {\it i.e.} admits an embedding
in $\bp^n$ for some $n\in\bn$.

In the notation of the previous section,
plainly the problem lies in the multipliers
$\nu_e$, {\it i.e.} their product around
a cycle may not be 1. To this end, consider
a singularity of the form III.4.1 (b), but
with $\lb\notin\bq$, and let us explore
an alternative strategy for effecting the
local residue computation. Specifically,
let $\a,\b\in\br_{>0}$ and put $|f|=\mx^\a\my^\b$,
with the notations of III.4.1 (b), with
$g=x^ay^b$, $a,b\in\bz_{\geq 0}$. Observe
that up to an irrelevant error given a
sufficiently high order of approximation,
$$x^{-a}y^{-b}\dfrac{dx}{x} \,=\,
(a+b\lb)^{-1} \dfrac{dg}{g^2} \dm$$
As such, since we require to calculate
the residue of the left hand side for
$a\leq p$, $b\leq q$, it will be 
sufficient to do the right hand side
for all such $g$ under the same hypothesis.
On the other hand for $X,Y$ some upper bounds
for $\mx,\my$,
\begin{align*}
-\int_{\substack{|f|=\e,\mzle\\ \mx\leq X,\my\leq Y}}
\dfrac{dg}{g^2}\dm\,=\,
\int_{\substack{|f|=\e,\mze\\ \mx\leq X,\my\leq Y}}
\dfrac{1}{g}\dm\,
+\,\int_{\substack{|f|=\e,\mzle\\ \mx= X}}
\dfrac{1}{g}\dm \, +\,
\int_{\substack{|f|=\e,\mzle\\ \my= Y}}
\dfrac{1}{g}\dm
\end{align*}
The first of these integrals may be
estimated in the usual way, {\it viz}:
\begin{align*}
\int_0^1 \dfrac{d\chi}{\chi}
\int_{\substack{|f|=\e,\mze\\ \mx\leq X,\my\leq Y}}
\dfrac{1}{|g|}\mmu\,=\,
\int_{\substack{\mze, |f|\leq\e\\ \mx\leq X,\my\leq Y}}
\Biggl\lvert\dfrac{1}{|g|}\dfrac{d|f|}{|f|}\Biggr\rvert\dm\,\leq\,
\szme
\end{align*}
while the other two integrals may be seen
to be negligible exactly as per III.2.4 (b), so:

\noindent{\bf IV.2.2 Fact} Let things be as in
III.4.1 (b) with an irrational eigenvalue, 
and suppose that the Segre class around
the exceptional divisor of our invariant
measure is zero, then for $\e>0$ outside
a set of finite $\lem$ measure, and $(X,Y,\log\chi)$
varying in a set (depending on $\e$) which
as close to full with respect to Lebesgue
measure as we please,
$$\lim_{\e\ra 0}\,
\int_{\substack{|f|=\chi^{-1}\e,\mzle\\ \mx\leq X,\my\leq Y}}
x^{-a}y^{-b}\dfrac{dx}{x}\dm\,=\,0,\,\,\, a\leq p,\, b\leq q$$

One should be aware that there is a hidden
subtlety here. More precisely, from the
perspective of switching to a local strategy
where one calculates all residues by Stokes,
rather than simply considerations of size on 
the boundary,
it is a priori required
that we can calculate all residues of
the form,
$$\dfrac{f(x,y)}{x^py^q}\dfrac{dx}{x}$$
and the possibility of doing this does
not exactly follow from IV.2.2, since
if one expands $f(x,y)$ in a Taylor
series, and applies IV.2.2 there could
be a small divisors issue. On the other
hand, one only has to do the residues
of IV.2.2 for $1\leq a\leq p$, and
$1\leq b\leq q$ together with,
$$x^{-a} f(y) \dfrac{dx}{x},\,\,\, 1\leq a\leq p$$
for any convergent $f$, and similarly
with $x$, $y$ interchanged. The relevant
rational combinations of $1$, $\lb$ which
might be a source of concern are, therefore,
$$ -a + n\lb,\,\,\, 1\leq a\leq p,\,\, n\in\bz_{\geq 0}$$
while $\Re(\lb) \leq 0$, so the modulus
of the combination in question is at
least 1, and there is no issue. 
Consequently, let us note where the real problem is
by way of,

\noindent{\bf IV.4.3 (a) Remark}
Evidently the same strategy works in the
alternative rational case of III.4.1 (b)
with $\lb=-k/l$ whenever $a+b\lb\neq 0$,
{\it i.e.} $x^ay^b$ not a power of the
approximately invariant function $x^ky^l$.
In such good cases, we can employ arbitrary
multipliers at the edges of our graph, and 
there will be no possibility of a continuity
problem occasioned by cycles in the graph.
In the situation that $x^ay^b=(x^ky^l)^d$,
some $d\in\bn$, there remains, however,
a risk of a problem since the only choice
that works is $k/l$.

In light of this improvement, let as
look again at the singularities III.3.1 (c),
and III.4.1 (a). In the latter case,
the required changes to argue as above
are negligible, but in the former case
there are several problems resulting
from the fact that the normal form doesn't
converge on completion in the singular locus,
albeit these can be dealt with as 
pre III.3.11, and a much thornier
problem about getting the right estimates
for a residue of the form $x^{-a}dx$.
Consequently, let us adopt an expedient
that is adequate for our purposes, and
consider only the possibility that,
$\a\leq p$, and $\b\leq q$, then,
supposing zero Segre class,
$$\oint_{\substack{|f|=\e\\ \mx\leq X, \my\leq Y}}\dfrac{df}{f}\dm \,=\,\e o(\e)$$
with the same proof, {\it i.e.} III.3.13,
\& III.4.2 respectively. This also holds
for $\a$ or $\b$ zero, and, in fact with
the better bound $\e^{p+r/\a}$ in the
case of nodes. Now such a bound may not
be adequate to compute any residue, but
supposing $\a$, $\b$ integers, it does suffice to deduce,

\noindent{\bf IV.4.4 (a) Fact} Say the edge is a node
or a singularity of the form III.4.1 (a), II.2.1 (b)
being trivial, then for $\o\in\kf$ vanishing to order
$p-\a$ along $E_1$, $q-\b$ along $E_2$ and $\e$
outside a set of finite $\lem$,
$$\lim_{\e\ra 0} \int_{\substack{|f|=\e\\ \mx\leq X, \my\leq Y}}\o\dm\,=\,0$$

As such, we may return to the difficulty
of rational resonances, and the last
remaining possibility III.4.1 (c). Here
the strategy of how to apply Stokes has
no essential difference with IV.2.2/IV.4.3 (a)
beyond the estimation at the
ends, which are no longer negligible. Indeed for 
$g=x^ay^b$, and $la\neq kb$ at the
end $\mx=X$ we find the mass bound,
$$
\int_{1/2}^1\int_{1/2}^1
\dfrac{d\chi}{\chi}\dfrac{d\phi}{\phi}
\int_{\substack{|f|=\chi\e\\ \my=\phi X\\ \mzle}}
\mmu \,\leq\,
\dfrac{l\tilde{i}-k\tilde{j}}{2}
\int_{\substack{\mx\leq \e^{1/\a}Y^{-\b}2^{\b/\a}\\
X/2\leq\my\leq X\\ \mzle}}
\Im(h)\dfrac{dy\otimes d\bar{y}}{\my^2}\dm
$$
where the notation is as per the normal
form III.4.1 (c), or adjusted as per
III.4.6, from which one finds that $h$
has the form,
$$h \,=\, (x^ky^l)^r + O( (x^ky^l)^{r+1})$$
thus the mass bound III.4.7 is exactly what
we require, so, again,

\noindent{\bf IV.4.4 (b) Fact} If the edge
has the form III.4.1 (c) and $ka-lb\neq 0$,
then IV.4.2 continues to hold.

It is instructive to observe,

\noindent{\bf IV.4.3 (b) Remark}
While this situation is a bit
better than the rather absolute
obstruction posed by the rational
case of III.4.1 (b), {\it i.e.}
for $w=x^ky^l$ we have up to an
irrelevant error, and irrelevant
homotheties,
$$w^{-m}\dfrac{dx}{x}\dm\,=\,
\biggl( \dfrac{dw}{w^{r+m+1}}+\nu \dfrac{dw}{w^{m+1}}\biggr)
\dm$$
for $1\leq m\leq n$, $n$ as per
III.4.1 (c). Integrating over the
boundary $|f|=\e$ as above, the
latter term is no problem since,
$\mw^m\geq\mx^p\my^q$, and the
bound III.4.7 on transversals is
adequate. For the initial term,
however, it is not, since this
needs the stronger bound of 
IV.2.2, which we may, or may 
not, have, and would, in any
case, be needed on both 
transversal $\mx$= const., and
$\my=$ const.

To improve the situation we
need to do a little surgery.
This will only be at points
where the induced foliation
in $\wh{Z}$ is smooth, so introduce
formal coordinates $x,y,z$ such
that the foliation is given by,
$$\pa \,=\, z\dz + x^p\dy$$
We wish to perform the weighted
blow up $\pi:\tu\ra U$ of
$(z,y,x^p)$ in the sense of 
champ de Deligne-Mumford so
as to preserve smoothness.
The only relevant \'etale patch of
the blow up has coordinates,
$\xi,\eta,\z$, where:
$$\eta^p=y,\, x=\eta\xi,\, \z=\eta^p\z,
\,\,\,\,
\eta\mpo \theta\eta,\,
\xi\mpo\theta^{-1}\xi,
\z\mpo\z
$$
and we have an open inclusion
in $\tu$ of the classifying
champ $[\D^3/\mu_p]$ for the
implied action of $\mu_p$ by the
$p$th root of unity $\theta$. 
As such the foliation is given by,
$$\pa\, =\, (1-\xi ^p)\z\dfrac{\pa}{\pa\z}
+\dfrac{\xi^p}{p}
\bigl(\eta\dfrac{\pa}{\pa\eta}-\xi\dfrac{\pa}{\pa\xi}\bigr)
$$
while on the other patches it is 
smooth or log flat, so absolutely
irrelevant for computing residues.
As to the measure itself this extends
(thanks to the finiteness of the
weighted Lelong number) by zero over 
the exceptional divisor
in a way that has finite mass with
respect to a smooth metric on $\tu$,
so it remains an invariant measure
on the same. It is, nevertheless,
true that we've replaced a smooth
point of the foliation in $\wh{Z}$
by a non-schematic singular one. 
This has, however, no adverse effect
on the residue calculation: we work
on a tubular neighbourhood $V$ of
the proper transform of $Y$, the
singularity in the formal centre
manifold in $V$ has no extra holonomy
since $\eta\xi\mpo \eta\xi$ under
$\mu_p$, and we can treat this singularity
as an extremely easy case of III.2.1 (b).

As to why we would want to do this:
observe that the exceptional divisor
$E$ defining the components gets
replaced by its proper transform
$\cE$ around some not wholly schematic
component $\cC$ of the singular locus,
so in the proper transform $\cZ$ of
the centre manifold,
$$N_{\cC/\cZ}\,=\,\cO_\cC (\pi^* E - F)\,=\,
\pi^* N_{C/\wh{Z}}(-F)$$
for $F$ the exceptional divisor.
As such, if we repeat this operation
in  enough points, the resulting
$\cZ$ will be (formally) convex.
Now let's bring the projectivity into
play, and take some very large multiple
$nH$, to be specified, of a very
ample divisor such that $D\in |nH|$
cuts $Y$ in some reduced set $B$ of
smooth points with tangency to order
$m_b$ to the divisor $y=0$ understood
in the above local coordinates around
$b\in B$. Consequently for $F_b$ the
exceptional divisor over $b\in B$,
the weighted blow up having been 
performed in all of these, and $\pi:V\ra U$
again a tubular neighbourhood around
the proper transform of $Y$,
$$n\pi^* H \,=\, \pi^* D\,=\, 
\sum_{b\in B} p(C) m_b F_b$$
for $p(C)$ the multiplicity of the
component $C$ through $b$. The 
good choice of $n$ is, therefore,
$$n\prod_C p(C)$$
and plainly $m_b=p(C)^{-1}n$, so
that:
$$\pi^* H -\sum_{b\in B} F_b \,=\, \cL\in\Pic(V)$$
is at worst $n$-torsion, and has
monodromy around each non-schematic
point $\mu_{p(C)}$ for $C$ the 
corresponding component. From which
the \'etale cover $\rho:\tu\ra V$
defined by $\cL$ is everywhere schematic,
and for $\tz$ the centre manifold
in $\tu$ with $\tc$ some component
over $C$,
$$N_{\tc/\tz} =\rho^*\pi^* N_{C/Z}(-H)$$
and, as it happens, $\rho^*\pi^* Z=
\tz + p(C)^{-1}n H$ on each component,
so $\tu$ is highly convex.

We may loose IV.3.1 (d), but in a
reasonably trivial way, {\it i.e.}
for $\cE$ an exceptional divisor
in $V$ defining a component $\rho^*\cE$
will still be smooth, but possibly
disconnected. In principle this is
a minor caveat, and in any case for
$E$ a component of some $\rho^*\cE$,
the multiplicity $p(E)$ is well
defined by way of $p(C)$ for $C$
the original component in $U$ 
corresponding to $\cE$. Next let
$S$ be the  singular
locus in $\tu$ with scheme
structure around the
components invariant by the
induced foliation in the centre
manifold. We have a total exceptional
divisor,
$$E_{{\mathrm{tot}}}\,=\, \sum p(E) E$$
and various divisors $E$ between $0$
and $E_{{\mathrm{tot}}}$ partially
ordered by increasing multiplicity.
In particlar, there are maps,
$$\kf\otimes\cO_{E'|_S}(-E)\longra
\kf\otimes\cO_{E''|_S}(-E)$$
for $E_{{\mathrm{red}}}\leq E''\leq E'\leq E_{{\mathrm{tot}}}$,
$E\geq 0$ any,  and the same canonical for
all of $U$, $V$, $\tu$, which we
emphasise by omitting $\pi^*$,
$\rho^*$ on the same. In addition,
as per I.1.5, the trace form yields
an isomorphism between $\cO_S$ and
$\kf|_S$, so if $E'=E'' +F$ we have
an exact sequence,
$$
0\ra \cO_{F|\tz}(-E''-E) \ra
\kf\otimes\cO_{E'|_S}(-E) \ra
\kf\otimes\cO_{E''|_S}(-E)
\ra 0
$$
The map $\rho\pi$ is acyclic
on the reduced components of
$S$, so taking $H$ sufficiently
large and inducting on the partial
ordering, we obtain,

\noindent{\bf IV.4.5 Lemma} In the
above notation all of the following
maps are surjections,
$$
\G( \kf\otimes\cO_{S}(-E)) \ra
\G(\kf\otimes\cO_{E'|_S}(-E))\ra
\G(\kf\otimes\cO_{E''|_S}(-E)) 
$$

\noindent{\bf proof} It only remains to
observe that the kernel of,
$$\cO_S\longra \cO_{E_{\mathrm{tot}}|_S}$$
is supported in dimension zero. $\Box$

Now applying the same considerations in
$V$ with $\cS$ the singular sous champ,
$\cE$ an exceptional divisor, {\it etc.},
we have exact sequences,
$$
\begin{CD}
0\ra \G( \kf\otimes\cO_{S}(-E_{{\mathrm{red}}}))@. \ra
\G(\kf\otimes\cO_{S})@. \ra
\G(\kf\otimes\cO_{E_{{\mathrm{red}}}|_S})=\bc
@.\ra 0\\ 
@AA{\rho^*}A @AA{\rho^*}A @| \\
0\ra \G( \kf\otimes\cO_{\cS}(-\cE_{{\mathrm{red}}}))@. \ra
\G(\kf\otimes\cO_{\cS})@. \ra
\G(\kf\otimes\cO_{\cE_{{\mathrm{red}}}|_\cS})=\bc
@.\ra 0 
\end{CD}
$$
such that the group in the top left
hand corner has lots of sections.
Nevertheless one should be cautious
since although its generated by
global sections, it may not, for
example, separate points, albeit
it separates pull-backs of points,
and similar,
from $\cS$, and this will be enough.

More precisely let $\o$ in the
bottom right group be given.
At each singularity in the
induced foliation it has a 
Taylor expansion. Apart from
III.4.1 (c), the normal forms
converge after completion in
the singular locus, and we're
only ever interested in this
modulo a large power of the
exceptional divisor on restriction
to the centre manifold, so the
expansion is, in fact, convergent
in such coordinates. In the cases
of interest, III.4.1 (b)/(c),
these coordinates are actually
unique up to homothety in the
irrational case, and otherwise
the unimportant modification III.4.6.
Let us look at this expansion more
closely in these cases writing it
as a Laurent expansion,
\begin{align*}
\o\,=\,& \sum_{\substack{1\leq a\leq p\\ 1\leq b\leq q}} w_{ab}x^{-a}y^{-b} \dfrac{dx}{x}
\,+\, \sum_{\substack{1\leq a\leq p\\ b\geq 0}} w_{ab}x^{-a}y^{b} \dfrac{dx}{x}\\
+ & \sum_{\substack{a\geq 0 \\ 1\leq b\leq q}} w_{ab}x^{a}y^{-b} \dfrac{dx}{x}
\,+\, \sum_{\substack{a,b \geq 0}} w_{ab}x^{a}y^{b} \dfrac{dx}{x}
\end{align*}
Each of these four regions has its
own structure. The final one is rather
trivial from the point of view of 
using Stokes for a residue calculation
along $|f|=\mx^\a\my^\b$, since,
$f^{-1}\pa(f)=(\a + \lb \b)x^{-1}\pa x$,
so everything here is trivially bounded
by the Segre class provided $\a + \lb \b\neq 0$,
and otherwise we won't be using Stokes, and,
will be keeping to  the previous strategy,
{\it i.e.} using the
approximately invariant perturbation of $|f|$.
The two middle regions need not be
in the image of $\pa$, but, {\it cf.}
post IV.2.2, they are modulo the good
region already discussed. Consequently,
the obstruction to being able to write
$\o$ as $d$ of meromorphic with a suitable
pole (IV.4.3 (b) being an example of what
is unsuitable), and whence apply Stokes
on a boundary of our choice, is wholly
in the leading region, and it is finite
dimensional. On the further hypothesis
that $\o$, viewed as a section of $\kf$,
vanishes on $\cE_{{\mathrm{red}}}$, {\it i.e.}
belongs to the group on bottom left of the diagram,
we can find $w_E\in\bc$, and, $\o_E\in\G(\kf\otimes\cO_S(-E))$
such that,
\begin{itemize}
\item At every singularity as above, 
$$\tilde{\o}:=\rho^*\o-\sum_{E\geq E_{{\mathrm{red}}}} w_E\o_E\in \kf(x^p,y^q)$$
equivalently: no Taylor coefficients in the
obstructed sector.
\item The form of the Taylor expansion of
$\o_E$ is what we shall call {\it simple},
{\it i.e.} for $E$ locally $x^cy^d$, and
$a=p-c, b=q-d$,
$$
\o_E=\begin{cases}
d(x^{-a}y^{-b})\,\, \mod\, \kf(x^p,y^q)& \text{non-resonant},\\
x^{-a}y^{-b} \dfrac{dx}{x} \,\, \mod\, \kf(x^p,y^q)& \text{otherwise}
\end{cases}
$$
\end{itemize}
where amongst the singularities of the
form III.4.1 (b)/(c), non-resonant means:
not a rational eigenvalue $k/l$ with
$ka=lb$. In either case we can effect the
local calculation of 
the residue of $\o_E$ using the boundary
$\mx^a\my^b=\e$ by way of Stokes in the
non-resonant case, and by the approximately
invariant function commensurate to it otherwise.
A similar statement is not quite possible
for arbitrary $\o$, but we can achieve,
$$\rho^*\o \,=\, \o_0 + 
\sum_{E\geq E_{{\mathrm{red}}}} w_E\o_E
+\tilde{\o}\,\,\,\,\,\,\,\o_E\in\G(\kf\otimes\cO_S(-E)),\, w_E\in\bc 
$$
where all of the above are satisfied,
except possibly the simplicity 
of $\o_0$ (here $a=p$, $b=q$)
for singularities of type III.4.1 (c),
or, indeed the rational case of III.4.1 (b),
should, notation as per op. cit., one of
$i$, or $j$ be non-zero, but not both.

In order to apply these considerations,
let us introduce an easy version of
something already encountered in \S I.4, {\it viz:}

\noindent{\bf IV.4.6 Definition} Let $G$
be a graph. By a (directed) multiplier
is to be understood the assignment to
each possible direction, say, $+$, $-$,
of an edge multipliers $\nu_e^-, \nu_e^+
\in \br^{\times}_+$ such that $\nu_e^-\nu_e^+=1$.
The multiplier is said to be continuous
if the product around every directed
cycle is 1.

Alternatively, as in \S I.4, loops $\g$ in $G$
are just sequences of directed edges 
returning to whence they came,
so
we have a representation,
$$\nu:\pi_1 (G) \longra \br^{\times}_+:
\g\mpo \prod_{e\in\g} \nu_e^{\pm}$$
or equivalently a class $\nu\in\H^1(G,\br^{\times}_+)$.
Consequently, the continuity condition
is the triviality of this class, or what
amounts to the same thing assigning
weights $v\mpo a(v)$ at each vertex,
so that for a directed edge from 
$v$ to $w$ the multiplier is
$a(w)a(v)^{-1}$.

Now plainly, post \S IV.3, the only thing that we
have left to worry about 
in our residue calculation
is the
continuity of our multipliers.
In the decomposition of an arbitrary
$\o$ as above, the term $\tilde{\o}$
can be done in a myriad of ways.
For example just take as multipliers
the trivial co-cycle corresponding
to taking the actual multiplicities
$p(v)$ as the weights for the
vertices. In this way the local
strategy at singularities is
unchanged except in
the cases
III.4.1 (b)/(c), and a non-resonant
boundary, {\it i.e.} $p+\lb q\neq 0$,
for $\lb$ the eigenvalue of the
singularity normalised as per op.
cit. Should this case occur,
we switch to the strategy of
using Stokes. The dichotomy of
IV.3.4 ensures that this causes
no additional continuity problems
in the case of III.4.1 (b),
with the procedure to be employed
for gluing the boundary outside
the singularity to that close by
being as per IV.3.4 bis. Pretty
much the same is true in the case
III.4.1 (c), but one may have to
add discontinuities about extra
homology classes/loops around the
singularity in either component
even when local/global agreement occurs
in IV.3.4 to account for the
change in strategy. This is, however,
no worse than the total variation
of $w(\e)$ which we already bounded
appropriately in the proof of III.4.4.

Essentially identical remarks 
apply to computing the residues
of $\o_E$ for $E\geq E_{{\mathrm{red}}}$,
or for that matter $\o_0$ were
all the Taylor expansions at every
singularity of type III.4.1 (c) to
be simple, {\it i.e.} both $i$, $j$
non-zero or both zero. Indeed for
weights on the vertices, one takes
the multiplicities of $E_{{\mathrm{tot}}}
-E$, at the vertex which corresponds
to a component of the same, so, as 
for $\tilde{\o}$, the multiplicities
of $E_{{\mathrm{tot}}}$ in the case
of $\o_0$. The change in local
strategy is again, in the
cases of III.4.1 (b)/(c), only
at non-resonant boundaries, which
for $a,b$ the weights at the $x=0$,
$y=0$ vertices, means $a+b\lb\neq 0$,
which occasions no further problems
beyond those already discussed. We
also need to change strategy in the
cases III.3.1 (c), III.4.1 (a) when
$E\neq 0$, but this just means use
IV.4.4 (a), and nothing of substance
changes.

This leaves us to compute $\o_0$ in
general, or equivalently in a series
of im-probabilities of increasing
ludicrousness. In any case we view
this as the calculation of $\rho_*\o_0$,
and one improvement we can 
make a priori is to perform sufficient 
 blowing up so as to have 
only one edge between any pair of
vertices. As such $G$ is the dual
graph of the singular locus invariant
by the induced foliation in $\wh{Z}$.
Irrespectively of whether only one
rather than both of $i,j$ in III.4.1 (b)/(c)
is zero, such edges have only one
possible multiplier that permits
the computation to be done locally,
{\it i.e.} the eigenvalue $\lb_e$
of the directed edge in the notation
of \S IV.3,
and we'll refer to such edges as
{\it rigid}.
 If both $i$, $j$ are zero
then this is the multiplier that would
arise by taking the multiplicities
in $E_{{\mathrm{tot}}}$ as weights
for the vertices. Otherwise we can
make  further a priori improvement
by blowing up in a problematic 
singularity when one of these is
zero, say, $j$. Our initial situation
is therefore in the normalisation
post IV.3.5, a directed edge of the form,
$$
\begin{CD}
v_-:(x_e=0)@>{\lb_e=l/k}>e>v_+: (y_e=0)
\end{CD}
$$ 
where $\lb_e$ is the eigenvalue at the
singularity, and we'll
say that it is of type $(0,i_0)$,
$i=i_0$- so, slightly
confusingly, ``$i$'' follows the $x$-axis,
not $x=0$. Plainly it is rigid. The
effect of blowing this up is,
$$
\begin{CD}
v @>{\lb_e+1}>{e_1}> F_1 @>{\lb_e/(\lb_e +1)}>{e_0}>w
\end{CD}
$$
where $F_1$ is the exceptional divisor,
and the numbers over the edges $e_{\bullet}$
the new eigenvalues.
Now the above edge $e_0$ is again of type
$(0,i_0)$, and there is absolutely no improvement.
The edge $e_1$ is more interesting. More
precisely,
\begin{itemize}
\item[(a)] In the obvious notation,
it is never of type $(j,0)$ for any
$j\in\bn\cup\{0\}$.
\item[(b)] It could be of type $(0,i_1)$. This 
happens if $(k+l)| i_0$, in which case,
$$i_1\,=\, \dfrac{li_0}{(k+l)}\, < \,i_0$$
\item[(c)] Otherwise it's of type $(j_1, i_1)$,
$i_1j_1\neq 0$.
\end{itemize}
Evidently (b) can not repeat itself ad nauseum,
so blowing up in $e_1$, to get a new chain in
which we always blow up in the leftmost edge,
with exceptional divisors $F_{\bullet}$,
we eventually obtain for some $d\in\bn$,
$$
\begin{CD}
v @>{\lb_d}>{e_d}> F_d @>{\lb_{d-1}}>{e_{d-1}}> \hdots  @>{\lb_1}>{e_1}> F_1 @>{\lb_0}>{e_0}>w
\end{CD}
$$
where $\lb_{\bullet}$ is the eigenvalue of the 
singularity normalised as per \S IV.3 according
to the direction $e_{\bullet}$, and $d$ is the
first integer for which (c) holds at $e_d$. 
We perform this operation a priori in $U$
before doing anything else, so this is our
graph $G$, and when we
construct $\o_0$ it will be simple at $e_d$.
Furthermore by (a), we have the non-resonance
condition $p(F_d)-\lb_d p(v)$, for $p$
the multiplicity of $E_{{\mathrm{tot}}}$, so
by IV.2.2/4(b) we can use a local strategy
of Stokes type at $e_d$ for an arbitrary
multiplier $\nu_d$. The above has
no effect on the topology of the original
graph $G_0$, and if an edge $e:v\ra w$
is in both $G$ and $G_0$, then we take
$\nu_e=p(w)p(v)^{-1}$, otherwise we 
replace it by the sequence of edges
$e_{\bullet}$ as above, on which,
$$
\nu(e_m) =\begin{cases}
\dfrac{k\lb_d}{l}p(w)p(v)^{-1}& \text{if $m=d$},\\
\lb_m & \text{if $0\leq m< d$}
\end{cases}
$$
This yields a trivial co-cycle, 
so, for example to keep ourselves
consistent with \S IV.3, choose
a root $R$ in $G$, and define
a multiplier at a component $C$,
by way of,
$$\lb_C =\prod_{e\in\{R,C\}} \nu_e^{\pm}$$
then proceed exactly as in the proof
of IV.3.6 with the same minor
caveats for edges where a change in
local strategy takes place
(which, beyond those already 
encountered for $\tilde{\o}$,
are uniquely
those of the form $e_d$ as found above)
which were already noted for $\o_E$ and
$\tilde{\o}$. Consequently, we have proved
the essential in,

\noindent{\bf IV.4.7 Fact} Let things be as in
IV.4.1, but without the projectivity
assumption, 
with $\iy$ the ideal of the singular
locus, and $\dm$ an invariant measure with
support outside the exceptional divisor and
zero Segre class around the same, then:
$$\RES:\Ext^2_U (\iy,\kxf)\longra\Ext^3(\cO/\iy,\kxf)$$
is zero on $\dm$.

\noindent{\bf proof} It only remains to remove
the hypothesis of projectivity. As already observed
the moduli of $V$, {\it i.e.} weighted blow up
in sufficiently many points is formally convex.
Indeed the exceptional divisor $F$, in the notation
of the construction of $V$, satisfies the
co-homological criterion for ampleness in
its completion $\wh{V}$ in the singular locus.
Whence, the singular locus in $\wh{V}$ is
formally contractible, and by \cite{artin} this
contraction converges. Whence $F|_V$ is ample,
while the schematic covering $\tu\ra V$ can 
always be constructed for topological reasons-
$V\bsh F$ homotopic to the singular locus minus $B$,
which we may suppose to contain at least 3 points
in every component. $\Box$

\newpage

\end{document}